\documentclass[10pt]{amsart}

\usepackage[all]{xy}
\usepackage{latexsym}
\usepackage{amsmath}
\usepackage{amsthm}
\usepackage{amssymb}
\usepackage{a4wide}
\usepackage[utf8]{inputenc}
\usepackage{dsfont}
\usepackage{multirow}
\usepackage{slashbox}

\swapnumbers
\newtheorem{thm}[subsubsection]{Theorem}
\newtheorem{dei}[subsubsection]{Definition}
\newtheorem{prop}[subsubsection]{Proposition}
\newtheorem{lem}[subsubsection]{Lemma}
\newtheorem{cor}[subsubsection]{Corollary}

\theoremstyle{definition}
\newtheorem*{rem}{\sc Remark}
\newtheorem*{rems}{\sc Remarks}
\newtheorem*{pf}{\sc Proof}

\newtheorem*{ex}{\sc Examples}
\newtheorem*{ex1}{\sc Example}

\newcommand{\B}{\mathrm{B}}
\newcommand{\C}{\mathcal{C}}

\newcommand{\F}{\mathcal{F}}

\newcommand{\Lc}{\mathcal{L}}
\newcommand{\R}{\mathcal{R}}

\newcommand{\Q}{\mathcal{Q}}
\newcommand{\Po}{\mathcal{P}}
\newcommand{\Poa}{\Po^{\ash}}
\newcommand{\So}{\mathbb{S}}
\newcommand{\K}{\mathbb{K}}

\newcommand{\cqfd}{\ \hfill \square}

\newcommand{\Com}{\mathcal{C}om}
\newcommand{\BV}{\mathcal{BV}}

\newcommand{\ucFrob}{uc\mathcal{F}rob}

\newcommand{\ash}{\textrm{!`}}
\newcommand{\q}{\mathrm{q}}

\newcommand{\qPs}{\q\Po^{\textrm{!`}}}
\newcommand{\epi}{\twoheadrightarrow}
\newcommand{\mono}{\rightarrowtail}

\newcommand{\qiso}{\xrightarrow{\sim}}
\newcommand{\bp}{\, {}^{\backprime} \!}

\newcommand{\As}{\mathcal{A}s}
\newcommand{\uAs}{u\mathcal{A}s}
\newcommand{\quAs}{\q \uAs^{\ash}}
\newcommand{\uAsinf}{uA_{\infty}}
\newcommand{\suAinf}{\texttt{su}A_{\infty}}
\newcommand{\muns}{\mu_n^S}
\newcommand{\munsc}{\overline{\mu}_n^S}

\newcommand{\id}{\operatorname{id}}
\newcommand{\End}{\operatorname{End}}

\newcommand{\thetaop}{\operatorname{\theta}}
\newcommand{\inc}{\operatorname{in}}
\newcommand{\Tw}{\operatorname{Tw}}
\newcommand{\Hom}{\operatorname{Hom}}

\newcommand{\un}{\vcenter{
\xymatrix@M=0pt@R=3pt@C=3pt{&&\\
& \ar@{{*}}[u] \ar@{-}[d]\\
& &}}}

\newcommand{\ass}{\vcenter{ \xymatrix@M=0pt@R=4pt@C=4pt{
      & & & & \\ \ar@{-}[ddrr] & & \ar@{-}[dl] & & \ar@{-}[ddll]\\
      &  & &  & \\
      & & \ar@{-}[d] & &\\
      & & }}
- 
\vcenter{
    \xymatrix@M=0pt@R=4pt@C=4pt{& & & & \\
      \ar@{-}[drdr] & &\ar@{-}[dr] & & \ar@{-}[dldl]  \\
      & & & & \\
      & &\ar@{-}[d] & & \\
      & & }}}

\newcommand{\mun}{\vcenter{ \xymatrix@M=0pt@R=4pt@C=4pt{
      & & & & \\
      & & & &\\
      & \ar@{{*}}[u] \ar@{-}[dr] & & \ar@{-}[dl] & \\
      & & \ar@{-}[d] & &\\
      & & }}}

\newcommand{\munsplit}{\vcenter{ \xymatrix@M=0pt@R=4pt@C=4pt{
      & & & & \\
      &\ar@{{*}}[u] & & &\\
      & \ar@{-}[dr] & & \ar@{-}[dl] & \\
      & & \ar@{-}[d] & &\\
      & & }}}

\newcommand{\lmun}{\vcenter{\xymatrix@M=0pt@R=7pt@C=10pt{
      & & & & \\
      & & & &\\
      & \ar@{{*}}[u] \ar@{-}[dr] & & \ar@{-}[dl] & \\
      & & \ar@{-}[d] & &\\
      & & }}}

\newcommand{\unm}{\vcenter{
    \xymatrix@M=0pt@R=4pt@C=4pt{& & & & \\
      & & & &\\
      & \ar@{-}[dr] & & \ar@{{*}}[u] \ar@{-}[dl] & \\
      & &\ar@{-}[d] & & \\
      & & }}}

\newcommand{\unmsplit}{\vcenter{
    \xymatrix@M=0pt@R=4pt@C=4pt{& & & & \\
      & & & \ar@{{*}}[u] &\\
      & \ar@{-}[dr] & & \ar@{-}[dl] & \\
      & &\ar@{-}[d] & & \\
      & & }}}

\newcommand{\mtwo}{\vcenter{\xymatrix@M=0pt@R=7pt@C=10pt{
    &\ar@{-}[dr] & & \ar@{-}[dl] \\
    & &\ar@{-}[d] & \\
    & & & }}}

\newcommand{\smtwo}{\vcenter{\xymatrix@M=0pt@R=4pt@C=4pt{
    &\ar@{-}[dr] & & \ar@{-}[dl] \\
    & &\ar@{-}[d] & \\
    & & & }}}

\newcommand{\mthree}{\vcenter{
    \xymatrix@M=0pt@R=4pt@C=4pt{
      \ar@{-}[dr] & \ar@{-}[d]& \ar@{-}[dl] &\\
      &\ar@{-}[dd] &  &\\
      & & & \\
      & & & }}}
\newcommand{\mthreec}{\vcenter{
    \xymatrix{
       *{}\ar@{-}[rd] &*{} \ar@{-}[d]& *{}\ar@{-}[ld] & \\
       &*{}\ar@{-}[d] & &  \\
       & & &
    }}}
\newcommand{\mfivectwo}{\vcenter{ \xymatrix@M=0pt@R=7pt@C=10pt{
    & & & \\
    & & & & \\
    \ar@{{*}}[u] \ar@{-}[drr] & \ar@{-}[dr] & \ar@{-}[d] &
    \ar@{-}[dl] \ar@{{*}}[u]  & \ar@{-}[dll] \\
    & & \ar@{-}[d] \\
    & &
    }}}

\newcommand{\smfivectwo}{\vcenter{\xymatrix@M=0pt@R=4pt@C=4pt{
    & & & \\
    & & & & \\
    \ar@{{*}}[u] \ar@{-}[drr] & \ar@{-}[dr] & \ar@{-}[d] &
    \ar@{-}[dl] \ar@{{*}}[u]  & \ar@{-}[dll] \\
    & & \ar@{-}[d] \\
    & &
    }}}

\newcommand{\ident}{\vcenter{\xymatrix@M=0pt@R=7pt@C=10pt{
  \ar@{-}[dd]& \\
  & \\
  &
}}}

\newcommand{\draftnote}[1]{}

\author{Joseph Hirsh and Joan Mill\`es}
\thanks{The first author was supported by a
  National Science Foundation Graduate Research Fellowship.}
\thanks{The second is supported by the ANR grant JCJC06 OBTH}

\title{Curved Koszul duality theory}

\begin{document}

\maketitle

\begin{abstract}
We extend the bar-cobar adjunction to operads and properads, not necessarily augmented. Due to the default of augmentation, the objects of the dual category are endowed with a curvature. We handle the lack of augmentation by extending the category of coproperads to include objects endowed with a curvature. As usual, the bar-cobar construction gives a (large) cofibrant resolution for any properad, such as the properad encoding unital and counital Frobenius algebras, a notion which appears in 2d-TQFT. We also define a curved Koszul duality theory for operads or properads presented with quadratic, linear and constant relations, which provides the possibility for smaller relations. We apply this new theory to study the homotopy theory and the cohomology theory of unital associative algebras.
\end{abstract}

\section*{Introduction}

In \cite{Hochschild}, Hochschild introduced a (co)homology theory for associative algebras and in \cite{Stasheff}, Stasheff introduced the homotopy theory for associative algebras. Nowadays, we know how to describe these theories in operadic terms, but this approach does not encode the units in unital associative algebras. In order to define a homotopy theory and a cohomology theory for unital associative algebras, we refine the operadic theory and more precisely its Koszul duality theory.\\

In representation theory, we interpret algebras as operations with one input and one output, $\mathrm{Hom}(A,\, \mathrm{End}(V))$. To encode operations with several inputs and one output, oneuses the notion of an operad \cite{May, BoardmanVogt}. More generally, one uses the notion of properads to encode operations with several inputs and several outputs \cite{Vallette2}. An associative algebra is a special kind of operad and an operad is a special kind of properad, and theories about properads generalize those of operads and associative algebras. For example a bar-cobar adjunction defined in a properadic setting generalizes one defined for operads and algebras. In \cite{Vallette2}, the bar construction (denoted by $\B$) assigned a coaugmented dg coproperad to an augmented dg properad and the cobar construction (denoted $\Omega$) assigns an augmented dg properad to a coaugmented dg coproperad, and the two constructions are adjoint. An important property of the adjunction is that the bar-cobar composition $\Omega \B \Po$ defines a cofibrant resolution of an augmented dg properad $\Po$.

In this paper, we extend the bar-cobar adjunction $(\Omega, \B)$ to non-augmented properads. We generalize the notion of dg coproperad to involve curvature ($d^2 \neq 0$, but its deviation from 0 is controlled by a term we call ``curvature"). Our bar construction assigns a coaugmented curved coproperad to a (not necessarily augmented) dg properad. We then extend the cobar construction of coaugmented coproperads to include coaugmented curved coproperads, resulting in a (not necessarily augmented) dg properad (with no curvature). The composition bar-cobar provides a cofibrant resolution $\Omega \B \Po$ of a properad $P$. For example, we obtain a cofibrant resolution for the properad encoding unital and/or counital Frobenius algebras. Since the datum of a $2$-dimensional topological quantum field theory, 2d-TQFT for short, is equivalent to a unital and counital Frobenius algebra structure \cite{Abrams, Kock}, this provides homotopy tools to study 2d-TQFT. With our model, the methods of \cite{Wilson} apply to show that the differential forms $\Omega (M)$ on a closed, oriented manifold $M$ bear a unital and counital Frobenius algebra structure up to homotopy.

The bar-cobar resolution $\Omega \B \Po$ is large and it is often desirable to have a smaller resolution. To this end, we develop a curved Koszul duality theory for properads generalizing the Koszul duality theory for properads \cite{Vallette2}, operads \cite{GetzlerJones, GinzburgKapranov}, and associative algebras \cite{Priddy}. One of the main object is the Koszul dual coproperad $\Poa$, which has, here, a curvature. It applies to properads with a quadratic, linear and constant presentation. The properads for which this theory apply are called Koszul properads. In this case, the cobar construction $\Omega \Poa$ is a resolution of $\Po$. We summarize the different generalizations of the Koszul duality theory in the following table:
$$\begin{tabular}{|l||c|c|p{2.7cm}|}
\hline
\multirow{2}{3.8cm}{\backslashbox{\bf Monoids}{\bf Relations}} & \multirow{2}{*}{Homogeneous quadratic} & \multirow{2}{*}{Quadratic and linear} & \multirow{2}{2.8cm}{Quadratic, linear and constant}\\
 &  &  & \\
\hline \hline
Associative algebras & \multicolumn{2}{|c|}{\cite{Priddy}} & \multicolumn{1}{c|}{\cite{Positselski, PolishchukPositselski}}\\
\hline
Operads & \cite{GetzlerJones, GinzburgKapranov} & \multirow{2}{*}{\cite{Ga-CaToVa}} & \multirow{2}{2.6cm}{{\it Section \ref{curvedKD}} of this paper}\\
\cline{1-2}
Properads & \cite{Vallette2} &  & \\
\hline
\end{tabular}$$\\

The operad $\uAs$ encoding unital associative algebras is an example of an operad with quadratic, linear and constant relations. It is an inhomogeneous Koszul operad in the previous sense. Hence we get a ``small'' cofibrant resolution $\uAsinf := \Omega \uAs^{\ash} \qiso \uAs$. This particularly simple resolution of the operad $\uAs$ allows us to define the notion of \emph{homotopy unital associative algebras}. This notion corresponds to the notion of homotopy unit for $A_{\infty}$-algebra which appears in \cite{Fukaya}. However, our presentation in terms of algebras over a cofibrant operad implies good homotopy properties for these algebras. With this approach, we also obtain functorial resolutions on the level of unital associative algebras. We use these other resolutions to study the cohomology theory of unital associative algebras.\\

We begin the paper with a survey of the results on homotopy unital associative algebras expressed in an internal language, explained without, for example, the words ``operad" or ``properad." This section corresponds to the results obtained in the last section of this paper. In Section 2, we recall definitions of associative algebras, operads and properads. In Section 3, we extend the bar and the cobar construction to the non-augmented framework and we define the notion of curved twisting morphims. In Section 4, we extend the Koszul duality theory for homogeneous quadratic properads to properads with quadratic, linear and constant relations. Section 5 is devoted to resolution of non-augmented properads as bimodules over themselves and to functorial resolutions of $\Po$-algebras. Section 6 studies the operad encoding unital associative algebras. We describe the homotopy theory and the cohomology theory for this category of algebras.\\

In this paper, we work over a field $\K$ of characteristic $0$.

\tableofcontents

\section{Results on unital associative algebras}
In this section, we develop the homotopy and cohomology theories of unital associative algebras. The definitions, proofs, techniques, and pictorial descriptions of the results are based on operad theory and can be found in Section \ref{resandtransfer}. However, this section does not contain the word ``operad'' and can be read independently from the rest of the paper. The comparison with the work of \cite{Fukaya} is referred to Section \ref{resandtransfer}.

\subsection{Unital associative algebra}
  A \emph{unital associative differential graded algebra} is a quadriple
  $(A,\, \mu,\, u,\, d_{A})$, where $(A,\, d_{A})$
  is a dg module, $\mu: A \otimes A \to A$, and
  $u: \K \to A$ are dg modules maps, such that the map $\mu$ is associative and such that the element
  $u(1_{\K})$ is a left and right unit for the associative product $\mu$.\\

The version of this structure ``up to homotopies''  is what we call a
\emph{$\uAsinf$-algebra}, for \emph{homotopy unital associative algebra}. Let $f : V \rightarrow W$ be a homogeneous $\K$-linear map of degree $|f|$. We denote its derivative by $\partial(f) := d_{W} \circ f - (-1)^{|f|}f \circ d_{V}$.

\subsection{Homotopy unital associative algebra} \label{uAsresults}
A \emph{homotopy unital associative algebra} or \emph{$\uAsinf$-algebra structure on a dg module $(A,\, d_{A})$} is given by a collection of maps $\{ \muns \}$, where the set $S$ runs over the set of
subsets of $\{ 1,\, \ldots,\, n\}$ for any integer $n \geq 2$ and where $S = \{ 1\}$ when $n = 1$. The $\muns$ are given pictorially by planar corollas with $n$ entries labelled by $1,\, \cdots ,\, n$ on which we put ``corks'' when the label is in $S$. For example, we have $\mu_{3}^{\{1\}} = \vcenter{\xymatrix@M=0pt@R=4pt@C=5pt{
        &&&\\
        & \ar@{-}[dr] \ar@{{*}}[u] & \ar@{-}[d]& \ar@{-}[dl]\\
        && \ar@{-}[d] & \\
        &&}}$. The maps $\muns: A^{\otimes (n-|S|)} \to A$ are of degree $n-2+|S|$ and satisfy the
  following identities:
$${\small \left\{ \begin{array}{lcl}
\partial \left(\hspace{.1cm} \vcenter{\xymatrix@M=0pt@R=5pt@C=5pt{
     &&\\
     \ar@{{*}}[u]&&\\
    & \ar@{-}[d] \ar@{-}[lu] \ar@{-}[ur] &\\
    &&}} \right) &=& \vcenter{\xymatrix@M=0pt@R=5pt@C=5pt{
        &&\\
        \ar@{{*}}[u] &&\\
        \ar@{-}[dr] && \\
        & \ar@{-}[d] \ar@{-}[ur] \\
        && }} - |\\
\partial \left(\vcenter{\xymatrix@M=0pt@R=5pt@C=5pt{
     &&\\
    &&\ar@{{*}}[u]\\
    & \ar@{-}[d] \ar@{-}[lu] \ar@{-}[ur] &\\
    &&}} \hspace{.1cm} \right) &=& \vcenter{\xymatrix@M=0pt@R=5pt@C=5pt{
    &&\\
        && \ar@{{*}}[u]\\
        \ar@{-}[dr] && \\
        & \ar@{-}[d] \ar@{-}[ur] \\
        && }} - |,   \end{array} \right. }$$
where the empty space between the corollas and the corks is the composition of applications and where $|$ is the identity of $A$ and, for $(n,\, S) \neq (2,\, \{ 1\})$ and  $(n,\, S) \neq (2,\, \{ 2\})$,
  $$\partial(\muns ) = \sum_{\substack{p+q+r=n\\
          p+1+r=m}}
      (-1)^{q(r+|S_{1}|)+|S_{2}||S_{1}'| + p+1} \mu_{m}^{S_{1}} \circ
      (\underbrace{\id, \ldots, \id}_{p-|S'_{1}|},\, \mu_{q}^{S_{2}},\, \underbrace{\id, \ldots,
          \id}_{r-|S''_{1}|}),$$
or pictorially,
$${\small \partial \left(\vcenter{\xymatrix@M=0pt@R=7pt@C=10pt{
    & & & & \\
    & & & & & \\
     \ar@{-}[drrr] & \ar@{{*}}[u] \ar@{-}[drr] & \ar@{-}[dr] & \ar@{-}[d] & \ar@{-}[dl] \ar@{{*}}[u]  \
& \ar@{-}[dll] \\
    & & & \ar@{-}[d] & \\
    & & &
    }}\right) = \sum_{\substack{p+q+r=n\\
          p+1+r=m}}
      (-1)^{q(r+|S_{1}|)+|S_{2}||S_{1}'| + p+1} \vcenter{\xymatrix@M=0pt@R=4pt@C=8pt{
        && \ar@{-}[dd] & \ar@{-}[ddl] \\
        & \ar@{-}[dr] \ar@{{*}}[u] &&\\
        && \ar@{-}[d] & \\
        &&&&&&\\
        \ar@{-}[ddrr] && \ar@{-}[dd] &&& \ar@{-}[ddlll]\\
        &&& \ar@{-}[dl] \ar@{{*}}[u] &&&\\
        && \ar@{-}[d] && .& \\
         & & }}}$$

\begin{ex}
  \begin{enumerate}
  \item[]
  \item Every unital associative differential graded algebra, dga for short, $(A,\, \mu,\, u)$ is
    naturally a $\uAsinf$-algebra by
    \begin{equation*}
      \mu_{n}^{S} = \left\{
        \begin{array}{rl}
          \mu & \text{if } n=2 \text{ and } S=\emptyset \\
          u & \text{if } n=1 \text{ and } S = \{ 1 \} \\
          0 & \text{otherwise} \\
        \end{array}
      \right.
    \end{equation*}
  \item A \emph{strictly unital $A_{\infty}$-algebra $(A,
      \{\mu_n\}_{n\geq 1},\, u)$} is an $A_{\infty}$-algebra $(A,
    \{\mu_n\}_{n\geq 1})$ \label{strictAinf} with $u \in A$ so that $d_A (u)=0$ and $u$ is a left and
    right unit for $\mu_2$, and $u$ annihilates
    $\mu_n$ for $n \geq 3$ \cite{Kontsevich}. Every strictly unital $A_{\infty}$-algebra is naturally a
    $\uAsinf$-algebra by
    \begin{equation*}
      \mu_{n}^{S} = \left\{
        \begin{array}{ll}
          \mu_n & \text{if } S = \emptyset \\
          u & \text{if } n = 1 \text{ and } S = \{1\} \\
          0 & \text{otherwise}
        \end{array}
      \right.
    \end{equation*}
  \end{enumerate}
\end{ex}
Every $\uAsinf$-algebra contains an $A_{\infty}$-algebra if we take $\mu_n := \mu_{n}^{\emptyset}$ for all $n \geq 1$. The algebraic structure given by a $\uAsinf$-algebra provides homotopies for the ``unital'' relations, along with the ``associativity'' homotopies.

\subsection{Infinity-morphism}

We define the notion of infinity-morphism between two $\uAsinf$-algebras $A$ and $B$ by a collection of maps $f_{n}^{S}: A^{\otimes (n - |S|)} \to B$ of degree $n-1+|S|$, represented graphically by planar trees with ``corks'' as the $\uAs_{\infty}$-algebra structures but with a triangle $\bigtriangledown$ as vertex. For example, we have $f_{3}^{\{1\}} = \vcenter{\xymatrix@M=0pt@R=4pt@C=5pt{
        &&&\\
        & \ar@{-}[dr] \ar@{{*}}[u] & \ar@{-}[d]& \ar@{-}[dl]\\
        && \bigtriangledown \ar@{-}[d] & \\
        &&}}$. The $f_{n}^{S}$ satisfy the relations: 
$${\small \partial \left(\vcenter{\xymatrix@M=0pt@R=8pt@C=10pt{
    & & & & \\
    & & & & & \\
    \ar@{-}[drrr] & \ar@{{*}}[u] \ar@{-}[drr] & \ar@{-}[dr] & \ar@{-}[d] & \ar@{-}[dl] \ar@{{*}}[u]  \
& \ar@{-}[dll] \\
    & & & \bigtriangledown \ar@{-}[d] & \\
    & & &
    }}\right) = \sum \pm\ \vcenter{\xymatrix@M=0pt@R=4pt@C=11pt{
        && \ar@{-}[dd] & \ar@{-}[ddl] \\
        & \ar@{-}[dr] \ar@{{*}}[u] &&\\
        && \ar@{-}[d] & \\
        \ar@{-}[ddrr] && \ar@{-}[dd] &&& \ar@{-}[ddlll]\\
        &&& \ar@{-}[dl] \ar@{{*}}[u] &&&\\
        && \bigtriangledown \ar@{-}[d] &&& \\
         & & }} - \sum \pm \vcenter{\xymatrix@M=0pt@R=4pt@C=11pt{
        & & & & \\
        \ar@{-}[ddr] &&& \ar@{-}[ddr] && \ar@{-}[ddl] && \ar@{-}[dd] \\
        && \ar@{-}[dl] \ar@{{*}}[u] &&&&&\\
       & \bigtriangledown \ar@{-}[dddrrr] &&& \bigtriangledown \ar@{-}[dddd] &&& \bigtriangledown \ar@{-}[d]\\
        &&&&&&& \ar@{-}[ddlll]\\
        &&&&& \ar@{-}[dl] \ar@{{*}}[u] &&\\
        &&&&&&&\\
        &&&&&&&,}}}$$
where the planar trees with ``corks'' and no triangle represent the $\uAsinf$-algebra structure of $A$ on the top and the $\uAsinf$-algebra structure of $B$ on the bottom. With this definition of infinity-morphism, we prove a rectification theorem.

\begin{thm}[Rectification Theorem, Theorem \ref{universalrectification}]
Let $A$ be a unital associative algebra up to homotopies. We can rectify $A$: there is a unital associative algebra $A'$ such that $A$ is infinity-quasi-isomorphic to $A'$.
\end{thm}

Moreover, we have a tranfer theorem.

\begin{thm}[Homotopy Transfer Theorem, Theorem \ref{transferthm}]
Let $A$ be a homotopy unital associative algebra and let $V$ be a chain complex. Given a strong deformation retract%
$$\xymatrix{     *{ V \ \ } \ar@<.8ex>[r]^(.45){i} & *{\
A \quad} \ar@(dr,ur)[]_h \ar@<.8ex>[l]^(.43){p}},$$
i.e., $p$ and $i$ are chain maps, where $p \circ i = \id_V$ and $d_Ah+hd_A =
\id_A - i \circ p$, there is a natural $\uAsinf$-algebra
structure on $V$, and a natural extension of $i$ to an infinity-morphism.
\end{thm}

\subsection{Comparison with the literature}

In the literature, there have been definitions of ``weakly
unital'' (or ``homotopy unital'') $A_{\infty}$-algebras
\cite{Kontsevich, Lyubashenko2, Fukaya2}. Each of these definitions,
however, is a definition of a \emph{property} of
$A_{\infty}$-algebra. In \cite{Lyubashenko} these properties are
shown to be equivalent. Our notion of ``homotopy unital'' $A_{\infty}$-structure, what we are calling a $\uAsinf$-algebra, is an additional \emph{structure} on an $A_{\infty}$-algebra. \\

One has the following statement:
\begin{thm}[\cite{Kontsevich}]
Let $A$ be an $A_{\infty}$-algebra with a homotopy unit. Then there
exists a strictly unital $A_{\infty}$-algebra which is
quasi-isomorphic to $A$. 
\end{thm}
We prove an analogous theorem for $\uAsinf$-algebras.
\begin{thm}[Theorem \ref{qisostrict}]
Let $A$ be a $\uAsinf$-algebra. There is a strictly unital
$A_{\infty}$-algebra structure on the homology of $A$ which is infinity-quasi-isomorphic to $A$.
\end{thm}
We extend this theorem to a broad class of algebraic structures, including
Batalin-Vilkovisky algebras and commutative algebras.

\subsection{André-Quillen cohomology theory for unital associative algebra}

Following the ideas of Quillen, we define a cohomology theory associated to any unital associative dga $A$ with coefficients in a $A$-bimodule $M$, denoted $\mathrm{H}_{u\As}^{\bullet}(A,\, M)$. We prove that this cohomology theory is an Ext-functor and that it is equal to the Hochschild cohomology theory of the associative algebra $A$.
\begin{thm}[Theorem \ref{hochschild}]
Let $A$ be a unital associative dga. We have
$$\mathrm{H}_{\uAs}^{\bullet} (A,\, M) \cong \mathrm{HH}^{\bullet +1} (A,\, M).$$
\end{thm}

\section{Properad}\label{properad}

In this section, we recall the notion of algebra, operad and properad
as successive generalizations. We refer to the book of Loday and
Vallette \cite{LodayVallette} for a complete and modern exposition
about algebras and operads in \textsf{dg mod}, to the book of \cite{MarklShniderStasheff} for another presentation and to the thesis of Vallette \cite{Vallette2} for properads.

\subsection{Algebra}

Let \textsf{$\K$-mod} denote the monoidal category $(\K \textrm{-mod},\, \otimes_{\K},\, \K)$ of $\K$-modules. A \emph{unital associative algebra} is a monoid $(A,\, \mu,\, u)$ in this monoidal category. The product $\mu : A\otimes_{\K} A \rightarrow A$ is associative and $u : \K \rightarrow A$ is a \emph{unit} for the product.

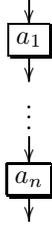
\begin{figure}[h]
$$\xymatrix@R=10pt@C=10pt@M=3pt{\ar[d]\\
*+[F-,]{a_1} \ar[d]\\
\vdots \ar[d]\\
*+[F-,]{a_n} \ar[d]\\
&}$$
\caption{\label{fig1} Representation of the product $a_{1} \cdots a_{n}$}
\end{figure}

As in representation theory, the elements of $A$ are seen as operations with one input and one output. Then we represent the product $a_{1} \cdots a_{n}$ by a vertical bar whose vertices are indexed by the $a_{i}$, see Figure \ref{fig1}.

\subsection{Operad}

An \emph{$\So$-module} $\Po = \{ \Po(n)\}_{n \geq 0}$ is a collection of vector spaces $\Po(n)$ endowed with right action of the symmetric group $\So_{n}$. One defines from \cite{May} the monoidal product $\circ$ on the category of $\So$-modules by
$$(\Po \circ \Q)(n) := \bigoplus_{k \geq 0} \left(\Po(k) \otimes \left(\bigoplus_{i_1 + \cdots + i_k = n} (\Q(i_1) \otimes \cdots \otimes \Q(i_k)) \otimes_{\So_{i_{1}}\times \cdots \times \So_{i_{k}}} \K[\So_{n}] \right) \right)_{\So_{k}},$$
where the notation $\otimes_{\So_{k}}$ stands for the space of coinvariants under the (diagonal) action of the symmetric group $\So_{k}$:
$$(p \otimes q_{1} \otimes \cdots \otimes q_{k} \otimes \sigma) \cdot \nu := p \cdot \nu \otimes q_{\nu^{-1}(1)} \otimes \cdots \otimes q_{\nu^{-1}(k)} \otimes \bar{\nu}^{-1} \cdot \sigma$$
for any $p \in \Po(k)$, $q_{j} \in \Q(i_{j})$, $\sigma \in \So_{n}$ and $\nu \in \So_{k}$ with $\bar{\nu} \in \So_{n}$ the induced block-wise permutation. This monoidal product encodes the composition of multilinear operations and we represent it by 2-levels trees as shown in Figure \ref{fig2}.

\begin{figure}[h]
$$\xymatrix@R=10pt@C=10pt{ 5 \ar[dr] & 3 \ar[d] & 4 \ar[dl]& 1 \ar[dr]& & 6 \ar[dl]
& 8 \ar[dr] & 2 \ar[d] & 7 \ar[dl] \\
\ar@{.}[r] & *+[F-,]{q_1}\ar@{.}[rrr] \ar[drrr]& & &
*+[F-,]{q_2}\ar@{.}[rrr] \ar[d]
 & & & *+[F-,]{q_3} \ar@{.}[r]
\ar[dlll]& \\
\ar@{.}[rrrr] & & & & *+[F-,]{p} \ar[d]\ar@{.}[rrrr] & & & & \\
 & & & & & & & &  }$$
\caption{\label{fig2} An element in $(\Po \circ \Q) (8)$}
\end{figure}

The unit for the monoidal product is $I := (0,\, \K,\, 0,\, \ldots)$ where the $\K$ is in arity $1$ and represent the identity element modeled by the tree $|$. It forms a monoidal category denoted by \textsf{$\So$-Mod}.

An \emph{operad} is a monoid $(\Po,\, \gamma,\, u)$ in the monoidal category of $\So$-modules $\So \textsf{-Mod}$. The associative product $\Po \circ \Po \rightarrow \Po$ is called the \emph{composition product} and $u : I \rightarrow \Po$ is the \emph{unit} for the composition product.

\begin{ex1}
A unital associative algebra induces an operad by this injective map
$$\textsf{Unital associative algebras} \mono \textsf{Operads},\ \ \ A \mapsto (0,\, A,\, 0,\, \ldots)$$
\end{ex1}

\subsection{Properad}

Algebras encode operations with one input and one output. Operads encode operations with several inputs and one output. To encode operations with multiple inputs and outputs, one uses the notion of \emph{properad}.

An \emph{$\So$-bimodule} $\Po$ is a collection $ \{ \Po(m,\, n) \}_{m,n\geq 0}$ of $\So_{m}^{\textrm{op}}$-$\So_{n}$-bimodules. One recalls from \cite{Vallette2} a monoidal product using 2-levels graphs as in Figure \ref{fig3}.

\begin{figure}[h]
$$\xymatrix@R=15pt@C=15pt{& \ar[dr] &  \ar[d] &  \ar[dl]& \ar[dr]& &\ar[dl] &\\
& \ar@{.}[r] & *+[F-,]{q_1}\ar@{.}[rrr] \ar[drrr] \ar@/_1pc/[d] & & & *+[F-,]{q_2}\ar@{.}[r] \ar@/^1pc/[d] \ar[dlll] |!{[d];[lll]}\hole & &\\
& \ar@{.}[r] & *+[F-,]{p_1} \ar@{.}[rrr] \ar[dl] \ar[dr] & & & *+[F-,]{p_{2}} \ar[d]\ar@{.}[r] & & \\
 & & & & & & & &}$$
\caption{\label{fig3} An element in $(\Po \boxtimes \Q) (3,\, 5)$}
\end{figure}

Let $a$ and $b$ the number of vertices on the first level and on the second level respectively. Let $N$ be the number of internal edges between the two levels. We associate to an $a$-tuple of integers $\bar{\imath} = (i_{1},\, \ldots ,\, i_{a})$ the sum $|\bar{\imath}| := i_{1} + \cdots + i_{a}$. To any pair of $a$-tuples $\bar{\imath}$ and $\bar{\jmath}$ we denote by $\Po(\bar{\jmath},\, \bar{\imath})$ the tensor product $\Po(j_{1},\, i_{1}) \otimes \cdots \otimes \Po(j_{a},\, i_{a})$ and by $\So_{\bar{\imath}}$ the image of $\So_{i_{1}} \times \cdots \times \So_{i_{a}}$ in $\So_{|\bar{\imath}|}$.

Let $\bar{k} = (k_{1},\, \ldots ,\, k_{b})$ be a $b$-tuple and let $\bar{\jmath} = (j_{1},\, \ldots ,\, j_{a})$ be an $a$-tuple such that $|\bar{k}| = |\bar{\jmath}| = N$. A \emph{$(\bar{k},\, \bar{\jmath})$-connected permutation} is a permutation $\sigma$ in $\So_{N}$ such that the graph of a geometric representation of $\sigma$ is connected when one connects the inputs labelled by $j_{1} + \cdots + j_{i} + 1,\, \ldots ,\, j_{1} + \cdots + j_{i + 1}$ for $0 \leq i \leq a-1$ and the outputs labelled by $k_{1} + \cdots + k_{i} + 1,\, \ldots ,\, k_{1} + \cdots + k_{i + 1}$ for $0 \leq i \leq b-1$. We denote by $\So^{c}_{\bar{k}, \bar{\jmath}}$ the set of $(\bar{k},\, \bar{\jmath})$-connected permutations.

We define the monoidal product $\boxtimes$ on the category of $\So$-bimodules by
$$(\Po \boxtimes \Q) (m,\, n) := \bigoplus_{N \in \mathbb{N}} \left( \bigoplus_{\bar{l}, \bar{k}, \bar{\jmath}, \bar{\imath}} \K[\So_{m}] \otimes_{\So_{\bar{l}}} \Po(\bar{l}, \bar{k}) \otimes_{\So_{\bar{k}}} \K[\So^{c}_{\bar{k}, \bar{\jmath}}] \otimes_{\So_{\bar{\jmath}}} \Q(\bar{\jmath},\, \bar{\imath}) \otimes_{\So_{\bar{\imath}}} \K[\So_{n}] \right)_{\So_{b}^{\textrm{op}} \times \So_{a}},$$
where the second direct sum runs over the $b$-tuples $\bar{l}$, $\bar{k}$ and the $a$-tuples $\bar{\jmath}$, $\bar{\imath}$ such that $|\bar{l}| = m$, $|\bar{k}| = |\bar{\jmath}| = N$, $|\bar{\imath}| = n$ and we consider the module of coinvariants with respect to the $\So_{b}^{\textrm{op}} \times \So_{a}$-action:
$$\rho \otimes p_{1} \otimes \cdots \otimes p_{b} \otimes \sigma \otimes q_{1} \otimes \cdots \otimes q_{a} \otimes \omega \sim \rho \cdot \tau_{\bar{l}}^{-1} \otimes p_{\tau(1)} \otimes \cdots \otimes p_{\tau(b)} \otimes \tau_{\bar{k}} \cdot \sigma \cdot \nu_{\bar{\jmath}} \otimes q_{\nu^{-1}(1)} \otimes \cdots \otimes q_{\nu^{-1}(a)} \otimes \nu_{\bar{\imath}}^{-1} \cdot \omega,$$
for $\rho \in \So_{m}$, $\omega \in \So_{n}$, $\sigma \in \So^{c}_{\bar{k}, \bar{\jmath}}$ and for $\tau \in \So_{b}$ with $\tau_{\bar{k}}$ the associated block-wise permutation, $\nu \in \So_{a}$ with $\tau_{\bar{\jmath}}$ the associated block-wise permutation. We write an element in $\Po \boxtimes \Q$ like this $\theta (p_{1},\, \ldots ,\, p_{b}) \sigma (q_{1},\, \ldots ,\, q_{a}) \omega$. The unit $I$ for the monoidal product is given by
$$\left\{ \begin{array}{lcll}
I (1,\, 1) & := & \K & \textrm{and}\\
I (m,\, n) & := & 0 & \textrm{otherwise.}
\end{array} \right.$$
It forms a monoidal category denoted by $\So$\textsf{-biMod}.

A \emph{properad} is a monoid $(\Po,\, \gamma,\, u)$ in the monoidal category $\So \textsf{-biMod}$ of $\So$-bimodules. The associative product $\gamma : \Po \boxtimes \Po \rightarrow \Po$ is called the \emph{composition product} and $u : I \rightarrow \Po$ is the \emph{unit} for the composition product.

\begin{ex1}
An operad induces a properad as follows
$$\textsf{Operads} \mono \textsf{Properads},\ \ \ \Po \mapsto \widetilde{\Po} \textrm{, where } \left\{ \begin{array}{lcll}
\widetilde{\Po} (1,\, n) & := & \Po(n) & \textrm{and}\\
\widetilde{\Po} (m,\, n) & := & 0 & \textrm{for $m \neq 1$.}
\end{array} \right.$$
\end{ex1}

Finally, we have the following inclusions:
$$\begin{array}{lccccc}
\textrm{\bf Monoidal category:} & (\K \textrm{-Mod},\, \otimes_{\K}) & \mono & (\So \textrm{-Mod},\, \circ) & \mono & (\So \textrm{-biMod},\, \boxtimes)\\
\textrm{\bf Monoid:} & \textrm{Associative algebras} & \mono & \textrm{Operads} & \mono & \textrm{Properads.}
\end{array}$$
So the results given in this paper for properads apply to operads and algebras as well.

One defines dually the notions of coalgebra, cooperad, coproperad. For example, a \emph{coproperad} is a comonoid $(\C,\, \Delta,\, \eta)$ in the monoidal category of $\So$-bimodules $\So$\textsf{-biMod}. The \emph{coproduct} $\Delta : \C \rightarrow \C \boxtimes \C$ is coassociative and admits the counit $\eta : \C \rightarrow I$. All these definitions extend to the \emph{differential graded} setting, or \emph{dg} setting for short. The differentials are compatible with the properad structure, resp. coproperad structure, in the sense that they are \emph{derivations}, resp. \emph{coderivations} (see \cite{LodayVallette} or \cite{Vallette2} for precise definitions). We call dg properads, dg coproperads, and so on, only by properads, coproperads, and so on.

\section{Curved twisting morphisms}

In this section, we recall the notion of twisting morphisms for
augmented properads and coproperads from \cite{Vallette} and \cite{MerkulovVallette} and the
associated bar-cobar adjunction. To extend these notions to the case
where the properad is not augmented: we introduce the new notion of
\emph{curved coproperad} and of \emph{curved twisting morphism}
between a curved coproperad and a non-necessarily augmented
properad. We also extend the bar and the cobar constructions to this
framework. This provides a functorial cofibrant replacement for properads. We emphasize the fact that the properad is not assumed to be augmented.

\subsection{Twisting morphisms}\label{twmorph}

We recall the classical theory of twisting morphisms between augmented coproperads augmented properads from \cite{MerkulovVallette}.

Let $M$ and $N$ be two $\So$-bimodules. By abuse of notation, we will denote by $M \otimes N$ the infinitesimal composite product of one element of $M$ with one element of $N$ grafted above, that is the space of linear combinations of connected graphs with two vertices, the first one labelled by an element of $M$ and the one above labelled by an element of $N$. This is equal to $M \boxtimes_{(1,1)} N$, with the notation of \cite{MerkulovVallette}. To an operad $\Po$, we associate the \emph{infinitesimal composition product} $\gamma_{(1,\, 1)} : \Po \boxtimes_{(1,\, 1)} \Po \rightarrow \Po$ with the help of $u$ and $\gamma$. Associated to a coproperad $\C$, we define the \emph{infinitesimal decomposition map} $\Delta_{(1,\, 1)} : \C \rightarrow \C \boxtimes_{(1,\, 1)} \C$ by the projection of $\Delta$ (with the help of $\eta$) on $\C \boxtimes_{(1,\, 1)} \C$, or with the above notation, on $\C \otimes \C$.

We recall the convolution product $\star$ on Hom$(\C,\, \Po) := \prod_{m,\, n \geq 0} \mathrm{Hom}_{\K}(\C(m,\, n),\, \Po(m,\, n))$ from \cite{MerkulovVallette}. Let $f,\, g \in$ Hom$(\C,\, \Po)$. We denote by $f \star g$ the composite
$$\C \xrightarrow{\Delta_{(1,\, 1)}} \C \boxtimes_{(1,\, 1)} \C \xrightarrow{f \boxtimes_{(1,\, 1)} g} \Po \boxtimes_{(1,\, 1)} \Po \xrightarrow{\gamma_{(1,\, 1)}} \Po.$$
We define the derivative $\partial$ of degree $-1$ on Hom$(\C,\, \Po)$ by
$$\partial(f) := d_{\Po} \circ f - (-1)^{|f|} f \circ d_{\C}.$$

The convolution product $\star$ on Hom$(\C,\, \Po)$ is a Lie-admissible product (see \cite{MerkulovVallette} for more details). It is stable on the space of equivariant maps from $\C$ to $\Po$ denoted by Hom$_{\So}(\C,\, \Po)$. Then the bracket $[f,\, g] := f \star g - (-1)^{|f||g|} g\star f$ is a Lie bracket on Hom$_{\So}(\C,\, \Po)$.

A morphism of $\So$-modules $\alpha : (\C,\, d_{\C}) \rightarrow (\Po,\, d_{\Po})$ of degree $-1$ in the Lie algebra Hom$_{\So}(\C ,\, \Po)$ is called a \emph{twisting morphism} if it is a solution to the \emph{Maurer-Cartan equation}
$$\partial(\alpha) + \alpha \star \alpha = \partial(\alpha) + \frac{1}{2}[\alpha,\, \alpha] = 0.$$
We denote by Tw$(\C,\, \Po)$ the set of twisting morphisms in Hom$_{\So}(\C,\, \Po)$.

We say that an operad $\Po$ is \emph{augmented} when there is a morphism $\Po \epi I$ of dg properads such that $I \xrightarrow{u} \Po \epi I$ is the identity. It is equivalent to $\Po \cong I \oplus \overline{\Po}$ as dg properads where $\overline{\Po} := \Po/I \cong \ker(\Po \epi I)$. Dually, we say that a coproperad $\C$ is \emph{coaugmented} when there is a morphism $I \mono \C$ of dg coproperads such that $I \mono \C \xrightarrow{\eta} I$ is the identity. It is equivalent to $\C \cong I \oplus \overline{\C}$ as dg coproperads where $\overline{\C} := \ker(\eta)$. When $\Po$ is augmented and $\C$ is coaugmented, we require the twisting morphisms $\alpha$ to satisfy the compositions $\C \xrightarrow{\alpha} \Po \epi I$ and $I \mono \C \xrightarrow{\alpha} \Po$ being equal to $0$. A coaugmented coproperad is called \emph{conilpotent} when for all $x \in \overline{\C}$, there exists an $n > 0$ such that $\overline{\Delta}_{(1,1)}^{n}(x) = 0$, where $\overline{\Delta}_{(1,1)} : \overline{\C} \rightarrow \overline{\C} \boxtimes_{(1,1)} \overline{\C}$ is the primitive part of $\Delta_{(1,1)}$ and where $\overline{\Delta}_{(1,1)}^{n} = (\overline{\Delta}_{(1,1)} \otimes id_{\C}^{\otimes (n-1)}) \circ \overline{\Delta}_{(1,1)}^{n-1}$ (see \cite{LodayVallette} for more details).

When $\Po$ is augmented and $\C$ is conilpotent, we recall from \cite{Vallette2} that the bifunctor Tw$(-,\, -)$ is representable on the left by the \emph{cobar construction} and on the right by the \emph{bar construction}, that is we have the following adjunction
$$\Omega : \textsf{conilpotent dg coprop.} \rightleftharpoons \textsf{augmented dg prop.} : \B$$
and there are natural correspondences
$$\mathrm{Hom}_{\textsf{aug. dg prop.}}(\Omega \C,\, \Po) \cong \mathrm{Tw}(\C,\, \Po) \cong \mathrm{Hom}_{\textsf{conil. dg coprop.}}(\C,\, \B \Po).$$

\subsection{Curved twisting morphism}

We refine the previous section to the case where $\Po$ is not necessarily
augmented. A \emph{curvature} has to be introduced on the level of dg coproperads to encode the default of augmentation. The associated notion is called a \emph{curved coproperad}. We define the notion of \emph{curved twisting morphism} between a curved coproperad and a dg properad as a solution of the \emph{curved Maurer-Cartan equation}.

\subsubsection{\bf Curved coproperad}\label{curvedcoprop}

A \emph{curved coproperad} is a triple $(\C,\, d_{\C},\, \theta)$, where $\C$ is a coproperad, the \emph{predifferential} $d_{\C}$ is a coderivation of $\C$ of degree $-1$ and the \emph{curvature} $\theta : \C \rightarrow I$ is a map of degree $-2$ such that:
\begin{enumerate}
\item[a)] $d_{\C}^{2} = (\theta \otimes id_{\C} - id_{\C} \otimes \theta) \circ \Delta_{(1,1)}$,
\item[b)] $\theta \circ d_{\C} = 0$.
\end{enumerate}

A \emph{morphism between curved coproperads} $(\C,\, d_{\C}, \theta) \rightarrow (\C',\, d_{\C'}, \theta')$ is a morphism of coproperads $f : \C \rightarrow \C'$ such that $d_{\C'} \circ f = f \circ d_{\C}$ and $\theta' \circ f = \theta$. We denote this category by \textsf{curved coprop.}.

\subsubsection{\bf The convolution curved Lie algebra}\label{curvedLiealg}

We define the new notion of \emph{curved Lie algebra} generalizing the notion of dg Lie algebra. A \emph{curved Lie algebra} is a quadriple $(\mathfrak{g},\, [-,\, -],\, d_{\mathfrak{g}},\, \theta)$, where $(\mathfrak{g},\, [-,\, -])$ is a Lie algebra, the predifferential $d_{\mathfrak{g}}$ is a derivation of $\mathfrak{g}$ of degree $-1$ and the curvature $\theta$ is an element of $\mathfrak{g}$ (or equivalently a map $\K \rightarrow \mathfrak{g}$) of degree $-2$ such that:
\begin{enumerate}
\item[a)] $d_{\mathfrak{g}}^{2} = [- ,\, \theta]$;
\item[b)] $d_{\mathfrak{g}} (\theta) = 0$.
\end{enumerate}

Let $(\C,\, d_{\C}, \theta)$ be a curved coproperad and let $(\Po,\, d_{\Po})$ be a dg properad. We fix the element
$$\varTheta := u \circ \theta : \C \xrightarrow{\theta} I \xrightarrow{u} \Po$$
of degree $-2$ in $\textrm{Hom}(\C,\, \Po)$.

\begin{prop}
When ${\bf \C}$ is a curved coproperad and $\Po$ is a dg properad, we have on $\mathrm{Hom}_{\So}(\C,\, \Po) = \prod_{m,\, n \geq 0} \mathrm{Hom}_{\So}(\C(m,\, n),\, \Po(m,\, n))$:
$$\left\{ \begin{array}{lclcl}
\partial^2 & = & [- ,\, \varTheta] & := & (\textrm{-} \star \varTheta) - (\varTheta \star \textrm{-})\\
\partial(\varTheta) & = & 0.&&
\end{array} \right.$$
Then $(\mathrm{Hom}_{\So}(\C,\, \Po),\, [-,\, -],\, \partial,\, \varTheta)$ is a curved Lie algebra, called the \emph{convolution curved Lie algebra}.
\end{prop}

\begin{pf}
We do the computations:
$$\begin{array}{ll}
\partial^2(f) &= d_{\Po} \circ \partial(f) - (-1)^{|\partial(f)|} \partial(f) \circ d_{\C}\\
& = {d_{\Po}}^{2} \circ f - (-1)^{|f|} d_{\Po} \circ f \circ d_{\C} + (-1)^{|f|} (d_{\Po} \circ f \circ d_{\C} - (-1)^{|f|} f \circ {d_{\C}}^2)\\
&= -f \circ {d_{\C}}^2 = - f \circ (\theta \otimes id_{\C} - id_{\C} \otimes \theta) \circ \Delta_{(1,1)} = f \star \varTheta - \varTheta \star f
\end{array}$$
and $\partial(\varTheta) = d_{\Po} \circ u \circ \theta -(-1)^{|\varTheta|} u \circ \theta \circ d_{\C} = 0$ since $d_{\Po} \circ u = 0$ and $\theta \circ d_{\C} = 0$.
$\cqfd$
\end{pf}

An element $\alpha : (\C,\, d_{\C},\, \theta) \rightarrow (\Po,\, d_{\Po})$ of degree $-1$ in the curved Lie algebra Hom$_{\So}(\C ,\, \Po)$ is called a \emph{curved twisting morphism} if it is a solution of the \emph{curved Maurer-Cartan equation}
$$\partial(\alpha) + \alpha \star \alpha = \varTheta.$$
We denote by Tw$(\C,\, \Po)$ the set of curved twisting morphisms in Hom$_{\So}(\C,\, \Po)$.

\begin{rem}
The words ``curved'' and ``curvature'' refer to the geometric context. In that setting, the Maurer-Cartan equation applied to a connection provides the curvature form. The flat case corresponds to the curvature equal to zero, that is to the classical case.
\end{rem}

\subsection{Bar and cobar constructions}

In this section, we extend the bar construction of augmented dg properads
to a \emph{curved bar construction} from dg properads with target in
curved coproperads. In the other way round, we extend the cobar
construction of coaugmented coproperads to coaugmented curved
coproperads. In the algebra case, the cobar construction generalizes the bar construction
of curved algebras given in \cite{PolishchukPositselski} and in
\cite{Positselski} to properads, though it is not immediate that our
constructions are the same, as \cite{PolishchukPositselski,
  Positselski} do not make use of coalgebras.

\subsubsection{\bf Semi-augmented dg properads}\label{dg}

A \emph{semi-augmented dg properad}, or \emph{sdg properad} for short, $(\Po,\, d_{\Po},\, \varepsilon)$ is a dg properad $\Po$ whose underlying $\So$-bimodule is endowed with an augmentation of $\So$-bimodules $\varepsilon : \Po \epi I$, non-necessarily dg or of properads, called \emph{semi-augmentation}. In other words, $\varepsilon$ is a retraction of $\So$-bimodules of the unit $u : I \rightarrow \Po$ and we have an isomorphism $u+inc : I \oplus \overline{\Po} \xrightarrow{\cong} \Po$ of $\So$-bimodules, where $\overline{\Po} := \ker \varepsilon$ and $inc$ is the inclusion $\overline{\Po} \mono \Po$. We denote $\rho := {(u+inc)^{-1}}^{|\overline{\Po}} : \Po \epi \overline{\Po}$. In the following, we do not write the inclusion $inc$ in the formulae. The map $\overline{\gamma} := \rho \circ \gamma : \overline{\Po} \boxtimes \overline{\Po} \rightarrow \overline{\Po}$ is not necessarily associative, even though the composition product $\gamma : \Po \boxtimes \Po \rightarrow \Po$ is associative.

\begin{rem}
The assumption for $\Po$ to have a semi-augmentation $\varepsilon$ is not restrictive
since we are working over a field $\K$ and since we just need to fix a section of $\Po(1,\, 1)$. When $\Po(1,\, 1) = I$, the choice is unique. This is often the case, as it is for the operad encoding unital associative algebras (see Section \ref{resandtransfer}).
\end{rem}

We define on $\overline{\Po}$ the map $d_{\overline{\Po}} := \rho \circ d_{\Po}$, which is a differential since $d_{\Po}$ is a differential and since the differential on $I$ is $0$. The differentials satisfy $\rho \circ d_{\Po} = d_{\overline{\Po}} \circ \rho$. However, we have $d_{\overline{\Po}} \neq d_{\Po}$ in general. 

A \emph{morphism between two sdg properads} $(\Po,\, d_{\Po},\, \varepsilon) \xrightarrow{f} (\Po',\, d_{\Po'},\, \varepsilon')$ is a morphism of dg properads $f : (\Po,\, d_{\Po}) \rightarrow (\Po',\, d_{\Po'})$ such that $\varepsilon' \circ f = \varepsilon$. We define $\bar f := \rho' \circ f : \overline{\Po} \rightarrow \overline{\Po'}$ and we remark that $d_{\overline{\Po'}} \circ \bar f = \bar f \circ d_{\overline{\Po}}$. We denote by \textsf{sdg prop.} the category of semi-augmented dg properads.

\subsubsection{\bf Coaugmented and conilpotent curved coproperads}\label{conilpotent}

When $\C$ is coaugmented, that is $\C$ has a coaugmentation $I \mono \C$ so that $\C \cong I \oplus \overline{\C}$ as coproperads, we require that any twisting morphism $\alpha$ satisfies the compositions $I \mono \C \xrightarrow{\alpha} \Po$ and $\C \xrightarrow{\alpha} \Po \xrightarrow{\varepsilon} I$ to be zero. We denote by \textsf{coaug. curved coprop.} the category of coaugmented curved coproperads and by \textsf{conil. curved coprop.} the category of conilpotent curved coproperads (see Section \ref{twmorph}).\\

We construct a pair of functors
$$\B : \textsf{sdg prop.} \rightleftharpoons \textsf{coaug. curved coprop.} : \Omega.$$
Let $M$ be an $\So$-bimodule. The notation $\F(M)$, resp. $\F^{c}(M)$, stands for the \emph{free properad} on $M$, resp. the \emph{cofree coproperad} on $M$. A derivation on $\F(M)$, resp. a coderivation on $\F^{c}(M)$, is characterized by its restriction on $M$, resp. by its image on $M$. The notation $sM$, resp. $s^{-1}M$, stands for the \emph{homological suspension}, resp. the \emph{homological desuspension}, of the $\So$-bimodule $M$. We refer to \cite{Vallette2} for more details.

\subsubsection{\bf Curved bar construction of a sdg properad}

The \emph{bar construction of the sdg properad $(\Po,\, d_{\Po},\, \varepsilon)$} is given by the conilpotent curved coproperad
$$\B \Po := (\F^{c}(s\overline{\Po}),\, d_{bar},\, \theta_{bar}).$$
The predifferential is defined by $d_{bar} := d_{1} + d_{2}$, where $d_{2}$ is the unique coderivation of degree $-1$ which extends the map
$$\F^{c}(s\overline{\Po}) \epi \F^{c}(s\overline{\Po})^{(2)} \cong s^{2} (\overline{\Po} \boxtimes_{(1,\, 1)} \overline{\Po}) \xrightarrow{s^{-1} \overline{\gamma}} s\overline{\Po}$$
where $\overline{\gamma} := \rho \circ \gamma : \overline{\Po} \boxtimes_{(1,\, 1)} \overline{\Po} \rightarrow \overline{\Po}$ and $d_{1}$ is the unique coderivation of degree $-1$ which extends the map
$$\F^{c}(s\overline{\Po}) \epi s\overline{\Po} \xrightarrow{id_{s} \otimes d_{\overline{\Po}}} s\overline{\Po}.$$
The curvature $\theta_{bar}$ is the map of degree $-2$
$$\F^{c}(s\overline{\Po}) \epi s\overline{\Po} \oplus \F^{c}(s\overline{\Po})^{(2)} \cong s\overline{\Po} \oplus s^{2} (\overline{\Po} \boxtimes_{(1,\, 1)} \overline{\Po}) \xrightarrow{s^{-1} d_{\Po} + s^{-2} \gamma} \Po \xrightarrow{\varepsilon} I.$$

\begin{lem}
The predifferential and the curvature satisfy
\begin{enumerate}
\item[a)] ${d_{bar}}^{2} = (\theta_{bar} \otimes id - id \otimes \theta_{bar}) \circ \Delta_{(1,\, 1)}$;
\item[b)] $\theta_{bar} \circ d_{bar} = 0$.
\end{enumerate}
\end{lem}

\begin{pf}
First we can restrict the proof of the equality $a)$ and $b)$ to $\F^{c}(s\overline{\Po})^{(\leq 3)}$ since ${d_{bar}}^2$ and $(\theta_{bar} \otimes id - id \otimes \theta_{bar}) \circ \Delta_{(1,\, 1)}$ are coderivations and since $\theta_{bar}$ is non zero only on $\F^{c}(s\overline{\Po})^{(2)}$.

The composite
$$\hspace{-3cm} \F^{c}(s\overline{\Po})^{(\leq 3)} \xrightarrow{{d_{bar}}^{|\F^{c}(s\overline{\Po})^{(\leq 2)}} - \left[(\theta_{bar} \otimes id - id \otimes \theta_{bar}) \circ \Delta_{(1,\, 1)} \right]^{|I \otimes s\overline{\Po} \oplus s\overline{\Po} \otimes I}} \F^{c}(s\overline{\Po})^{(\leq 2)} \oplus$$
$$\hspace{7cm} (I \otimes s\overline{\Po} \oplus s\overline{\Po} \otimes I) \xrightarrow{({d_{bar}}^{|s\overline{\Po}} - \theta_{bar}) + \gamma_{|I \otimes s\overline{\Po} \oplus s\overline{\Po} \otimes I}} I \oplus s\overline{\Po}$$
equals to $\left({d_{bar}}^{2} - (\theta_{bar} \otimes id - id \otimes \theta_{bar}) \circ \Delta_{(1,\, 1)} - \theta_{bar} \circ d_{bar}\right)^{|I \oplus s \overline{\Po}}$ and to ${(d_{\gamma + d_{\Po}})^2}^{|I \oplus s \overline{\Po}}$ where $d_{\gamma + d_{\Po}}$ is the unique coderivation of degree $-1$ on $\F^{c}(s\Po)$ which extends the map
$$\F^{c}(s\Po) \epi \overline{\F^{c}}(s\Po)^{(\leq 2)} \cong s \Po \oplus s^2 \Po \boxtimes_{(1,\, 1)} \Po \xrightarrow{id_{s} \otimes d_{\Po} + s^{-1} \gamma} s\Po.$$
Moreover, since $\gamma$ is associative and $d_{\Po}$ is a differential, we have ${d_{\gamma + d_{\Po}}}^{2} = 0$. Thus
$${d_{bar}}^{2} - (\theta_{bar} \otimes id - id \otimes \theta_{bar}) \circ \Delta_{(1,\, 1)} - \theta_{bar} \circ d_{bar} = 0,$$
that is, due to the degree
$$\left\{\begin{array}{l}
{d_{bar}}^{2} = (\theta_{bar} \otimes id - id \otimes \theta_{bar}) \circ \Delta_{(1,\, 1)}\\
\theta_{bar} \circ d_{bar} = 0.
\end{array} \right.$$
$\cqfd$
\end{pf}

\begin{lem}
The bar construction is a functor $\B : \mathsf{sdg\ prop.} \rightarrow \mathsf{conil.\ curved\ coprop.}$.
\end{lem}

\begin{pf}
Let $f : (\Po,\, d_{\Po},\, \varepsilon) \rightarrow (\Po',\, d_{\Po'},\, \varepsilon')$ be a morphism of sdg properads. It induces a morphism of dg $\So$-bimodules $\bar{f} : \overline{\Po} \rightarrow \overline{\Po'}$. The map $\F^{c}(\bar f) : \F^{c}(s\overline{\Po}) \rightarrow \F^{c}(s\overline{\Po'})$ is a map of coproperads by construction. The morphism $\bar f$ commutes with $\overline{\gamma}_{\Po}$ and $\overline{\gamma}_{\Po'}$, thus $\F^{c}(\bar f)$ commutes with the predifferentials. For a similar reason $\theta'_{bar} \circ \F^{c}(\bar f) = \theta_{bar}$.
$\cqfd$
\end{pf}

\subsubsection{\bf Cobar construction of a coaugmented curved coproperad}

The \emph{cobar construction of the coaugmented curved coproperad $(\C,\, d_{\C}, \theta)$} is given by the sdg properad
$$\Omega \C := (\F(s^{-1}\overline{\C}),\, d := d_{0} + d_{1} - d_{2},\, \varepsilon).$$
The term $d_{0}$ is the unique derivation of degree $-1$ which extends the map
$$s^{-1}\overline{\C} \xrightarrow{s\theta} I \mono \F(s^{-1}\overline{\C}).$$
The term $d_{1}$ is the unique derivation of degree $-1$ which extends the map
$$s^{-1}\overline{\C} \xrightarrow{id_{s^{-1}} \otimes d_{\overline{\C}}} s^{-1}\overline{\C} \mono \F(s^{-1}\overline{\C}).$$
The term $d_{2}$ is the unique derivation of degree $-1$ which extends the infinitesimal decomposition map of $\overline{\C}$, up to desuspension:
$$s^{-1}\overline{\C} \xrightarrow{s^{-1} \overline{\Delta}_{(1,1)}} s^{-2}\F^{c}(\overline{\C})^{(2)} \cong \F(s^{-1}\overline{\C})^{(2)} \mono \F(s^{-1}\overline{\C}).$$
The semi-augmentation $\varepsilon$ is the natural projection $\F(s^{-1} \overline{\C}) = I \oplus s^{-1}\overline{\C} \oplus \cdots \epi I$. It is an augmentation of properads but it is not an augmentation of dg properads in general.

\begin{lem}
The derivation $d$ on $\F(s^{-1}\overline{\C})$ satisfies $d^{2} = 0$.
\end{lem}

\begin{pf}
First of all, if we define the weight on $\C$ by $\C^{(0)} = I$, $\C^{(1)} = \overline{\C}$ and $\C^{(n)} = 0$ when $n \neq 0,\, 1$ and extend it to $\F (s^{-1}\overline{\C})$, we get that the map $d_{0}$ is of weight $-1$, the map $d_{1}$ is of weight $0$ and the map $d_{2}$ is of weight $1$. Thus, the term $d^{2}$ split in the following way
$$d^{2} = \underbrace{{d_{0}}^{2}}_{\textrm{weight} = -2} + \underbrace{d_{0} d_{1} + d_{1} d_{0}}_{\textrm{weight} = -1} + \underbrace{{d_{1}}^{2} - d_{0} d_{2} - d_{2} d_{0}}_{\textrm{weight} = 0} - \underbrace{(d_{1} d_{2} + d_{2} d_{1})}_{\textrm{weight} = 1} +\underbrace{{d_{2}}^{2}}_{\textrm{weight} = 2}.$$
So, we have to show that each group of terms is equal to zero. The term $d_{0}^{2}$ is zero by the Koszul sign rule. The sum $d_{0} d_{1} + d_{1} d_{0}$ is zero since $\theta \circ d_{\C} = 0$ and $d_{\C}$ is zero on $I$ and by the Koszul sign rule. The equality $d_{\C}^{2} = (\theta \otimes id_{\C} - id_{\C} \otimes \theta) \circ \Delta_{(1,1)}$ and the Koszul sign rule give ${d_{1}}^{2} - d_{0} d_{2} - d_{2} d_{0} = 0$. The equality $d_{1} d_{2} + d_{2} d_{1} = 0$ is due to the fact that $d_{\C}$ is a coderivation. Finally ${d_{2}}^{2} = 0$ by ``coassociativity'' of $\overline{\Delta}_{(1,1)}$ and by the Koszul sign rule.
$\cqfd$
\end{pf}

\begin{lem}
The cobar construction is a functor $\Omega : \mathsf{coaug.\ curved\ coprop.} \rightarrow \mathsf{sdg\ prop.}$.
\end{lem}

\begin{pf}
Let $f : (\C,\, d_{\C},\, \theta) \rightarrow (\C',\, d_{\C'},\, \theta')$ be a morphism between coaugmented curved coproperads. The map $\F(f) : \F(s^{-1}\overline{\C}) \rightarrow \F(s^{-1}\overline{\C'})$ is a map of properads by construction and $d'_{2} \circ \F(f) = \F(f) \circ d_{2}$ since $f$ is a morphism of coproperads. The equality $d_{\C'} \circ f = f \circ d_{\C}$ implies $d'_{1} \circ \F(f) = \F(f) \circ d_{1}$, the equality $\theta' \circ f = \theta$ implies $d'_{0} \circ \F(f) = \F(f) \circ d_{0}$ and then $\F(f)$ commutes with the differential.
$\cqfd$
\end{pf}

\subsection{Bar-cobar adjunction}\label{sectbarcobaradjunction}

The cobar construction on conilpotent curved coproperads and the bar construction on dg properads represent the bifunctor of curved twisting morphisms and form a pair of adjoint functors. The counit of adjunction provides a cofibrant replacement functor for dg properads.

\begin{thm}\label{barcobaradjunction}
For any conilpotent curved coproperad $\C$ and for any sdg properad $\Po$,
there is are natural correspondences
$$\mathrm{Hom}_{\mathsf{sdg\ prop.}}(\Omega \C,\, \Po) \cong \mathrm{Tw}(\C,\, \Po) \cong \mathrm{Hom}_{\mathsf{coaug.\ curved\ coprop.}}(\C,\, \B \Po).$$
\end{thm}

\begin{pf}
We make the first bijection explicit. A morphism of properads $f_{\alpha} : \F(s^{-1}\overline{\C}) \rightarrow \Po$ is uniquely determined by a map $s\alpha : s^{-1}\overline{\C} \rightarrow \Po$ of degree $0$ such that $s^{-1}\overline{\C} \xrightarrow{s\alpha} \Po \xrightarrow{\varepsilon} I$ is $0$ (by definition of dg properads' morphisms, see \ref{dg}) or equivalently, by a map $\alpha : \C \rightarrow \Po$ of degree $-1$ satisfying $I \mono \C \xrightarrow{\alpha} \Po$ and $\C \xrightarrow{\alpha} \Po \xrightarrow{\varepsilon} I$ are zero (condition for twisting morphisms when $\C$ is coaugmented, see \ref{conilpotent}).

Moreover, $f_{\alpha}$ commutes with the differentials if and only if the following diagram commutes
$$\xymatrix{s^{-1}\overline{\C} \ar[r]^{s\alpha} \ar[d]_{d_{0} + d_{1} - d_{2}} & \Po \ar[r]^{d_{\Po}} & \Po\\
\F(s^{-1}\overline{\C}) \ar[rr]_{\F(s\alpha)} && \F(\Po), \ar[u]_{\widetilde{\gamma}}}$$
where $\widetilde{\gamma}$ is induced by $\gamma$. We have
$$\begin{array}{l}
d_{\Po} \circ (s\alpha) = s(d_{\Po} \circ \alpha)\\
\widetilde{\gamma} \circ \F(s\alpha) \circ d_0 = u \circ (s\theta) = s(u \circ \theta)\\
\widetilde{\gamma} \circ \F(s\alpha) \circ d_1 = s\alpha \circ (id_{s^{-1}} \otimes d_{\C}) = - s(\alpha \circ d_{\C})\\
\widetilde{\gamma} \circ \F(s\alpha) \circ d_2 = \gamma \circ (s\alpha \boxtimes_{(1,\, 1)} s\alpha) \circ s^{-1} \Delta_{(1,\, 1)} = s(\gamma \circ (\alpha \boxtimes_{(1,\, 1)} \alpha) \circ \Delta_{(1,\, 1)}).
\end{array}$$
Thus the commutativity of the previous diagram is equivalent to the equality
$$u \circ \theta - \alpha \circ d - \gamma \circ (\alpha \boxtimes_{(1,\, 1)} \alpha) \circ \Delta_{(1,\, 1)} = d_{\Po} \circ \alpha,$$
that is $\partial(\alpha) + \alpha \star \alpha = \varTheta$.

We now make the second bijection explicit. A morphism of coaugmented coproperads $g_{\alpha} : \C \rightarrow \F^{c}(s\overline{\Po})$ is uniquely determined by a map $s\alpha : \C \rightarrow s\overline{\Po}$ which sends $I$ to $0$, that is by a map $\alpha : \C \rightarrow \Po$ of degree $-1$ satisfying $I \mono \C \xrightarrow{\alpha} \Po$ and $\C \xrightarrow{\alpha} \Po \xrightarrow{\varepsilon} I$ are zero.

Moreover, $g_{\alpha}$ commutes with the predifferential and with the curvature if and only if the following diagrams commute
$$\xymatrix@C=68pt@R=10pt{\C \ar[r]^{s\alpha + (s\alpha \otimes s\alpha) \circ \Delta_{(1,\, 1)} \hspace{1cm}} \ar[dd]_{d_{\C}} & s\overline{\Po} \oplus s\overline{\Po} \boxtimes_{(1,\, 1)} s\overline{\Po} \ar[dd]^{d_{bar} = d_{1} + d_{2}} &\\
& & \textrm{and}\\
\C \ar[r]_{s\alpha} & s\overline{\Po} &}
\hspace{2cm} \mbox{\xymatrix@R=10pt{\C \ar[dd]_{\theta} \ar[r]^{g_{\alpha}} & \B \Po \ar[ddl]^{\theta_{bar}}\\
&\\
I. &}}$$
Since $\alpha \star \alpha = -(s^{-1}inc) \circ d_{2} \circ (s\alpha \otimes s\alpha) \circ \Delta_{(1,\, 1)} + u \circ \theta_{bar} \circ g_{\alpha}$, the commutativity of the diagrams gives $\partial(\alpha) + \alpha \star \alpha = \varTheta$. Moreover, the projections of the curved Maurer-Cartan equation on $\overline{\Po}$ and on $I$ give the two commutative diagrams. This concludes the proof.
$\cqfd$
\end{pf}

\begin{ex}
\begin{itemize}
\item[]
\item To the identity morphism $id_{\B \Po} : \B \Po \rightarrow \B \Po$ of coaugmented curved coproperads corresponds the curved twisting morphism $\pi : \B \Po \rightarrow \Po$ defined by $\F^{c}(s\overline{\Po}) \epi s\overline{\Po} \cong \overline{\Po} \mono \Po$.
\item To the identity morphism $id_{\Omega \C} : \Omega \C \rightarrow \Omega \C$ of properads corresponds the curved twisting morphism $\iota : \C \rightarrow \Omega \C$ defined by $\C \rightarrow \overline{\C} \cong s^{-1}\overline{\C} \mono \F(s\overline{\C})$.
\end{itemize}
\end{ex}

\begin{lem}\label{universaltwmorph}
For any conilpotent curved coproperad $\C$ and for any sdg properad $\Po$, every curved twisting morphism $\alpha : \C \rightarrow \Po$ factors through the universal curved twisting morphisms $\pi$ and $\iota$:
$$\xymatrix{& \Omega \C \ar@{-->}[dr]^{f_{\alpha}} &\\
\C \ar[rr]^{\alpha} \ar[ur]^{\iota} \ar@{-->}[dr]_{g_{\alpha}} && \Po\\
& \B \Po, \ar[ur]_{\pi} &}$$
where $f_{\alpha}$ is a morphism of sdg properads and $g_{\alpha}$ is a morphism of conilpotent curved coproperads.
\end{lem}

\begin{pf}
The dashed arrows are just the images of $\alpha$ by the two bijections of Proposition \ref{barcobaradjunction}.
$\cqfd$
\end{pf}

\subsubsection{\bf Weight filtration}

We say that a dg $\So$-bimodule $M$ is \emph{weight filtered differential graded}, or \emph{wfdg} for short, when it is endowed with a filtration of dg $\So$-bimodules $F_{\omega}M$, $\omega \in \mathbb{N}$. When $M$ is a (co)properad, we assume that the (co)product preserves the filtration. In the weight filtered setting, we also assume that the twisting morphisms preserve the filtration. A wfdg properad $\Po$ is called \emph{connected} when $F_{0}\Po = I\ (= Im (u))$.

We endow any free properad $\F(V)$ with a weight grading given by the number of generators. This induces a weight filtration on any properad $\F(V)/(R)$ defined by generators and relations. Sub-coproperads of $\F^{c}(V)$ are also weight filtered by the number of generators. When $\Po$ is a wfdg properad, we also endow $\B \Po$ with a weight filtration. An element in $\B \Po$ is a connected graph whose vertices are labelled by elements $\mu_{i}$ of $\overline{\Po}$. It is in the component of weight $\omega$ of $\B \Po$ if there exist $\omega_{i}$ such that any $\mu_{i}$ is in the component of weight $\omega_{i}$ of $\overline{\Po}$ and $\sum \omega_{i} \leq \omega$. Similarly, we endow $\Omega \C$ with a weight filtering when $\C$ is weight filtered.

The curved twisting morphism $\pi$ preserves the weight filtration.

\begin{thm}\label{barcobarresprop}
Let $(\Po,\, d_{\Po},\, \varepsilon)$ be a connected wfdg semi-augmented properad. The counit of the bar-cobar adjunction is a quasi-isomorphism of wfdg semi-augmented properads, that is the \emph{bar-cobar construction} $\Omega \B \Po$ is a resolution of $\Po$
$$\Omega \B \Po \qiso \Po.$$
When $\Po$ is concentrated in non-negative degree, the bar-cobar construction is a cofibrant properad for the model category defined in Appendix A of \cite{MerkulovVallette2}.
\end{thm}

\begin{pf}
We work in the model category defined in Appendix A of \cite{MerkulovVallette2}. Since $\Omega \B \Po$ is quasi-free, the remark after Corollary 40 of \cite{MerkulovVallette2} gives that $\Omega \B \Po$ is cofibrant when we assume that $\Po$ is non-negatively homologically graded.

As explained in the previous section, $\Omega \B \Po = (\F (s^{-1} \overline{\F} (s \overline{\Po})),\, d = d_{0} + d_{1} - d_{2})$ is weight filtered by $F_{p}$ when $\Po$ is weight filtered. We have
$$d_{0} : F_{p} \rightarrow F_{p-1} \textrm{ and } d_{1} : F_{p} \rightarrow F_{p} \textrm{ and } d_{2} : F_{p} \rightarrow F_{p},$$
where $d_{0}$ is induced by $\theta_{bar}$, $d_{1}$ is induced by $d_{bar}$ and $d_{2}$ is induced by the coproduct on $\F^{c}(s \overline{\Po})$. So $F_{p}$ is a filtration of chain complexes, it is exhaustive and bounded below and we can apply the classical theorem of convergence of spectral sequences (cf. Theorem 5.5.1 of \cite{Weibel}) to obtain
$$E_{p,q}^{\bullet} \Rightarrow \mathrm{H}_{p+q}(\Omega \B \Po).$$
We endow $\Po$ with a filtration $F'_{p}$ induced by the weight. This is a filtration of chain complexes since $d_{\Po}$ preserves the weight filtration. The filtration $F'_{p}$ is exhaustive and bounded below so we can apply the classical theorem of convergence of spectral sequences (cf. Theorem 5.5.1 of \cite{Weibel}) to obtain
$${E'}^{ \bullet}_{p,q} \Rightarrow \mathrm{H}_{p+q}(\Po).$$
The counit of the bar-cobar adjunction preserves the filtration and induces a map of spectral sequences $E^{ \bullet}_{p,q} \rightarrow {E'}^{ \bullet}_{p,q}$. Moreover, $E_{\bullet, \bullet}^{0} = \Omega \B (gr\Po)$. The graded properad $gr\Po$ associated to the filtration $F'_{p}$ on $\Po$ is always augmented and connected (in the sense of \cite{Vallette}, that is $gr\Po$ is weight graded and $gr\Po^{(0)} = I$). By Theorem 5.8 of \cite{Vallette2}, we get that $E_{p,q}^{1} = gr^{(p)} (\mathrm{H}_{p+q}(gr\Po))$. Thus the counit of the bar-cobar adjunction induces an isomorphism of spectral sequences $E^{ \bullet}_{p,q} \rightarrow {E'}^{ \bullet}_{p,q}$ when $\bullet \geq 1$. Since ${E'}^{ \bullet}_{p,q} \Rightarrow \mathrm{H}_{p+q}(\Po)$, the same is true for $E^{\bullet}_{p,q}$ and the morphism $\Omega \B \Po \qiso \Po$ is a quasi-isomorphism.
$\cqfd$
\end{pf}

\begin{rems}
\begin{itemize}
\item[]
\item In \cite{Positselski2}, Positselski defined a bar construction and a cobar construction between curved dg algebras and curved dg coalgebras. The curvatures on both sides encode the default of augmentation or of coaugmentation. In this paper, we are interested only by the default of augmentation and the picture becomes asymmetric. When we reduce our bar construction and our cobar construction to semi-augmented algebras and curved coalgebras, we recover the particular case of \cite{Positselski2} where the curved coalgebras are coaugmented.
\item In \cite{Nicolas}, Nicol\`as proved a similar bar-cobar adjunction on the level of algebras and coalgebras. But the picture is dual. The bar construction goes from \emph{curved associative algebras} to \emph{conilpotent graded-augmented coalgebras} (see \cite{Nicolas} for the precise definitions) and the cobar construction goes the other way around. In his case, the curvature does not control the default of augmentation with respect to the composition product and with respect to the dg setting, but only with respect to the dg setting. In the spirit of \cite{Nicolas}, we should say the dual statement: the default of augmentation with respect to the dg setting measures the curvature.
\end{itemize}
\end{rems}

\subsubsection{\bf Homotopy Frobenius algebras}

A \emph{unital and counital Frobenius algebras} is a quintuple $(A,\, \mu,\, \Delta,\, u,\, \eta)$ where $A$ is a vector space, $\mu : A \otimes A \rightarrow A$ is a commutative and associative product, $\Delta : A \rightarrow A \otimes A$ is a cocommutative and coassociative coproduct, $u : \K \rightarrow A$ is a unit for the product and $\eta : A \rightarrow \K$ is a counit for the coproduct such that the product  $\mu = \vcenter{\xymatrix@M=0pt@R=4pt@C=4pt{\ar@{-}[dr] && \ar@{-}[dl]\\
           & \ar@{-}[d] &\\
           &&}}$ and the coproduct $\Delta = \vcenter{\xymatrix@M=0pt@R=4pt@C=4pt{&&\\
           & \ar@{-}[u] &\\
           \ar@{-}[ur] && \ar@{-}[ul]}}$ satisfy the \emph{Frobenius relation}
$$\vcenter{\xymatrix@M=0pt@R=6pt@C=6pt{& \ar@{-}[d] &  &\ar@{-}[d]  \\  & \ar@{-}[dl]
 \ar@{-}[dr] &  & \ar@{-}[dl]\\
\ar@{-}[d]&  & \ar@{-}[d] & \\ & & &}} = \vcenter{\xymatrix@M=0pt@R=6pt@C=6pt{ \ar@{-}[dr] &  &\ar@{-}[dl]  \\  &\ar@{-}[d] & \\
& \ar@{-}[dl]\ar@{-}[dr]& \\ & &}} = \vcenter{\xymatrix@M=0pt@R=6pt@C=6pt{ \ar@{-}[d] & &\ar@{-}[d] &  \\ \ar@{-}[dr] &
& \ar@{-}[dl] \ar@{-}[dr]& \\
& \ar@{-}[d]&  & \ar@{-}[d] \\ & & &}}\ .$$
In operadic terms, we get that $A$ is an algebra over the properad\\
$\ucFrob := $
$${\small{\F (\vcenter{
    \xymatrix@M=0pt@R=3pt@C=3pt{&&\\
      & \ar@{{*}}[u] \ar@{-}[d]\\
      &}},\, \vcenter{
    \xymatrix@M=0pt@R=3pt@C=3pt{&&\\
      & \ar@{{*}}[d] \ar@{-}[u]\\
      &}},\,
  \vcenter{\xymatrix@M=0pt@R=4pt@C=4pt{\ar@{-}[dr] && \ar@{-}[dl]\\
      & \ar@{-}[d] &\\
      &&}},\, \vcenter{\xymatrix@M=0pt@R=4pt@C=4pt{&&\\
      & \ar@{-}[u] &\\
      \ar@{-}[ur] && \ar@{-}[ul]
      }})/ (\vcenter{ \xymatrix@M=0pt@R=4pt@C=4pt{
      & & & & \\ \ar@{-}[ddrr] & & \ar@{-}[dl] & & \ar@{-}[ddll]\\
      &  & &  & \\
      & & \ar@{-}[d] & &\\
      & & }} - \vcenter{
    \xymatrix@M=0pt@R=4pt@C=4pt{& & & & \\
      \ar@{-}[drdr] & &\ar@{-}[dr] & & \ar@{-}[dldl]  \\
      & & & & \\
      & &\ar@{-}[d] & & \\
      & & }},\, \vcenter{ \xymatrix@M=0pt@R=4pt@C=4pt{
      & & & & \\ & & \ar@{-}[u] & &\\
      &  & &  & \\
      \ar@{-}[uurr] & & \ar@{-}[ul] & & \ar@{-}[uull]\\
      & & }} - \vcenter{
    \xymatrix@M=0pt@R=4pt@C=4pt{& & & & \\
      & &\ar@{-}[u] & & \\
      & & & & \\
      \ar@{-}[urur] & &\ar@{-}[ur] & & \ar@{-}[ulul]  \\
      & & }},\, \vcenter{ \xymatrix@M=0pt@R=4pt@C=4pt{
      & & & & \\
      & & & &\\
      & \ar@{{*}}[u] \ar@{-}[dr] & & \ar@{-}[dl] & \\
      & & \ar@{-}[d] & &\\
      & & }} - |, \vcenter{
    \xymatrix@M=0pt@R=4pt@C=4pt{& & & & \\
      & & & &\\
      & \ar@{-}[dr] & & \ar@{{*}}[u] \ar@{-}[dl] & \\
      & &\ar@{-}[d] & & \\
      & & }} - |,\, \vcenter{ \xymatrix@M=0pt@R=4pt@C=4pt{
      & & & & \\
      & & \ar@{-}[u] & &\\
      & \ar@{{*}}[d] \ar@{-}[ur] & & \ar@{-}[ul] & \\
      & & & &\\
      & & }} - |, \vcenter{
    \xymatrix@M=0pt@R=4pt@C=4pt{& & & & \\
      & & \ar@{-}[u] & &\\
      & \ar@{-}[ur] & & \ar@{{*}}[d] \ar@{-}[ul] & \\
      & & & & \\
      & & }} - |,\, \vcenter{\xymatrix@M=0pt@R=4pt@C=4pt{& \ar@{-}[d] &  &\ar@{-}[d]  \\  & \ar@{-}[dl]
 \ar@{-}[dr] &  & \ar@{-}[dl]\\
\ar@{-}[d]&  & \ar@{-}[d] & \\ & & &}} - \vcenter{\xymatrix@M=0pt@R=4pt@C=4pt{ \ar@{-}[dr] &  &\ar@{-}[dl]  \\  &\ar@{-}[d] & \\
& \ar@{-}[dl]\ar@{-}[dr]& \\ & &}},\, \vcenter{\xymatrix@M=0pt@R=4pt@C=4pt{ \ar@{-}[dr] &  &\ar@{-}[dl]  \\  &\ar@{-}[d] & \\
& \ar@{-}[dl]\ar@{-}[dr]& \\ & &}} - \vcenter{\xymatrix@M=0pt@R=4pt@C=4pt{ \ar@{-}[d] & &\ar@{-}[d] &  \\ \ar@{-}[dr] &
& \ar@{-}[dl] \ar@{-}[dr]& \\
& \ar@{-}[d]&  & \ar@{-}[d] \\ & & &}} )}.}$$
This properad is not augmented but the Theorem \ref{barcobarresprop} applies and we get as a corollary:

\begin{thm}
The bar-cobar resolution on $\ucFrob$ is a cofibrant resolution of the properad $\ucFrob$, that is
$$\Omega \B \, \ucFrob \qiso \ucFrob.$$
\end{thm}

We define a  \emph{$\ucFrob$-algebra up to homotopy} as an algebra over this resolution. As proved in \cite{Abrams, Kock}, the datum of a \emph{2-dimensional topological quantum field theory}, \emph{2d-TQFT} for short is equivalent to a unital and counital Frobenius algebra structure. Therefore, we should be able to use this to study 2d-TQFT with homotopy methods.\\

There is an interesting application in differential geometry. With the present resolution of $\ucFrob$ and with the methods of \cite{Wilson}, one endows the differential forms $\Omega (M)$ on a closed, oriented manifold $M$ with a structure of $\ucFrob$-algebra up to homotopy, which extends the commutative algebra structure on the $\Omega (M)$ and which induces the $\ucFrob$-algebra structure on the cohomology $\mathrm{H}^{\bullet} (M)$. 

\section{Curved Koszul duality theory} \label{curvedKD}

We extend the Koszul duality theory for homogeneous quadratic properads
\cite{Vallette2} and quadratic-linear properads \cite{Ga-CaToVa} to \emph{inhomogeneous quadratic properads} with a quadratic, linear and constant presentation. When the properad is inhomogeneous quadratic, it is not necessarily augmented. Therefore we introduce a Koszul dual coproperad endowed with a curvature, which measures this failure. As explained in Section \ref{properad}, an associative algebra is a particular case of properad. Hence this section applies to associative algebras as well to recover the construction given by \cite{Positselski} and \cite{PolishchukPositselski}. However, the
presentation given here is slightly different and more general: it works without any finiteness
assumption. We end the section with a Poincaré-Birkhoff-Witt theorem for properads.

\subsection{Inhomogeneous quadratic properad}\label{inhomoquadprop}

An \emph{inhomogeneous quadratic properad} is a properad $\Po$ which admits a presentation of the form $\Po = \F(V)/(R)$, where $V$ is a degree graded $\So$-bimodule and $(R)$ is the ideal generated by a degree graded $\So$-bimodule $R \subset I \oplus V \oplus \F(V)^{(2)}$. The superscript $(2)$ indicates the weight degree. We require that $R$ is a direct sum of (homological) degree homogeneous subspaces. Thus the properad $\Po$ is degree graded and has a weight filtration induced by the $\So$-bimodule of generators $V$. We assume further that the following conditions hold:
\begin{itemize}
\item[(I)] The space of generators is minimal, that is $R \cap \{I \oplus V\} = \{ 0\}$.
\item[(II)] The space of relations is maximal, that is $(R)\cap \{I \oplus V \oplus \F(V)^{(2)}\} = R$.
\end{itemize}

Let $\q : \F(V) \epi \F(V)^{(2)}$ be the canonical projection and let $\q R \subset \F(V)^{(2)}$ be the image under $\q$ of $R$. We consider the quadratic properad $\q\Po := \F(V)/(\q R)$. Since $R \cap \{I \oplus V\} = \{0\}$, there exists a map $\varphi : \q R \rightarrow I \oplus V$ such that $R$ is the graph of $\varphi$:
\begin{eqnarray*}
R & = & \{X - \varphi(X),\, X \in \q R\}\\
& = & \{X - \varphi_{1}(X) + \varphi_{0}(X),\, X\in \q R,\, \varphi_{1}(X)\in V,\, \varphi_{0}(X)\in \K\}.
\end{eqnarray*}
The weight grading on the free properad $\F(V)$ induces the following filtration on $\Po$
$$F_{p} := \pi \left(\oplus_{\omega \leq p} \F(V)^{(\omega)}\right),$$
where $\pi$ stands for the canonical projection $\F(V) \epi \Po$. We denote the associated graded properad by $gr(\Po)$. The relations $\q R$ hold in $gr(\Po)$. Therefore, there exists an epimorphism of graded properads
$$p : \q\Po \epi gr(\Po).$$
We assume in all the paper that every inhomogeneous quadratic properad is semi-augmented in the sense of Section \ref{dg}. We recall that $sV$ stands for the homological suspension of $V$. Recall that the \emph{Koszul dual coproperad of the homogeneous quadratic properad} $\q\Po$ is the coproperad generated by $sV$ with relations in $s^{2}\q R$ (see Section 2.2 of \cite{Vallette}) denoted
$$\q\Po^{\ash} := \C(sV,\, s^{2}\q R) = I \oplus sV \oplus s^{2}\q R \oplus \cdots.$$
It is a subcoproperad of the cofree coproperad $\F^{c}(sV)$ on $sV$. In the cofree coproperad $\F^{c}(V)$, the weight of an element corresponds to the number of generating elements from $V$ used to write it. There exists a unique coderivation $\tilde{d} : \qPs \rightarrow \F^{c}(sV)$ of degree $-1$ (see Section 3.2 in \cite{MerkulovVallette}) which extends the map
$$\qPs \epi s^{2}\q R \xrightarrow{s^{-1}\varphi_{1}} sV.$$
Moreover, we denote by $\theta : \qPs \rightarrow I$ the map of degree $-2$
$$\qPs \epi s^{2}\q R \xrightarrow{s^{-2}\varphi_{0}} I.$$
There exists a unique coderivation $D : \qPs \rightarrow \F^{c}(sV)$ of degree $-2$ which extends the map
$$\qPs \epi {\qPs}^{(3)} \xrightarrow{(s^{-2}\varphi_{0} \otimes id_{sV} - id_{sV} \otimes s^{-2}\varphi_{0})\circ \Delta_{(1, 1)}} sV.$$

\begin{lem}\label{Dcoder}
The coderivation $D$ is equal to $(\theta \otimes id_{\qPs} - id_{\qPs} \otimes \theta) \circ \Delta_{(1, 1)}$.
\end{lem}

\begin{pf}
By definition, $\theta$ is eventually non-zero only on $s^{2}\q R$, that is on elements of the form
$$\vcenter{\xymatrix@R=8pt@C=4pt{& \ar[d] &\\
*{} \ar@{--}[rr] \ar@{--}[ddd] & *{} & *{} \ar@{--}[ddd]\\
& *+[F-,]{\hspace{.3cm}} &\\
& *+[F-,]{\hspace{.3cm}} &\\
*{} \ar@{--}[rr] & *{} \ar[d] & *{}\\
&&}} =: \mu \in \Po(1,\, 1).$$
Consider an element $\nu$ in $\qPs$ containing the dashed box $\mu$ in its graphical representation. When there is something on the top of $\mu$, the term with $\theta (\mu)$ instead of $\mu$ in $\nu$ contributes to the coderivation applied to $\nu$ and when there is something on the bottom of $\mu$, the term with $\theta (\mu)$ instead of $\mu$ in $\nu$ contributes to the coderivation applied to $\nu$ with a coefficient $-1$. These are the only possibilities for the contribution of $\theta (\mu)$. Therefore, when there is something on the top and on the bottom of $\mu$, the contributions vanish. This shows that the coderivation $D$ is equal to $(\theta \otimes id_{\qPs} - id_{\qPs} \otimes \theta) \circ \Delta_{(1, 1)}$.
$\cqfd$
\end{pf}

\begin{lem}\label{curvedpartprop}
Let $\Po = \F(V)/(R)$ be an inhomogeneous quadratic properad. Condition (II) implies that:
\begin{itemize}
\item The coderivation $\tilde{d}$ on $\F^{c}(sV)$ restricts to a coderivation $d_{\Poa}$ of degree $-1$ on the subcoproperad $\q\Po^{\ash} = \C(sV,\, s^{2}\q R)$;
\item The coderivation $d_{\Poa}$ satisfies $d_{\Poa}^{2} = (\theta \otimes id_{\qPs} - id_{\qPs} \otimes \theta) \circ \Delta_{(1,1)}$;
\item The coderivation $d_{\Poa}$ satisfies $\theta \circ d_{\Poa} = 0$.
\end{itemize}
\end{lem}

\begin{pf}
Condition (II) implies in particular
$$\{V\otimes R + R\otimes V\} \cap \{I \oplus V \oplus V^{\otimes 2}\} \subset R.$$
Projecting on each direct summand, we can rewrite this inclusion as the system of equations
\begin{enumerate}
\item $(s^{-1}\varphi_{1} \otimes id_{sV} + id_{sV} \otimes s^{-1}\varphi_{1}) \circ \Delta_{(1,1)}({\qPs}^{(3)}) \subset {\qPs}^{(2)}$ (projection on $V^{\otimes 2}$);
\item $\big(s^{-1}\varphi_{1} \circ (s^{-1}\varphi_{1} \otimes id_{sV} + id_{sV} \otimes s^{-1}\varphi_{1}) - (s^{-2}\varphi_{0} \otimes id_{sV} - id_{sV} \otimes s^{-2}\varphi_{0})\big) \circ {\Delta_{(1,1)}}_{|{\qPs}^{(3)}} = 0$ (projection on $V$);
\item $s^{-2}\varphi_{0} \circ (s^{-1}\varphi_{1} \otimes id_{sV} + id_{sV} \otimes s^{-1}\varphi_{1}) \circ {\Delta_{(1,1)}}_{|{\qPs}^{(3)}} = 0$ (projection on $I$).
\end{enumerate}
By the universal property which defines $\qPs = \C(sV,\, s^{2}\q R)$, it is enough to check that $\tilde{d}({\qPs}^{(3)}) \subset {\qPs}^{(2)}$ to restrict $\tilde{d}$ to a coderivation of degree $-1$ on $\qPs$, this is exactly the meaning of equation $(1)$. The equation $(2)$ corresponds to the second point of the lemma restricted to ${\qPs}^{(3)}$. The equality extends to $\qPs$ since ${d_{\Poa}}^{2} = \frac{1}{2} [d_{\Poa},\, d_{\Poa}]$ and $(\theta \otimes id_{\qPs} - id_{\qPs} \otimes \theta) \circ \Delta_{(1,1)}$ are coderivations (see Lemma \ref{Dcoder}). The equation $(3)$ corresponds to the third point of the lemma since $\theta$ is zero outside of ${\qPs}^{(2)}$.
$\cqfd$
\end{pf}

\subsection{Koszul dual coproperad} \label{defKoszuldual}
Let $\Po$ be an inhomogeneous quadratic properad with a quadratic, linear and constant presentation $\Po = \F(V)/(R)$ (such that Conditions (I) and (II) hold). The \emph{Koszul dual coproperad} of $\Po$ is the weight graded curved coproperad
$$\Po^{\ash} := (\qPs,\, d_{\Poa},\, \theta).$$

\subsection{Koszul properad}\label{Koszuldef}
A properad is called a \emph{Koszul properad} if it admits an inhomogeneous quadratic presentation $\Po = \F(V)/(R)$ such that Conditions (I) and (II) hold and such that its associated quadratic properad $\q\Po := \F(V)/(\q R)$ is Koszul in the classical sense.\\

Since the underlying $\So$-bimodule of $\Poa$ is $I \oplus sV \oplus s^{2}\q R \oplus \cdots$, we define the map of coproperads $g_{\kappa} : \Poa \mono \F^{c}(sV) \mono \B \Po$. This map is in fact a morphism of curved coproperads, so by Lemma \ref{universaltwmorph}, there is a curved twisting morphism $\kappa : \Poa \mono \B \Po \xrightarrow{\pi} \Po$. It is explicitly equal to $\Poa \epi sV \xrightarrow{s^{-1}} V \mono \Po$. By Theorem \ref{barcobaradjunction}, we also obtain the map of dg properads $\Omega \Poa \epi \Omega \B \Po \rightarrow \Po$.

\begin{thm}\label{resolution}
Let $\Po$ be a Koszul properad. The cobar construction on the Koszul dual curved coproperad $\Po^{\ash}$ is a cofibrant resolution of $\Po$:
$$\Omega \Po^{\ash} \qiso \Po.$$
\end{thm}

\begin{pf}
We work in the model category defined in the Appendix A of \cite{MerkulovVallette2}. Since we are working in the non-negatively graded case and $\Omega \Poa$ is quasi-free, the remark after Corollary 40 gives that $\Omega \Poa$ is cofibrant.

Let $\C := s^{-1}\overline{\qPs}$ be the desuspension of the augmentation coideal of the coproperad $\qPs$. So, the underlying $\So$-bimodule of $\Omega \Poa$ is $\F(\C)$. Let us consider the new ``homological'' degree induced by the weight of elements of $\qPs$, given by the weight in $\F^{c}(V)$, minus $1$. As in the proof of the Appendix $A$ of \cite{Ga-CaToVa}, Theorem $30$, we call this grading the \emph{syzygy degree}. Therefore, the syzygy degree of an element in $\F(\C)$ is given by the sum of the weight of the elements which label its vertices minus the numbers of vertices. Since the weight of an element in $\C$ is greater than $1$, the syzygy degree on $\F(\C)$ is non-negative.\\
The term $d_{0}$, induced by $\theta$, the term $d_{1}$, induced by $d_{\Poa}$ and the term $d_{2}$, induced by the infinitesimal decomposition map on $\C$, lower the syzygy degree by $1$. Hence, we get a well-defined non-negatively graded chain complex.\\
We endow $\Omega \Poa = \F(\C)$ with a filtration given by the total weight, that is the weight of an element in $\F(\C)$ is the sum of the weight of the elements which label the vertices. We have
$$d_{0} : F_{p} \rightarrow F_{p-2} \textrm{ and } d_{1} : F_{p} \rightarrow F_{p-1} \textrm{ and } d_{2} : F_{p} \rightarrow F_{p}.$$
This filtration is exhaustive and bounded below so we can apply the classical theorem of convergence of spectral sequences (cf. Theorem $5.5.1$ of \cite{Weibel}) to obtain that
$$E_{p,\, q}^{\bullet} \Rightarrow \mathrm{H}_{p+q}(\Omega \Poa).$$
The filtration $F_{p}$ induces a filtration $F_{p}$ on the homology of $\Omega \Poa$ such that
$$E_{p,\, q}^{\infty} \cong F_{p}(\mathrm{H}_{p+q}(\Omega \Poa))/F_{p-1}(\mathrm{H}_{p+q}(\Omega \Poa)) =: gr^{(p)}(\mathrm{H}_{p+q}(\Omega \Poa)).$$
Moreover, we have $E_{p,\, q}^{0} = F_{p}(\F(\C)_{p+q})/F_{p-1}(\F(\C)_{p+q}) = \F(\C)_{p+q}^{(p)}$, that is the elements of syzygy degree equal to $p+q$ and of weight $p$. The differential $d^{0}$ on the first term of the spectral sequence is given by $d_{2}$. Hence, since $\q\Po$ is Koszul and concentrated in syzygy degree $0$, we have $E_{p,\, q}^{1} = \q\Po_{p+q}^{(p)}$, concentrated in the line $p+q = 0$ and the spectral sequence collapses at rank $1$. We have
$$\left\{
\begin{array}{lcccccl}
E_{p,\, \textrm{-}p}^{1} & = & \q\Po^{(p)} & \cong & E_{p,\, \textrm{-}p}^{\infty} & \cong & gr^{(p)}(\mathrm{H}_{0}(\Omega \Poa))\\
E_{p,\, q}^{1} & = & 0 & = & E_{p,\, q}^{\infty} & \cong & gr^{(p)}(\mathrm{H}_{p+q}(\Omega \Poa)) \textrm{ when } p+q \neq 0.
\end{array}
\right.$$
For the syzygy degree, we have
$$\mathrm{H}_{0}(\Omega \Poa) \cong \F(V)/Im(d_{0}+d_{1}-d_{2}) \cong \Po.$$
So, the quotient $gr^{(p)}(\mathrm{H}_{0}(\Omega \Poa))$ is equal to $gr^{(p)}(\Po)$. Finally, the morphism $\Omega \Poa \qiso \Po$ is a quasi-isomorphism.
$\cqfd$
\end{pf}

\begin{thm}[Poincaré-Birkhoff-Witt theorem]\label{PBW}
When $\Po$ is a Koszul properad, the natural epimorphism of properads $\q\Po \epi gr\Po$ is an isomorphism of bigraded properads, with respect to the weight grading and the homological degree. Therefore, the following $\So$-bimodules, graded by the homological degree, are isomorphic
$$\Po \cong gr(\Po) \cong \q\Po.$$
\end{thm}

\begin{pf}
It is a direct corollary of the previous proof.
$\cqfd$
\end{pf}

To show Condition (II), that is $(R) \cap \{I \oplus V \oplus \F(V)^{(2)}\} = R$, can be difficult. The following proposition shows that we do not have to compute the full $(R)$ but only the part $\{R\otimes V + V \otimes R\}$.

\begin{prop}
A properad $\Po$ is Koszul if and only if it admits a presentation $\Po = \F(V)/(R)$ such that $R \subset I \oplus V \oplus \F(V)^{(2)}$ satisfying the following conditions
\begin{enumerate}
\item[(I)] $R \cap \{I \oplus V\} = \{ 0\}$;
\item[(II')] $\{R\otimes V + V \otimes R\} \cap \{V \oplus \F(V)^{(2)}\} \subset R$;
\item[(III)] the associated quadratic properad $\q\Po := \F(V)/(\q R)$ is Koszul in the classical sense.
\end{enumerate}
\end{prop}

\begin{pf}
Definition \ref{Koszuldef} always implies conditions (I), (II') and (III). First, we have to remark that the property (II') instead of (II) is enough to show Lemma \ref{curvedpartprop} and to define $\Poa$. Moreover, Theorem \ref{resolution} and Theorem \ref{PBW} are still true. Then we can apply the Poincaré-Birkhoff-Witt Theorem which gives in weight $2$ that $\q R = \q ((R) \cap \{I \oplus V \oplus \F(V)^{(2)}\})$. This last equality is equivalent to (II) under the condition (I).
$\cqfd$
\end{pf}

\section{Resolution of algebras}\label{resofalg}

We now give a resolution of semi-augmented dg properad $(\Po,\, d_{\Po},\, \varepsilon)$ as a $\Po$-bimodule this time. In the operadic case, this provides functorial cofibrant resolutions for $\Po$-algebras. We use such resolutions to define a cohomology theory associated to unital associative algebras in the next section.

\subsection{Resolutions of properads as bimodule}

We generalize the resolution given by the \emph{bar construction with coefficients} to properads (non-necessarily augmented). Moreover, for an inhomogeneous properad which is Koszul, we get a smaller resolution of it called the \emph{Koszul complex}.

\subsubsection{\bf Dg composite product}

Let $(M,\, d_{M})$ and $(N,\, d_{N})$ be two dg $\So$-bimodules. Recall from \cite{MerkulovVallette2} the differential on the monoidal product $\boxtimes$ of two $\So$-bimodules. Let $id_{M} \boxtimes' d_{N} : M \boxtimes N \rightarrow M \boxtimes N$ be the morphism of $\So$-bimodules defined by
$$(id_{M} \boxtimes' d_{N}) (\rho (m_{1},\, \ldots ,\, m_{b}) \sigma (n_{1},\, \ldots ,\, n_{a}) \omega) := \sum_{j=1}^{a} \pm \rho (m_{1},\, \ldots ,\, m_{b}) \sigma (n_{1},\, \ldots ,\, d_{N}(n_{j}),\, \ldots ,\, n_{a}) \omega$$
and let $d_{M} \bp \boxtimes id_{N} : M \boxtimes N \rightarrow M \boxtimes N$ be the morphism of $\So$-bimodules defined by
$$(d_{M} \bp \boxtimes id_{N}) (\rho (m_{1},\, \ldots ,\, m_{b}) \sigma (n_{1},\, \ldots ,\, n_{a}) \omega) := \sum_{i=1}^{b} \pm \rho (m_{1},\, \ldots ,\, d_{M}(m_{i}),\, \ldots ,\, m_{b}) \sigma (n_{1},\, \ldots ,\, n_{a}) \omega.$$
This gives a differential on $M \boxtimes N$ by $d_{M \boxtimes N} := d_{M} \bp \boxtimes id_{N} + id_{M} \boxtimes' d_{N}$.

\subsubsection{\bf Twisted composite product}

In this section, we study the free dg $\Po$-bimodules over a curved coproperad $(\C,\, d_{\C},\, \theta)$. To any map $\alpha : \C \rightarrow \Po$ of degree $-d$, we associate the unique derivation (see Section 3.2 of \cite{MerkulovVallette2} for precise definitions) of left $\Po$-module $d_{\alpha}^{l} : \Po \boxtimes \C \rightarrow \Po \boxtimes \C$ of degree $-d$ which extends the map
$$\C \xrightarrow{\Delta_{(1,\, 1)}} \C \boxtimes_{(1,\, 1)} \C \xrightarrow{\alpha \otimes id_{\C}} \Po \boxtimes_{(1,\, 1)} \C.$$
By symmetry, we define also the derivation of right $\Po$-modules $d_{\alpha}^{r} : \C \boxtimes \Po \rightarrow \C \boxtimes \Po$ of degree $-d$. We endow the free $\Po$-bimodule $\Po \boxtimes \C \boxtimes \Po$ with the following derivation of $\Po$-bimodules:
$$d_{\alpha} := d_{\Po \boxtimes \C \boxtimes \Po} - d_{\alpha}^{l} \bp \boxtimes id_{\Po} + id_{\Po} \boxtimes' d_{\alpha}^{r},$$
where $d_{\Po \boxtimes \C \boxtimes \Po} := d_{\Po} \bp \boxtimes id_{\C \boxtimes \Po} + (id_{\Po} \boxtimes' d_{\C}) \bp \boxtimes id_{\Po} + id_{\Po \boxtimes \C} \boxtimes' d_{\Po}$ with $(id_{\Po} \boxtimes' d_{\C}) \bp \boxtimes id_{\Po} = id_{\Po} \boxtimes' (d_{\C} \bp \boxtimes id_{\Po})$ by associativity of the composite product.

\begin{lem}
On the $\Po$-bimodule $\Po \boxtimes \C \boxtimes \Po$, the derivation $d_{\alpha}$ satisfies
$${d_{\alpha}}^{2} = -d^{l}_{\partial(\alpha) + \alpha \star \alpha - \Theta} \bp \boxtimes id_{\Po} + id_{\Po} \boxtimes' d^{r}_{\partial(\alpha) + \alpha \star \alpha - \Theta}.$$
Thus, when $\alpha \in \mathrm{Tw}(\C,\, \Po)$, we have ${d_{\alpha}}^{2} = 0$ and the derivation $d_{\alpha}$ defines a differential on the chain complex
$$\Po \boxtimes_{\alpha} \C \boxtimes_{\alpha} \Po := (\Po \boxtimes \C \boxtimes \Po, d_{\alpha} = d_{\Po \boxtimes \C \boxtimes \Po} - d_{\alpha}^{l} \bp \boxtimes id_{\Po} + id_{\Po} \boxtimes' d_{\alpha}^{r}).$$
\end{lem}

\begin{pf}
We do the computation for $d_{\Po} = 0$, the general case follows immediately. We have
$$\begin{array}{lcl}
{d_{\alpha}}^{2} & = & ((id_{\Po} \boxtimes' d_{\C}) \bp \boxtimes id_{\Po} - d_{\alpha}^{l} \bp \boxtimes id_{\Po} + id_{\Po} \boxtimes' d_{\alpha}^{r})^{2}\\
& = & (id_{\Po} \boxtimes' d_{\C}^{2}) \bp \boxtimes id_{\Po} + (d_{\alpha}^{l})^{2} \bp \boxtimes id_{\Po} + id_{\Po} \boxtimes' (d_{\alpha}^{r})^{2}\\
& & - ((id_{\Po} \boxtimes' d_{\C}) \circ d_{\alpha}^{l} + d_{\alpha}^{l} \circ (id_{\Po} \boxtimes' d_{\C})) \bp \boxtimes id_{\Po} + id_{\Po} \boxtimes' ((d_{\C} \bp \boxtimes id_{\Po}) \circ d_{\alpha}^{r} + d_{\alpha}^{r} \circ (d_{\C} \bp \boxtimes id_{\Po}))\\
& & - (d_{\alpha}^{l} \bp \boxtimes id_{\Po}) \circ (id_{\Po} \boxtimes' d_{\alpha}^{r}) - (id_{\Po} \boxtimes' d_{\alpha}^{r}) \circ (d_{\alpha}^{l} \bp \boxtimes id_{\Po}).
\end{array}$$
Since ${d_{\C}}^{2} = (\theta \otimes id_{\C} - id_{\C} \otimes \theta) \circ \Delta_{(1,\, 1)}$, we have $(id_{\Po} \boxtimes' {d_{\C}}^{2}) \bp \boxtimes id_{\Po} = d^{l}_{\varTheta} \bp \boxtimes id_{\Po} - id_{\Po} \boxtimes' d^{r}_{\varTheta}$. Moreover, we have $(d_{\alpha}^{l})^{2} = -d^{l}_{\alpha \star \alpha}$ and $(d_{\alpha}^{r})^{2} = d^{r}_{\alpha \star \alpha}$ since $\alpha$ has degree $-1$. Then $(id_{\Po} \boxtimes' d_{\C}) \circ d_{\alpha}^{l} + d_{\alpha}^{l} \circ (id_{\Po} \boxtimes' d_{\C}) = d^{l}_{\alpha \circ d_{\C}}$ and $(d_{\C} \bp \boxtimes id_{\Po}) \circ d_{\alpha}^{r} + d_{\alpha}^{r} \circ (d_{\C} \bp \boxtimes id_{\Po}) = d^{r}_{\alpha \circ d_{\C}}$ since $d_{\C}$ is a coderivation. Finally, $(d_{\alpha}^{l} \bp \boxtimes id_{\Po}) \circ (id_{\Po} \boxtimes' d_{\alpha}^{r}) + (id_{\Po} \boxtimes' d_{\alpha}^{r}) \circ (d_{\alpha}^{l} \bp \boxtimes id_{\Po}) = 0$ since $\alpha$ has degree $-1$. This gives the result.
$\cqfd$
\end{pf}

\subsubsection{\bf Koszul morphism}

A curved twisting morphism $\alpha : (\C,\, d_{\C},\, \theta) \rightarrow (\Po,\, d_{\Po},\, \varepsilon)$ is called a \emph{Koszul morphism} when the map $\xi$ defined by $\Po \boxtimes_{\alpha} \C \boxtimes_{\alpha} \Po \epi \Po \boxtimes I \boxtimes \Po \cong \Po \boxtimes \Po \xrightarrow{\gamma} \Po$ is a resolution of $\Po$, that is
$$\xi : \Po \boxtimes_{\alpha} \C \boxtimes_{\alpha} \Po \qiso \Po.$$

\begin{prop}\label{augbarres}
Let $\Po$ be a wfdg semi-augmented properad. The curved twisting morphism $\pi : \B \Po \rightarrow \Po$ is a curved Koszul morphism, that is the twisted composite product $\Po \boxtimes_{\pi} \B \Po \boxtimes_{\pi} \Po$ is a resolution of the properad $\Po$ called the \emph{augmented bar resolution}
$$\xi : \Po \boxtimes_{\pi} \B \Po \boxtimes_{\pi} \Po \qiso \Po.$$
\end{prop}

\begin{pf}
The method is the same as in the proof of Theorem \ref{barcobarresprop}. The weight filtration on $\Po$ induces a filtration on $\B \Po$ given by the total weight. This gives a filtration $F_{p}$ by the weight on $\Po \boxtimes_{\pi} \B \Po \boxtimes_{\pi} \Po$ and a filtration $F'_{p}$ by the weight on $\Po$. These filtrations are filtrations of chain complexes since the differentials might only decrease the weight. The filtrations are exhaustive and bounded below and the map $\xi$ preserves the filtration. We apply the classical theorem of convergence of spectral sequences (cf. Theorem 5.5.1 of \cite{Weibel}) to obtain
$$\left\{ \begin{array}{l}
E_{p,q}^{\bullet} \Rightarrow \mathrm{H}_{p+q}(\Po \boxtimes_{\pi} \B \Po \boxtimes_{\pi} \Po)\\
{E'}_{p,q}^{\bullet} \Rightarrow \mathrm{H}_{p+q}(\Po).
\end{array} \right.$$
The study of the differential show that $E_{\bullet, \bullet}^{0} = gr \Po \boxtimes_{\pi} \B (gr\Po) \boxtimes_{\pi} gr \Po$. Since $gr \Po$ is an augmented properad, we can apply Theorem 4.17 of \cite{Vallette2} to $gr \Po$ with $R = gr \Po$ to get that $E^{1}_{p,q} = \mathrm{H}_{p+q}(gr^{(p)} \Po) = {E'}_{p,q}^{1}$. Then $E_{p,q}^{r}$ and ${E'}_{p,q}^{r}$ coincide for $r \geq 1$ and $\xi$ induces an isomorphism between $E_{p,q}^{\infty}$ and ${E'}_{p,q}^{\infty} \cong gr^{(p)}\mathrm{H}_{p+q}(\Po)$. This concludes the proof.
$\cqfd$
\end{pf}

Let $\Po$ be an inhomogeneous properad, $\Poa$ its Koszul dual cooperad and $\kappa : \Poa \rightarrow \Po$ the associated curved twisting morphism. The chain complex $\Po \boxtimes_{\kappa} \Poa \boxtimes_{\kappa} \Po$ is called the \emph{total Koszul complex}.

\begin{prop}\label{koszulres}
Let $\Po$ be an inhomogeneous properad and $\Poa$ be its Koszul dual coproperad. When $\Po$ is Koszul, the curved twisting morphism $\kappa : \Poa \rightarrow \Po$ is a curved Koszul morphism, that is the total Koszul complex $\Po \boxtimes_{\kappa} \Poa \boxtimes_{\kappa} \Po$ is a resolution of the properad $\Po$
$$\xi : \Po \boxtimes_{\kappa} \Poa \boxtimes_{\kappa} \Po \qiso \Po.$$
\end{prop}

\begin{pf}
The proof is similar to the proof of Proposition \ref{augbarres}. The differences are the following. Since $\Po$ is Koszul, the Poincaré-Birkhoff-Witt Theorem \ref{PBW} gives $E_{\bullet, \bullet}^{0} = gr\Po \boxtimes_{\kappa} \q\Poa \boxtimes_{\kappa} gr\Po \cong \q\Po \boxtimes_{\kappa} \q\Poa \boxtimes_{\kappa} \q\Po$. So $E_{p, q}^{1} = \mathrm{H}_{p+q}(\q\Po^{(p)}) \cong \mathrm{H}_{p+q}(gr^{(p)}\Po)$ (by Theorem 7.8 and by Theorem 5.4 of \cite{Vallette2} with $L = \q\Po \boxtimes_{\kappa} \q\Poa \boxtimes_{\kappa} \q\Po$, $L' = \q\Po$, $\Po' = \q\Po$, $M = \q\Po \boxtimes_{\kappa} \q\Poa$ and $M' = I$ and by the Poincaré-Birkhoff-Witt Theorem).
$\cqfd$
\end{pf}

\subsection{Resolution of algebras}

From now on, we consider only operads and cooperads since there is in general no notion of free algebra over a properad. In this section, we use the resolutions of $\Po$ as a $\Po$-bimodule of the previous section to provide functorial resolutions for algebras over $\Po$ as, for example, for unital associative algebras (see Section \ref{resandtransfer}).

\subsubsection{\bf Coalgebra over a curved cooperad}

Let $(\C,\, d_{\C},\, \theta)$ be a curved cooperad. A \emph{$(\C,\, d_{\C},\, \theta)$-coalgebra} $(C,\, \Delta_{C}, d_{C})$ is a $\C$-coalgebra such that
$$\left\{ \begin{array}{l}
\Delta_{C} \circ d_{C} = (d_{\C} \circ id_{C} + id_{\C} \circ' d_{C}) \circ \Delta_{C}\\
{d_{C}}^{2} = (\theta \circ id_{C}) \circ \Delta_{C},
\end{array} \right.$$
where the $\circ'$ corresponds to the previous $\boxtimes'$ but for operads.

A \emph{morphism of $(\C,\, d_{\C},\, \theta)$-coalgebras} $f : (C,\, \Delta_{C}, d_{C}) \rightarrow (C',\, \Delta_{C'}, d_{C'})$ is a morphism $f : C \rightarrow C'$ of $\C$-coalgebras which commutes with the predifferentials $d_{C}$ and $d_{C'}$.

\subsubsection{\bf Relative composition product}

Let $(\Po,\, d_{\Po},\, \varepsilon)$ be a sdg operad. A \emph{right $\Po$-module} $(\Lc,\, \rho)$ is an $\So$-module endowed with a map $\rho : \Lc \circ \Po \rightarrow \Lc$ compatible with the product and the unit of the operad. We define similarly the notion of \emph{left $\Po$-module}. We define the relative composite product $\Lc \circ_{\Po} \R$ of a right $\Po$-module $(\Lc,\, \rho)$ and a left $\Po$-module $(\R,\, \lambda)$ by the coequalizer diagram
$$\xymatrix{\mathcal{L} \circ \Po \circ \R  \ar@<0.5ex>[r]^{\hspace{0.2cm} \rho \circ id_{\R}} \ar@<-0.5ex>[r]_{\hspace{0.2cm} id_{\mathcal{L}} \circ \lambda} & \mathcal{L} \circ \R \ar@{->>}[r] & \mathcal{L} \circ_{\Po} \R}.$$
These definitions extend to the dg setting.

\subsubsection{\bf Bar construction of $\Po$-algebras}\label{barconstruction}

To any curved twisting morphism $\alpha : \C \rightarrow \Po$ from a curved cooperad $(\C,\, d_{\C},\, \theta)$ to an operad $(\Po,\, d_{\Po},\, \varepsilon)$, we associate a functor
$$\B_{\alpha} : \textsf{dg }(\Po,\, d_{\Po},\, \varepsilon) \textsf{-algebras} \rightarrow \textsf{dg }(\C,\, d_{\C},\, \theta) \textsf{-coalgebras}.$$
For a $\Po$-algebra $(A,\, \gamma_{A})$, we define on $\C(A) = (\C \circ \Po) \circ_{\Po} A$ the maps
$$\left\{ \begin{array}{l}
d_{1} : \C(A) \xrightarrow{d_{\C} \circ id_{A} + id_{\C} \circ' d_{A}} \C (A)\\
d_{2} := d_{\alpha}^{r} \circ_{\Po} id_{A} : \C(A) \xrightarrow{\Delta_{(1)} \circ id_{A}} \C \circ_{(1)} \C (A) \xrightarrow{(id_{\C} \otimes \alpha) \circ id_{A}} \C \circ \Po (A) \xrightarrow{id_{\C} \circ \gamma_{A}} \C (A),
\end{array} \right.$$
where $(\C \circ \Po) \circ A \xrightarrow{d_{\alpha}^{r} \circ id_{A}} (\C \circ \Po) \circ A \epi \C (A)$ factors through $\C (A)$ to give $d_{\alpha}^{r} \circ_{\Po} id_{A}$ since $\gamma_{A}$ is a dg map. (Here, $\Delta_{(1)}$ corresponds to the infinitesimal decomposition map $\Delta_{(1,1)}$ and $\C \circ_{(1)} \C$ corresponds to $\C \boxtimes_{(1,1)} \C$ when we restrict to cooperads.)

\begin{lem}\label{forcodiff}
Since $\alpha$ is a curved twisting morphism, we have
$$(d_{1} + d_{2})^{2} = (\theta \circ id_{\C(A)}) \circ \Delta_{\C(A)}.$$
\end{lem}

\begin{pf}
We compute
$$\left\{ \begin{array}{cl}
{d_{1}}^{2} & = {d_{\C}}^{2} \circ id_{A} = ((\theta \otimes id_{\C} - id_{\C} \otimes \theta) \circ \Delta_{(1)}) \circ id_{A}\\
& = (\theta \circ id_{\C(A)}) \circ \Delta_{\C(A)} - d_{\Theta}^{r} \circ_{\Po} id_{A}\\
{d_{2}}^{2} & = d_{\alpha \star \alpha}^{r} \circ_{\Po} id_{A}\\
d_{1} d_{2} + d_{2} d_{1} & = d_{\partial(\alpha)}^{r} \circ_{\Po} id_{A}.
\end{array} \right.$$
Thus $(d_{1} + d_{2})^{2} = d_{\partial(\alpha) + \alpha \star \alpha - \varTheta}^{r} \circ_{\Po} id_{A} + (\theta \circ id_{\C(A)}) \circ \Delta_{\C (A)}$ and we get the result since $\alpha$ is a curved twisting morphism.
$\cqfd$
\end{pf}

The \emph{bar construction on $A$} is the $(\C,\, d_{\C},\, \theta)$-coalgebra $\B_{\alpha}A := (\C(A),\, d := d_{1} + d_{2})$.

\subsubsection{\bf Cobar construction of a $\C$-coalgebra}\label{cobarconstruction}

Similarly to the previous section, to any curved twisting morphism $\alpha : (\C,\, d_{\C},\, \theta) \rightarrow (\Po,\, d_{\Po},\, \varepsilon)$, we associate a functor
$$\Omega_{\alpha} : \textsf{dg }(\C,\, d_{\C},\, \theta) \textsf{-coalgebras} \rightarrow \textsf{dg }(\Po,\, d_{\Po},\, \varepsilon) \textsf{-algebras}.$$
For any $(\C,\, d_{\C},\, \theta)$-coalgebra $(C,\, d_{C},\, \Delta_{C})$, we define on $\Po (C)$ the maps
$$\left\{ \begin{array}{l}
d_{1} : \Po(C) \xrightarrow{d_{\Po} \circ id_{C} + id_{\Po} \circ' d_{C}} \Po (C)\\
d_{2} : \Po(C) \xrightarrow{id_{\Po} \circ' \Delta_{C}} \Po \circ_{(1)} \C (C) \xrightarrow{(id_{\Po} \otimes \alpha) \circ id_{C}} \Po \circ \Po (C) \xrightarrow{\gamma \circ id_{C}} \Po (C).
\end{array} \right.$$

\begin{lem}
Since $\alpha$ is a curved twisting morphism, we have
$$(d_{1} - d_{2})^{2} = 0.$$
\end{lem}

\begin{pf}
We compute
$$\left\{ \begin{array}{cl}
{d_{1}}^{2} & = id_{\Po} \circ' {d_{C}}^{2} = id_{\Po} \circ' ((\theta \circ id_{C}) \circ \Delta_{C})\\
{d_{2}}^{2} & = - (\gamma \circ id_{C}) \circ (id_{\Po} \circ (\alpha \star \alpha) \circ id_{C}) \circ (id_{\Po} \circ' \Delta_{C})\\
- d_{1} d_{2} - d_{2} d_{1} & = - (\gamma \circ id_{C}) \circ (id_{\Po} \circ \partial(\alpha) \circ id_{C}) \circ (id_{\Po} \circ' \Delta_{C}).
\end{array} \right.$$
Thus $(d_{1} - d_{2})^{2} = - (\gamma \circ id_{C}) \circ (id_{\Po} \circ (\partial(\alpha) + \alpha \star \alpha - \varTheta) \circ id_{C}) \circ (id_{\Po} \circ' \Delta_{C}) = 0$ since $\alpha$ is a curved twisting morphism.
$\cqfd$
\end{pf}

The \emph{cobar construction on $C$} is the dg $\Po$-algebra $\Omega_{\alpha} C := (\Po (C),\, d_{\Omega_{\alpha}C} := d_{1} - d_{2})$.

\subsubsection{\bf The bar-cobar resolution}

The bar-cobar construction on a $\Po$-algebra provides a functorial cofibrant resolution of any $\Po$-algebra when the curved twisting morphism $\alpha$ is Koszul.

\begin{prop}\label{barcobarres}
Let $\alpha : (\C,\, d_{\C},\, \theta) \rightarrow (\Po,\, d_{\Po},\, \varepsilon)$ be a curved Koszul morphism between a curved cooperad $(\C,\, d_{\C},\, \theta)$ and a sdg operad $(\Po,\, d_{\Po},\, \varepsilon)$. Then the bar-cobar resolution $\Omega_{\alpha} \B_{\alpha} A$ is a cofibrant construction of the $\Po$-algebra $A$, that is
$$\Omega_{\alpha} \B_{\alpha} A = \Po \circ_{\alpha} \C \circ_{\alpha} A \qiso A.$$
\end{prop}

\begin{pf}
In the (semi-)model category structure on $\Po$-algebras defined in \cite{Fresse3}, cofibrant $\Po$-algebras are retracts of quasi-free $\Po$-algebras endowed with a good filtration (Proposition 12.3.8 in \cite{Fresse3}). This is the case here since the chain complexes are non-negatively graded.

The bar-cobar construction $\Omega_{\alpha} \B_{\alpha}A$ is isomorphic to the relative composite product $(\Po \circ_{\alpha} \C \circ_{\alpha} \Po) \circ_{\Po} A$. So it is defined by a short exact sequence
$$0 \rightarrow (\Po \circ_{\alpha} \C \circ_{\alpha} \Po) \circ \Po \circ A \rightarrow (\Po \circ_{\alpha} \C \circ_{\alpha} \Po) \circ A \rightarrow (\Po \circ_{\alpha} \C \circ_{\alpha} \Po) \circ_{\Po} A \rightarrow 0$$
which induces a long exact sequence in homology. Since $\K$ is a field of characteristic $0$, the rings $\K [\So_{n}]$ are semi-simple by Maschke's Theorem, that is every $\K[\So_{n}]$-module is projective. So K\"unneth formula and the fact that $\alpha$ is Koszul imply that $\mathrm{H}_{\bullet}((\Po \circ_{\alpha} \C \circ_{\alpha} \Po) \circ \Po \circ A) \cong \mathrm{H}_{\bullet}(\Po) \circ \mathrm{H}_{\bullet}(\Po) \circ \mathrm{H}_{\bullet}(A)$ and that $\mathrm{H}_{\bullet}((\Po \circ_{\alpha} \C \circ_{\alpha} \Po) \circ A) \cong \mathrm{H}_{\bullet}(\Po) \circ \mathrm{H}_{\bullet}(A)$. Finally, this gives
$$\mathrm{H}_{\bullet}((\Po \circ_{\alpha} \C \circ_{\alpha} \Po) \circ_{\Po} \circ A) \cong \mathrm{H}_{\bullet}(\Po) \circ_{\mathrm{H}_{\bullet}(\Po)} \mathrm{H}_{\bullet}(A) \cong \mathrm{H}_{\bullet}(A).$$
$\cqfd$
\end{pf}

\begin{thm}
Let $(\Po,\, d_{\Po},\, \varepsilon)$ be a sdg operad. The curved Koszul morphism $\pi : \B \Po \rightarrow \Po$ gives a cofibrant resolution
$$\Omega_{\pi} \B_{\pi} A = \Po \circ_{\pi} \B \Po \circ_{\pi} A \qiso A.$$
When $\Po$ is a Koszul operad, the total Koszul complex gives a smaller cofibrant resolution
$$\Omega_{\kappa} \B_{\kappa} A = \Po \circ_{\kappa} \Poa \circ_{\kappa} A \qiso A.$$
\end{thm}

\begin{pf}
It is a direct corollary of Proposition \ref{barcobarres} and Propositions \ref{augbarres} and \ref{koszulres}.
$\cqfd$
\end{pf}

\section{Homotopy and cohomology theories for unital associative algebras}\label{resandtransfer}

In this section we describe a simple  resolution of the operad which
encodes unital associative algebras, $\uAs$, obtained by the methods described in section
\ref{curvedKD}. In fact, many of the theorems in this section can be generalized in a straightforward way to any (inhomogeneous) Koszul properad. Algebras over the resolution $\uAsinf$ are called \emph{homotopy
  unital $A_{\infty}$-algebras}, or \emph{$\uAsinf$-algebras}, for
short. We use some nice properties of our
resolution to prove that $\uAsinf$-algebras may be replaced
up to equivalence
by strictly unital associative algebras. Using our explicit transfer
formulae, we show that a unital associative algebra may be transferred
to homology as a strictly
unital $A_{\infty}$-algebra  (see Definition \ref{suainf}). This gives a proof that one may always choose a minimal
model for a $\uAsinf$-algebra which is actually a strictly unital
$A_{\infty}$-algebra. In this sense, it is ``enough'' to resolve only
the associative relation of $\uAs$, obtaining the
operad $A_{\infty}$, and then adjoin a unit, giving the operad which encodes strictly unital $A_{\infty}$-algebras. As a corollary of our discussion, we provide
sufficient conditions so that: ``When trying to find
resolutions of algebraic structures with units, it is `good
enough' to resolve the structure (without its units) first, and then
append the units to that resolution.'' The notion of $\uAsinf$-algebras is exactly the notion of ``$A_{\infty}$-algebras with a homotopy unit'' of \cite{Fukaya}. Concerning the notion of $\infty$-morphism and the nice properties of $\uAsinf$-algebras, we still have to compare them with the theory presented in \cite{Fukaya}.

\subsection{Homotopy unital associative algebras}

We give a presentation for the operad encoding unital associative
algebras. This presentation is an inhomogeneous quadratic presentation
and we can apply the theory of the previous sections to compute its
Koszul dual cooperad, and hence an explicit resolution.

We use the notation $\underline{n} := \{1,\, \ldots,\, n\}$, the notation $\overline{\mu}$ stands for an element in a cooperad and the notation $\mu$ stands for an element in an operad.

\subsubsection{\bf The operad encoding unital associative algebras}
We denote by $u\As$ the operad whose representations in
the category of dg modules are precisely differential graded unital
associative algebras. We consider the following
presentation
 $$\uAs = {\small{\F \left(\vcenter{
         \xymatrix@M=0pt@R=3pt@C=3pt{&&\\
           & \ar@{{*}}[u] \ar@{-}[d]\\
           &}}, \,
       \vcenter{\xymatrix@M=0pt@R=4pt@C=4pt{\ar@{-}[dr] && \ar@{-}[dl]\\
           & \ar@{-}[d] &\\
           &&}}\right) / (\vcenter{ \xymatrix@M=0pt@R=4pt@C=4pt{
           & & & & \\ \ar@{-}[ddrr] & & \ar@{-}[dl] & & \ar@{-}[ddll]\\
           &  & &  & \\
           & & \ar@{-}[d] & &\\
           & & }} - \vcenter{
         \xymatrix@M=0pt@R=4pt@C=4pt{& & & & \\
           \ar@{-}[drdr] & &\ar@{-}[dr] & & \ar@{-}[dldl]  \\
           & & & & \\
           & &\ar@{-}[d] & & \\
           & & }}, \, \vcenter{ \xymatrix@M=0pt@R=4pt@C=4pt{
           & & & \\
           & & &\\
           & \ar@{{*}}[u] \ar@{-}[dr] & & \ar@{-}[dl] \\
           & & \ar@{-}[d] &\\
           & & }} - |, \, \vcenter{
         \xymatrix@M=0pt@R=4pt@C=4pt{& & & \\
           & & &\\
           & \ar@{-}[dr] & & \ar@{{*}}[u] \ar@{-}[dl]\\
           & &\ar@{-}[d] & \\
           & & }} - | ).}}$$

\begin{rem}
  We fix this presentation to make our computations of the Koszul
  dual, $\uAs^{\ash}$ and ultimately $\uAsinf$. Note that this presentation for
  $u\As$ is an inhomogeneous quadratic presentation (see
  \ref{inhomoquadprop} for a definition).
\end{rem}

To make the Koszul dual cooperad, $\uAs^{\ash}$ of $\uAs$ explicit, we compute
its associated quadratic operad
\begin{equation*}
\quAs = {\small{\F (\vcenter{
    \xymatrix@M=0pt@R=3pt@C=3pt{&&\\
      & \ar@{{*}}[u] \ar@{-}[d]\\
      &}}, \,
  \vcenter{\xymatrix@M=0pt@R=4pt@C=4pt{\ar@{-}[dr] && \ar@{-}[dl]\\
      & \ar@{-}[d] &\\
      &&}})/ (\vcenter{ \xymatrix@M=0pt@R=4pt@C=4pt{
      & & & & \\ \ar@{-}[ddrr] & & \ar@{-}[dl] & & \ar@{-}[ddll]\\
      &  & &  & \\
      & & \ar@{-}[d] & &\\
      & & }} - \vcenter{
    \xymatrix@M=0pt@R=4pt@C=4pt{& & & & \\
      \ar@{-}[drdr] & &\ar@{-}[dr] & & \ar@{-}[dldl]  \\
      & & & & \\
      & &\ar@{-}[d] & & \\
      & & }}, \, \vcenter{ \xymatrix@M=0pt@R=4pt@C=4pt{
      & & & & \\
      & & & &\\
      & \ar@{{*}}[u] \ar@{-}[dr] & & \ar@{-}[dl] & \\
      & & \ar@{-}[d] & &\\
      & & }}, \vcenter{
    \xymatrix@M=0pt@R=4pt@C=4pt{& & & & \\
      & & & &\\
      & \ar@{-}[dr] & & \ar@{{*}}[u] \ar@{-}[dl] & \\
      & &\ar@{-}[d] & & \\
      & & }} )
= \un \oplus \As .
}}
\end{equation*}

Let's take a moment to explain the notation on the left-hand side of
the equation above.

\begin{dei}
Let $\Po$, $\Q$ be augmented operads. Then the \emph{direct sum operad} $\Po \oplus \Q$ is defined
to be $\mathcal{F} \left( \overline{\Po}, \overline{\Q} \right) / \left( R_{\Po}, R_{\Q},
  R_{\Po \Q} \right)$, where $R_{\Po}, R_{\Q}$ are the relations in
$\Po, \Q$ respectively, and $R_{\Po \Q}$ is the collection of all
compositions of a pair of elements, one in  $\overline{\Po}$, one in $\overline{\Q}$.
\end{dei}

\begin{rem}
The direct sum operad is the product in the category of operads.
\end{rem}

\begin{prop}
If $\Po$ and $\Q$ are both quadratic augmented operads, then $\Po \oplus \Q$ is a quadratic augmented operad.
\end{prop}

\begin{pf}
For any two presented operads, $\Po = \mathcal{F} \left(V_1 \right) /
\left( R_1 \right), \Q = \mathcal{F} \left( V_2 \right) / \left( R_2
\right)$, the direct sum operad $\Po \oplus \Q$ is naturally presented
by $\mathcal{F} \left( V_1, V_2 \right) / \left( R_1, R_2, R_{V_1 V_2}
\right)$. If $(V_1, R_1)$ and $(V_2, R_2)$ are both quadratic
presentations, then so is the natural presentation for $\Po \oplus \Q$.
$\cqfd$
\end{pf}

We will make use of the identification $\quAs = \un \oplus \As$ to
compute the Koszul dual cooperad of $\quAs$ (see
\ref{Koszuldualproduct}). Before we compute the resulting cooperad,
$\quAs$, we first describe it.

Linearly, we have an isomorphism $\q\uAs^{\ash} \cong
\K \cdot \big\{ \munsc \big\}_{n \geq 1, S \subset
  \underline{n}}$. The element $\munsc \in \q\uAs^{\ash}$
represents a (co)operation with $n-|S|$ inputs: however, we draw this
operation as a corolla with $n$ leaves, and a \emph{cork} covering each of the leaves in the set $S$. For example,
$\overline{\mu}_{5}^{\{1,4\}}$ represents $\smfivectwo.$ We point out here that
the space of $n$-to-1 operations is infinite dimensional for every $n
\geq 0$. To see this, note that every $n$-to-1 corolla is an
$n$-ary operation, and by adding a corked leaf, we get a new
$n$-to-1 operation. Continuing to add corked leaves gives
infinitely many new $n$-ary operations.

Also notice that $\overline{\mu}_{n}^{\emptyset} = \vcenter{\xymatrix@M=0pt@R=4pt@C=4pt{
    & & & & \\
    \ar@{-}[drr] & \ar@{-}[dr] & \ar@{-}[d] & \ar@{-}[dl] & \ar@{-}[dll] \\
    & & \ar@{-}[d] \\
    & & }}$ for $n \geq 1$ spans the subcooperad corresponding to
$\As^{\ash}$ and $\big\{ \overline{\mu}_{1}^{\emptyset} = |,\, \overline{\mu}_{1}^{\{1\}} =
\vcenter{
  \xymatrix@M=0pt@R=3pt@C=3pt{&&\\
    & \ar@{{*}}[u] \ar@{-}[d]\\
    &}} \big\}$ spans the subcooperad corresponding to $\un^{\ash}$
(with $\overline{\mu}_{1}^{\emptyset}$ corresponding to the identity cooperation in
both cases).

Using this basis, the infinitesimal decomposition $\Delta_{(1)}$ is
given by summing over all possible ways to split the corolla into two, preserving the number of leaves
and the number and positions of the corks. Pictorially:
$$\mfivectwo \mapsto \Sigma \pm \vcenter
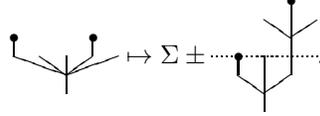
{\xymatrix@M=0pt@R=7pt@C=10pt{
   & & & & \\
   & & \ar@{-}[dr] & \ar@{-}[d] \ar@{{*}}[u] & \ar@{-}[dl] \\
   & & & \ar@{-}[d] & \\
   \ar@{.}[rrrr] & & & & \\
   & \ar@{{*}}[u] \ar@{-}[dr] &\ar@{-}[d] \ar@{-}[u] &\ar@{-}[dl] \ar@{-}[u] & \\
   & &\ar@{-}[d] & & \\
   & & & & }}.$$ For example,
$$\Delta_{(1)}\left(\ \vcenter{\xymatrix@M=0pt@R=4pt@C=6pt{
      & & \\
      \ar@{-}[dr] \ar@{{*}}[u] & \ar@{-}[d] & \ar@{-}[dl] \\
      & \ar@{-}[d] & \\
      & & }}\ \right) = \vcenter{\xymatrix@M=0pt@R=4pt@C=6pt{
    & & & &\\
    \ar@{.}[rrrr] & \ar@{-}[dr] \ar@{{*}}[u] & \ar@{-}[d] & \ar@{-}[dl] &\\
    & &\ar@{-}[d] & &\\
    & & & &}} - \vcenter{\xymatrix@M=0pt@R=4pt@C=6pt{
    & &  & &\\
    \ar@{{*}}[u] \ar@{-}[dr] & &\ar@{-}[dl] & &\\
    & \ar@{-}[d] & & &\\
    \ar@{.}[rrrr] & \ar@{-}[dr] & & \ar@{-}[dl] &\\
    & & \ar@{-}[d] & &\\
    & & & &}} - \vcenter{\xymatrix@M=0pt@R=4pt@C=6pt{
    & & \ar@{-}[dr] & & \ar@{-}[dl] \\
    & & & \ar@{-}[dd] & \\
    \ar@{.}[rrrr]& & & & \\
    & \ar@{{*}}[u] \ar@{-}[dr] & &\ar@{-}[dl] & \\
    & &\ar@{-}[d] & & \\
    & & & }}.$$

We compute the Koszul dual cooperad, $q\uAs^{\ash}$ by
the following proposition.

\begin{prop} \label{Koszuldualproduct} Let $\Po = \F(V)/(R)$ be a
  quadratic operad where $V$ is finite-dimensional. We denote by $\vcenter{
    \xymatrix@M=0pt@R=3pt@C=3pt{&&\\
      & \ar@{{*}}[u] \ar@{-}[d]\\
      &}} \oplus \Po$ the operad given by $(\vcenter{
    \xymatrix@M=0pt@R=3pt@C=3pt{&&\\
      & \ar@{{*}}[u] \ar@{-}[d]\\
      &}} \oplus \Po)(0) := (\K \cdot \vcenter{
    \xymatrix@M=0pt@R=3pt@C=3pt{&&\\
      & \ar@{{*}}[u] \ar@{-}[d]\\
      &}}) \oplus \Po(0)$ and $(\vcenter{
    \xymatrix@M=0pt@R=3pt@C=3pt{&&\\
      & \ar@{{*}}[u] \ar@{-}[d]\\
      &}} \oplus \Po)(n) := \Po(n)$ for all $n\neq 0$ and endowed with the operadic structure given by the structure on $\vcenter{
    \xymatrix@M=0pt@R=3pt@C=3pt{&&\\
      & \ar@{{*}}[u] \ar@{-}[d]\\
      &}}$ and the structure on $\Po$. The Koszul dual
  cooperad of $\vcenter{
    \xymatrix@M=0pt@R=3pt@C=3pt{&&\\
      & \ar@{{*}}[u] \ar@{-}[d]\\
      &}} \oplus \Po$ is given by the coaugmented cooperad
$$(\vcenter{
  \xymatrix@M=0pt@R=3pt@C=3pt{&&\\
    & \ar@{{*}}[u] \ar@{-}[d]\\
    &}} \oplus \Po)^{\ash} \cong \K \cdot \{ \overline{\mu}_{S}, \textrm{
  where } \overline{\mu} \in \Poa(n),\, S \subset \underline{n} \textrm{ and
} |\overline{\mu}_{S}| = |\overline{\mu}| + |S|\}.$$ The set $S$ is the set of the
positions of the ``corks'' $\vcenter{
  \xymatrix@M=0pt@R=3pt@C=3pt{&&\\
    & \ar@{{*}}[u] \ar@{-}[d]\\
    &}}$. Let $\overline{\xi} \in \Poa(n)$ such that $\Delta_{(1)}(\overline{\xi}) =
\sum (\overline{\mu};\, \underbrace{id,\, \ldots ,\, id}_{p},\, \overline{\nu},\,
\underbrace{id,\, \ldots,\, id}_{r}),$ where $\overline{\mu} \in \Poa(m)$, $\overline{\nu}
\in \Poa(q)$, $p+1+r = m$ and $p+q+r = n$. Then the infinitesimal decomposition map
on $\overline{\xi}_{S} \in (\vcenter{
  \xymatrix@M=0pt@R=3pt@C=3pt{&&\\
    & \ar@{{*}}[u] \ar@{-}[d]\\
    &}} \oplus \Po)^{\ash}$, where $S \subset \underline{n}$, is given
by
$$\Delta_{(1)}(\overline{\xi}_{S}) = \sum (-1)^{|\overline{\nu}| |S_{1}| + |S_{2}| |S''_{1}|} (\overline{\mu}_{S_{1}}; \underbrace{id,\, \ldots ,\, id}_{p-|S'_{1}|},\, \overline{\nu}_{S_{2}},\, \underbrace{id,\, \ldots,\, id}_{r-|S''_{1}|}),$$
where $\overline{\mu}_{S_{1}} \in \Poa(m-|S_{1}|)$, $\overline{\nu}_{S_{2}} \in
\Poa(q-|S_{2}|)$ and $\left\{ \begin{array}{lll} S'_{1} & \subset &
    \underline{p}\\ S_{2} & \subset & \underline{q}\\ S''_{1} &
    \subset & \{p+2,\, \ldots,\, p+1+r\} \end{array} \right.$ such
that $S = S'_{1} \sqcup (S_{2} + p) \sqcup (S''_{1} + q)$ and $S_{1} =
S'_{1} \sqcup S''_{1}$.
\end{prop}

\begin{pf}
  The operad $\vcenter{
    \xymatrix@M=0pt@R=3pt@C=3pt{&&\\
      & \ar@{{*}}[u] \ar@{-}[d]\\
      &}} \oplus \Po$ is a quadratic operad given by $\F(\vcenter{
    \xymatrix@M=0pt@R=3pt@C=3pt{&&\\
      & \ar@{{*}}[u] \ar@{-}[d]\\
      &}} \oplus V)/(R \oplus V \otimes \vcenter{
    \xymatrix@M=0pt@R=3pt@C=3pt{&&\\
      & \ar@{{*}}[u] \ar@{-}[d]\\
      &}})$ where
$$V \otimes \vcenter{
  \xymatrix@M=0pt@R=3pt@C=3pt{&&\\
    & \ar@{{*}}[u] \ar@{-}[d]\\
    &}} := \{ \mu^{\{ k\}}, \textrm{ with } \mu \in V(n) \textrm{ and
} \{k\} \subset \underline{n}\}.$$ We follow the Appendix B of
\cite{Loday2} defining $\Po^{\textrm{!}} := \F(V^{\vee})/(R^{\bot})$,
where $V^{\vee} := V^{*} \otimes (sgn)$ with the signature representation $(sgn)$ and $R^{\bot}$ is the
orthogonal space for the natural pairing $\langle -,- \rangle : V^{\vee} \otimes V
\rightarrow \K$. Since $(V \otimes \vcenter{
  \xymatrix@M=0pt@R=3pt@C=3pt{&&\\
    & \ar@{{*}}[u] \ar@{-}[d]\\
    &}})^{\bot} = \F_{(2)}(V^{\vee})$, we get
$$(\vcenter{
  \xymatrix@M=0pt@R=3pt@C=3pt{&&\\
    & \ar@{{*}}[u] \ar@{-}[d]\\
    &}} \oplus \Po)^{\mathrm{!}} = \F(\vcenter{
  \xymatrix@M=0pt@R=3pt@C=3pt{&&\\
    & \ar@{{*}}[u] \ar@{-}[d]\\
    &}}^{\vee} \oplus V^{\vee})/(R^{\bot} \cap \F_{(2)}(V^{\vee}))
\cong \{ \mu^{S}, \textrm{ where } \mu \in \Po^{\textrm{!}}(n)
\textrm{ and } S \subset \underline{n}\}$$ and the product is induced,
up to signs, by the product on $\Po^{\textrm{!}}$.

Following \cite{LodayVallette}, we have $\Poa :=
{\mathcal{S}^{-1}}^{c} \otimes_{H} (\Po^{\textrm{!}})^{*}$ where
${\mathcal{S}^{-1}}^{c}$ is the operadic desuspension. Then the Koszul dual cooperad of $\vcenter{
  \xymatrix@M=0pt@R=3pt@C=3pt{&&\\
    & \ar@{{*}}[u] \ar@{-}[d]\\
    &}} \oplus \Po$ is equal to
$$(\vcenter{
  \xymatrix@M=0pt@R=3pt@C=3pt{&&\\
    & \ar@{{*}}[u] \ar@{-}[d]\\
    &}} \oplus \Po)^{\ash} \cong \{ \overline{\mu}^{S}, \textrm{ where }
\overline{\mu} \in \Poa(n),\, S \subset \underline{n} \textrm{ and }
|\overline{\mu}^{S}| = |\overline{\mu}| + |S|\}$$ and the (infinitesimal) coproduct is given, up to
signs, by the (infinitesimal) coproduct of $\Poa$. To compute the signs, we recall
that the corks $\vcenter{
  \xymatrix@M=0pt@R=3pt@C=3pt{&&\\
    & \ar@{{*}}[u] \ar@{-}[d]\\
    &}}$ have degree $-1$ and we apply the Koszul rule. The sign
$(-1)^{|\overline{\nu}| |S_{1}|}$ in the formula of the proposition comes
from the fact that $\overline{\nu}$ passes through the corks indexed by
$S_{1}$ and the sign $(-1)^{|S_{2}| |S''_{1}|}$ comes from the fact
that the corks indexed by $S_{2}$ pass through the corks indexed by
$S''_{2}$.  $\cqfd$
\end{pf}

\begin{cor}\label{quAs}
  The Koszul dual cooperad associated to $\q\uAs$ is equal to
$$\quAs = (\vcenter{
  \xymatrix@M=0pt@R=3pt@C=3pt{&&\\
    & \ar@{{*}}[u] \ar@{-}[d]\\
    &}} \oplus \As)^{\ash} \cong \K \cdot \{ \munsc \},$$
$\textrm{ where } \overline{\mu}_{n} \in \As^{\ash}(n),\, S \subset
\underline{n} \textrm{, so } \munsc \in \uAs^{\ash}(n-|S|) \textrm{
  and } | \munsc | = n-1 + |S|.$ The infinitesimal decomposition map is given by
$$\Delta_{(1)} ( \munsc ) = \sum_{\substack{p+q+r=n\\ p+1+r=m}} (-1)^{(q+1)(r+|S_{1}|)+|S_{2}||S_{1}''|} (\overline{\mu}_{m}^{S_{1}}; \underbrace{\textrm{id, \ldots, id}}_{p-|S'_{1}|},\, \overline{\mu}_{q}^{S_{2}},\, \underbrace{\textrm{id, \ldots, id}}_{r-|S''_{1}|}),$$
where $\left\{ \begin{array}{lll} S'_{1}
    & \subset & \underline{p}\\ S_{2} & \subset & \underline{q}\\
    S''_{1} & \subset & \{p+2,\, \ldots,\, p+1+r\} \end{array}
\right.$ such that $S = S'_{1} \sqcup (S_{2} + p) \sqcup (S''_{1} +
q)$ and $S_{1} = S'_{1} \sqcup S''_{1}$.

Moreover, the coproduct is given by
$$\Delta( \munsc ) = \sum_{i_{1}-|T_{1}|+ \cdots +i_{m-|T|}-|T_{m-|T|}|=n-|S|} (-1)^{\epsilon} (\overline{\mu}_{m}^{T}; \overline{\mu}_{i_{1}}^{T_{1}},\, \ldots ,\, \overline{\mu}_{i_{m-|T|}}^{T_{m-|T|}}),$$
where $\left\{ \begin{array}{lll} T & \subset & \underline{m}\\ T_{j} & \subset & \underline{i_{j}} \end{array} \right.$ such that $T = R_{0} \sqcup \ldots \sqcup R_{m-|T|}$ and
$$\begin{array}{lll} S & = & R_{0} \sqcup (T_{1}+|R_{0}|) \sqcup (R_{1}+i_{1}) \sqcup \ldots \sqcup \\ && (T_{m-|T|}+|R_{0}|+\cdots +|R_{m-|T|-1}|+i_{1}+\cdots +i_{m-|T|-1}) \sqcup (R_{m-|T|}+i_{1}+\cdots +i_{m-|T|}) \end{array}$$
and where
$$\epsilon := |T|(n-m)+\sum_{j=1}^{m-|T|} \big[(i_{j}-1)(k-j+|T_{1}|+\cdots +|T_{j-1}|)+|R_{j}|(|T_{1}|+\cdots +|T_{j}|)\big].$$
\end{cor}

\begin{pf}
  Provided that the degree of $\overline{\mu}_{n} \in \As^{\ash}(n)$ is $n-1$
  and provided the formula for the coproduct in $\As^{\ash}$ given in
  \cite{LodayVallette}, chapter 8, where we include the decomposition involving
  $\overline{\mu}_{m} = |$ or $\overline{\mu}_{q} = |$,
$$\Delta_{(1)}(\overline{\mu}_{n}) = \sum_{\substack{p+q+r=n\\ p+1+r=m}} (-1)^{(q+1)r} (\overline{\mu}_{m}; \underbrace{\textrm{id, \ldots, id}}_{p},\, \overline{\mu}_{q},\, \underbrace{\textrm{id, \ldots, id}}_{r}),$$
the Proposition \ref{Koszuldualproduct} gives the description of $\q\uAs^{\ash}$ and of the infinitesimal decomposition map. The coproduct is given in the same way as explained in the proof of Proposition \ref{Koszuldualproduct} thanks to the coproduct in $\As^{\ash}$ given in \cite{LodayVallette} by
$$\Delta(\overline{\mu}_{n}) = \sum_{i_{1}+ \cdots + i_{m}=n} (-1)^{\epsilon'}(\overline{\mu}_{m}; \overline{\mu}_{i_{1}},\, \ldots ,\, \overline{\mu}_{i_{m}}),$$
where $\epsilon' := \sum_{j=1}^{m} (i_{j}-1)(k-j)$.
$\cqfd$
\end{pf}

\begin{prop}\label{KosquAs}
The operad $\q\uAs$ is Koszul, that is
$$\q\uAs \circ_{\kappa} \q\uAs^{\ash} \qiso I.$$
\end{prop}

\begin{pf}
  We remark that
$$\q\uAs \circ_{\kappa} \q\uAs^{\ash} \cong \vcenter{
  \xymatrix@M=0pt@R=4pt@C=3pt{&&\\
    & \ar@{{*}}[u] \ar@{-}[d]\\
    &}}
\oplus \left\{\bigoplus_{S \subseteq \underline{n}} (\As \circ_{\kappa}
  \As^{\ash}) (n) \right\}_{n \geq 1}.$$
Since $\vcenter{
  \xymatrix@M=0pt@R=4pt@C=4pt{
    & & & & \\
    & & & &\\
    & \ar@{{*}}[u] \ar@{-}[dr] & & \ar@{-}[dl] & \\
    & & \ar@{-}[d] & &\\
    & & }} = 0 = \vcenter{
  \xymatrix@M=0pt@R=4pt@C=4pt{& & & & \\
    & & & &\\
    & \ar@{-}[dr] & & \ar@{{*}}[u] \ar@{-}[dl] & \\
    & &\ar@{-}[d] & & \\
    & & }}$ in $\q\uAs$, the differential on $(\As \circ_{\kappa} \As^{\ash})
(n)$ is given by the usual differential on $\As \circ_{\kappa} \As^{\ash}$
except for $d(\vcenter{
  \xymatrix@M=0pt@R=4pt@C=3pt{&&\\
    & \ar@{{*}}[u] \ar@{-}[dd] &\\
    \ar@{.}[rr] &&\\
    &&}}) = \vcenter{
  \xymatrix@M=0pt@R=4pt@C=3pt{&&\\
    & \ar@{{*}}[u] \ar@{-}[d]\\
    &}}$. Moreover, we know that $(\As \circ_{\kappa} \As^{\ash}) (n) \qiso
I(n)$. Thus $\q\uAs \circ_{\kappa} \q\uAs^{\ash} \qiso I$.  $\cqfd$
\end{pf}

\begin{lem}
  The curved cooperad $u\As^{\ash}$ is equal to the curved cooperad
$$u\As^{\ash} = \left( \q\uAs^{\ash}, \Delta_{\q\uAs^{\ash}}, 0,
  \theta \right),$$ where $\Delta_{\q\uAs^{\ash}}$ was made explicit in Corollary
\ref{quAs} and
\begin{equation*}
  \thetaop (\muns) = \left\{
    \begin{array}{rl}
      -1  \cdot | &  \text{if } n = 2 \text{ and, } S = \{1 \} \text{ or } S = \{2\}\\
      0 \quad & \text{otherwise}
    \end{array}
  \right.
\end{equation*}

\end{lem}

\begin{pf}
 For the definitions given in \ref{inhomoquadprop}, we remark that the space of generators defining $\uAs$ satisfies Conditions (I) and (II) of Section \ref{inhomoquadprop}. According to the definition \ref{defKoszuldual}, we just have to compute the predifferential $d_{\uAs^{\ash}}$ and its curvature $\theta$. Since
  the relations in $u\As$ have no linear terms, the predifferential $d_{\uAs^{\ash}}
  = 0$. To compute $\theta$, we find the elements of weight 2, which
  correspond to the relations in $\q\uAs$. We identify each cooperation with the
  corresponding leading quadratic term of a relation in $\uAs$, and
  then assign to that operation the opposite of the corresponding constant term of the
  relation:

  \begin{equation*}
    \begin{split}
      \ass &\longleftrightarrow \mthree \mapsto 0 \\
      \mun &\longleftrightarrow \mun \mapsto -1 \cdot | \\
      \unm &\longleftrightarrow \unm \mapsto -1 \cdot |
    \end{split}
  \end{equation*}
  $\cqfd$
\end{pf}

\begin{thm}
  The cobar construction on the Koszul dual curved cooperad associated
  to $\uAs$ provides a cofibrant resolution of $\uAs$
$$\uAsinf := \Omega \uAs^{\ash} \qiso \uAs.$$
\end{thm}

\begin{pf}
  By Proposition \ref{KosquAs}, $\q\uAs$ is Koszul, and then Theorem
  \ref{resolution} gives the result.  $\cqfd$
\end{pf}

We now make the operad $\uAsinf$ more explicit.

The underlying operad of the dg operad $\Omega \, {\uAs}^{\ash}$ is
the free operad
$\mathcal{F}\left(s^{-1}\overline{\q\uAs^{\ash}}\right) = \mathcal{F}
\left(s^{-1} \left\{ \muns \right\} \right)$, $n \geq 2$, $S
\subset \underline{n}$ and $n=1, S=\{1\}$, giving a free generating
set for $\Omega \, {\uAs}^{\ash}$. As a derivation of the composition
structure, the differential $d = d_0 + 0 -d_2$ is completely
determined by its action on the generators:

\begin{equation}\label{bigD}
 \left\{ \begin{split}
    \mun &\mapsto \munsplit - |\\
    \unm &\mapsto \unmsplit - |\\
    \smfivectwo &\mapsto%
    \Sigma (-1)^{\epsilon} \, \vcenter{\xymatrix@M=0pt@R=7pt@C=7pt{
        & & & \\
        & \ar@{-}[dr] & \ar@{-}[d] \ar@{{*}}[u] & \ar@{-}[dl] \\
        & & \ar@{-}[d] & \\
        & & & \\
        \ar@{-}[dr] \ar@{{*}}[u] &\ar@{-}[d] &\ar@{-}[dl] & \\
        &\ar@{-}[d] & & \\
        & & & }}
  \end{split}\right.
\end{equation}
where the last line means: for $(n,\, S) \neq (2,\, \{ 1\})$ and  $(n,\, S) \neq (2,\, \{ 2\})$, we have
\begin{equation*}
  d(\munsc) = \sum_{\substack{p+q+r=n\\
      p+1+r=m}} (-1)^{q(r+|S_{1}|)+|S_{2}||S_{1}'| + p+1} (\overline{\mu}_{m}^{S_{1}}; \underbrace{\textrm{id, \ldots, id}}_{p-|S'_{1}|},\,
  \overline{\mu}_{q}^{S_{2}},\, \underbrace{\textrm{id, \ldots,
      id}}_{r-|S''_{1}|}),
\end{equation*}

\begin{rem}
  On the right-hand side of equation (\ref{bigD}), the two-level trees
  now represent the compositions in the free operad.
\end{rem}

We obtain the following description for a $\uAsinf$-algebras structure.

\begin{prop}\label{assinfstructure}
A \emph{$\uAsinf$-algebra structure on a dg module $(A,\, d_{A})$} is given by a
collection of maps, $\mu_{1}^{\{1\}},\, \{ \muns \}_{n \geq 2,\, S \subset
  \underline{n}}$ where each $\muns $ is a map $
A^{\otimes (n-|S|)} \to A$ of degree $n+|S|-2$ which together satisfy the following identities:
  $$\left\{ \begin{array}{lcl}
      \partial(\mu_{2}^{\{1 \}}) &=& \mu_{2}^{\emptyset} \circ(\mu_{1}^{\{1\}}, -) - \id_A \\
      \partial(\mu_{2}^{\{2 \}}) &= & \mu_{2}^{\emptyset} \circ (- ,\mu_{1}^{\{1\}}) - \id_A
     \end{array} \right.$$
and for $(n,\, S) \neq (2,\, \{ 1\})$ and  $(n,\, S) \neq (2,\, \{ 2\})$
  $$\partial(\muns ) = \sum_{\substack{p+q+r=n\\
          p+1+r=m}}
      (-1)^{q(r+|S_{1}|)+|S_{2}||S_{1}'| + p+1} \mu_{m}^{S_{1}} \circ
      (\underbrace{\id, \ldots, \id}_{p-|S'_{1}|},\, \mu_{q}^{S_{2}},\, \underbrace{\id, \ldots,
          \id}_{r-|S''_{1}|}).$$
\end{prop}

\begin{pf}
Since $\uAsinf$ is a quasi-free operad, a map $\mu_{A} : \uAsinf (A) \rightarrow A$ of degree $0$ is determined by a collection of maps, $\mu_{1}^{\{1\}},\, \{ \muns \}_{n \geq 2,\, S \subset \underline{n}}$ where each $\muns $ is a map $A^{\otimes (n-|S|)} \to A$ of degree $n+|S|-2$. The fact that the map $\mu_{A}$ is a dg map gives the relations with the notations $\mu_{n}^{S} (a_{1} \otimes \cdots \otimes a_{n-|S|}) := \mu_{A}(\overline{\mu}_{n}^{S} \otimes a_{1} \otimes \cdots \otimes a_{n-|S|})$.
$\cqfd$
\end{pf}

\begin{rem}
This notion of $\uAsinf$-algebra corresponds to the notion of homotopy unit for an $A_{\infty}$-algebra given in \cite{Fukaya}.
\end{rem}

\subsection{Infinity-morphisms} \label{inftymorphisms} 
Following the classical case, we describe the
\emph{infinity-morphisms} of algebras over the Koszul resolution of a Koszul
inhomogeneous quadratic operad. We give explicit formulae for infinity-morphism of $\uAsinf$-algebras.

Unless we indicate otherwise, for the rest of this section, $\Po$ will
denote a Koszul inhomogeneous quadratic operad, $\Poa$ its curved Koszul dual cooperad and $\Po_{\infty} := \Omega \Poa$ denotes the Koszul resolution of $\Po$ (see Section \ref{curvedKD}).

Let $A$ be a $\Po_{\infty}$-algebra, and denote its structure map by
  $\mu_A \in \Hom_{\textsf{dg op}} (\Po_{\infty},
    \End_A).$ Then by the bar-cobar adjunction
  \ref{barcobaradjunction}, we have 
\begin{equation*}
\Hom_{\textsf{dg operads}}(\Omega \Po^{\ash}, \End_A) \cong
\Tw ( \Po^{\ash}, \End_A).
\end{equation*}
By classical $\Hom$-tensor duality, we have the bijection
$$\begin{array}{ccc}
\Hom_{\mathbb{S}\textsf{-Mod}}(\Po^{\ash},\End_A) & \cong & \Hom_{\textsf{dg mod}}(\Po^{\ash}(A), A)\\
\mu_{A} & \longmapsto & d_{\mu_{A}}.
\end{array}$$
We recall the classical lemma, that we can find for example in \cite{LodayVallette}

\begin{lem}
A coderivation of $\Poa(A)$ is completely characterized by its corestriction to the generators
$$\begin{array}{ccc}
\Hom_{\mathsf{mod.}}(\Po^{\ash}(A), A) & \cong & \mathrm{Coder}(\Po^{\ash}(A))\\
d_{\mu_{A}} & \longmapsto & D_{\mu_{A}}^{r}.
\end{array}$$
\end{lem}

We call a \emph{curved codifferential} any coderivation $D$ of degree $-1$ which satisfies
$$D^{2} = (\theta \circ \id_{\Poa(A)}) \circ \Delta_{\Poa(A)}.$$

We have the following extension of a classical result about codifferentials:

\begin{lem}
A $\Po_{\infty}$-algebra structure on $A$ is equivalent to a codifferential on $\Poa(A)$
$$\begin{array}{ccl}
\Tw (\Poa,\, \End_{A}) & \cong & \mathrm{curCodiff}(\Po^{\ash}(A))\\
\mu_{A} & \longmapsto & \hspace{0.5cm} D_{\mu_{A}} := d_{\Poa(A)} + D_{\mu_{A}}^{r}.
\end{array}$$
\end{lem}

\begin{pf}
The predifferential $d_{\Poa}$ is a coderivation so the map $D_{\mu_{A}} := d_{\Poa(A)} + D_{\mu_{A}}^{r}$ is a coderivation. The construction here is the same as the construction in Section \ref{barconstruction} with $D_{\mu_{A}}^{r} = d_{\mu_{A}}^{r} \circ_{\Po} \id_{A}$, so $\mu_{A} \in \Tw (\Poa,\, \End_{A})$ implies $D_{\mu_{A}}^2 = (\theta \circ \id_{\Po^{\ash}(A)}) \circ \Delta_{\Poa(A)}$.

According to the proof of Lemma \ref{forcodiff}, we only have to remark that $D^{r}_{\partial(\mu_{A})+\mu_{A} \star \mu_{A} - \varTheta} = D_{\mu_{A}}^{2} - (\theta \circ id_{\Poa(A)}) \circ \Delta_{\Poa (A)} = 0$ implies $d_{A} \circ d_{\mu_{A}} + d_{\mu_{A}} \circ d_{\Poa(A)} + d_{\mu_{A} \star \mu_{A}} - d_{\varTheta} = (D^{r}_{\partial(\mu_{A})+\mu_{A} \star \mu_{A} - \varTheta})^{|A} = 0$. Since $d_{A} \circ d_{\mu_{A}} + d_{\mu_{A}} \circ d_{\Poa(A)} + d_{\mu_{A} \star \mu_{A}} - d_{\varTheta}$ is sent to $\partial(\mu_{A})+\mu_{A} \star \mu_{A} - \varTheta$ and $0$ is sent to $0$ by reversing the bijection in the Hom-tensor duality, we get the result.
$\cqfd$
\end{pf}

\subsubsection{\bf Infinity-morphism of $\Po_{\infty}$-algebras}\label{infmorph}
Let $A$ and $B$ be two $\Po_{\infty}$-algebras, with structure maps $\mu_A$ and
$\mu_B$. A \emph{$\infty$-morphism $A \rightsquigarrow B$ of $\Po_{\infty}$-algebras} is
a dg $\Po^{\ash}$-coalgebras map
$$F: (\Po^{\ash}(A),\, D_{\mu_{A}}) \to (\Po^{\ash}(B),\, D_{\mu_{B}}).$$

This description of $\infty$-morphisms makes it clear that
$\Po_{\infty}$-algebras, $\infty$-morphisms, and composition given by
composition of dg $\Po^{\ash}$-coalgebra maps forms a category.

A $\uAs^{\ash}$-coalgebras map $F : \uAs^{\ash} (A) \to \uAs^{\ash}(B)$ is characterized by its corestriction to $B$, that is $F$ is determined by a collection of maps $f_{n}^{S}: A^{\otimes (n - |S|)} \to B$. The fact that $F$ commutes with the differentials is equivalent to a family of equations on the $f_{n}^{S}$.
Pictorially, the collection of maps $f_{n}^{S}$ satisfy: 
$${\small \partial \left(\vcenter{\xymatrix@M=0pt@R=8pt@C=10pt{
    & & & & \\
    & & & & & \\
    \ar@{-}[drrr] & \ar@{{*}}[u] \ar@{-}[drr] & \ar@{-}[dr] & \ar@{-}[d] & \ar@{-}[dl] \ar@{{*}}[u]  \
& \ar@{-}[dll] \\
    & & & f_{n}^{S} \ar@{-}[d] & \\
    & & &
    }}\right) = \sum \pm\ \vcenter{\xymatrix@M=0pt@R=4pt@C=11pt{
        && \ar@{-}[dd] & \ar@{-}[ddl] \\
        & \ar@{-}[dr] \ar@{{*}}[u] &&\\
        && \mu_{q}^{S_{2}}(A) \ar@{-}[d] & \\
        \ar@{-}[ddrr] && \ar@{-}[dd] &&& \ar@{-}[ddlll]\\
        &&& \ar@{-}[dl] \ar@{{*}}[u] &&&\\
        && f_{m}^{S_{1}} \ar@{-}[d] &&& \\
         & & }} - \sum \pm \vcenter{\xymatrix@M=0pt@R=4pt@C=11pt{
        & & & & \\
        \ar@{-}[ddr] &&& \ar@{-}[ddr] && \ar@{-}[ddl] &&& \ar@{-}[dd] \\
        && \ar@{-}[dl] \ar@{{*}}[u] &&&&&&\\
        & f_{i_{1}}^{T_{1}} \ar@{-}[dddrrr] &&& f_{i_{2}}^{T_{2}} \ar@{-}[d] &&&& f_{i_{m-|T|}}^{T_{m-|T|}} \ar@{-}[dddllll]\\
        &&&& \ar@{-}[dd] &&&&\\
        &&&&& \ar@{-}[dl] \ar@{{*}}[u] &&&\\
        &&&& \mu_{m}^{T}(B) \ar@{-}[d] &&&&\\
        &&&&&&&&.}}}$$

\begin{prop}\label{homotopyunit}
Let $A$, $B$ be two $\uAsinf$-algebras, and let $\muns(A),
\muns (B)$ be the respective structure maps. An $\infty$-morphism between $A$ and $B$ is a collection of maps
$$\{f_{n}^{S}: A^{\otimes (n - |S|)} \to B\}_{n\geq 1,\, S \subset \underline{n}} \textrm{ of degree } n+|S|-1,$$
satisfying: for $n=1$, $d_{A} \circ f_{1}^{\emptyset} = f_{1}^{\emptyset} \circ d_{A}$, that is $f_{1}^{\emptyset}$ is a chain map, and for $n>2$,\\
$\partial(f_{n}^{S}) =$
$$\sum_{\substack{p+q+r=n\\ p+1+r=m}} (-1)^{p+q(r+|S_{1}|)+|S_{2}||S''_{1}|} f_{m}^{S_{1}} \circ (\underbrace{id_{A},\, \ldots ,\, id_{A}}_{p-|S'_{1}|},\, \mu_{q}^{S_{2}}(A),\, \underbrace{id_{A},\, \ldots ,\, id_{A}}_{r-|S''_{1}|})$$
$$- \hspace{-0.5cm} \sum_{\substack{i_{1}-|T_{1}|+ \cdots
  +i_{m-|T|}-|T_{m-|T|}|=n-|S|}} \epsilon (-1)^{(m+|T|-1)(n-m+|S|-|T|)} \mu_{m}^{T}(B) \circ \big(f_{i_{1}}^{T_{1}},\, \ldots ,\, f_{i_{m-|T|}}^{T_{m-|T|}}\big),$$
where $\left\{ \begin{array}{lll} S'_{1}
    & \subset & \underline{p}\\ S_{2} & \subset & \underline{q}\\
    S''_{1} & \subset & \{p+2,\, \ldots,\, p+1+r\} \end{array}
\right.$ such that $S = S'_{1} \sqcup (S_{2} + p) \sqcup (S''_{1} +
q)$ and $S_{1} = S'_{1} \sqcup S''_{1}$, where $\left\{ \begin{array}{lll} T & \subset & \underline{m}\\ T_{j} & \subset & \underline{i_{j}} \end{array} \right.$ such that $T = R_{0} \sqcup \ldots \sqcup R_{m-|T|}$ and
$$\begin{array}{lll} S & = & R_{0} \sqcup (T_{1}+|R_{0}|) \sqcup (R_{1}+i_{1}) \sqcup \ldots \sqcup \\ && (T_{m-|T|}+|R_{0}|+\cdots +|R_{m-|T|-1}|+i_{1}+\cdots +i_{m-|T|-1}) \sqcup (R_{m-|T|}+i_{1}+\cdots +i_{m-|T|}) \end{array}$$ and where $\epsilon := |T|(n-m)+\sum_{j=1}^{m-|T|} \big[(i_{j}-1)(k-j+|T_{1}|+\cdots +|T_{j-1}|)+|R_{j}|(|T_{1}|+\cdots +|T_{j}|)\big]$.
\end{prop}

\begin{pf}
An $\infty$-morphism $A \rightsquigarrow B$ is a $\uAs^{\ash}$-coalgebras morphism $F : \uAs^{\ash}(A) \rightarrow \uAs^{\ash}(B)$. Such a morphism is completely determined by its image on the cogenerators of $\uAs^{\ash}(B)$, that is by a map $f : \uAs^{\ash}(A) \rightarrow B$ (of degree $0$), or equivalently by a collection of maps $\{f_{n}^{S}: A^{\otimes (n - |S|)} \to B\}_{n\geq 1,\, S \subset \underline{n}}$ of degree $n+|S|-1$. The fact that $F$ commutes with the predifferential is equivalent to the following commutative diagram
$$\xymatrix{\uAs^{\ash}(A) \ar[r]^{\hspace{-0.5cm} \Delta \circ id_{A}} \ar[d]_{d_{1}+d_{2}} & \uAs^{\ash} \circ \uAs^{\ash}(A) \ar[r]^{\hspace{0.5cm} id \circ f} & \uAs^{\ash}(B) \ar[d]^{d_{B}+d_{2}^{|B}}\\
\uAs^{\ash}(A) \ar[rr]_{f} && B.}$$
Making this diagram explicit gives exactly the formulae of the Proposition.
$\cqfd$
\end{pf}

\begin{ex1}
For $n = 1$ and $S = \{1\}$, the formula gives
$${\small \partial \left(\vcenter{\xymatrix@M=0pt@R=7pt@C=10pt{
     \\
     f_{1}^{\{1\}} \ar@{-}[d] \ar@{{*}}[u]\\
    &}} \hspace{-.35cm}\right) = \ \vcenter{\xymatrix@M=0pt@R=7pt@C=9pt{
        &\\
        \mu_{1}^{\{1\}}(A) \ar@{-}[d] \ar@{{*}}[u]\\
        f_{1}^{\emptyset} \ar@{-}[d] \\
        & }} \hspace{-.3cm} - \vcenter{\xymatrix@M=0pt@R=7pt@C=10pt{
        \\
        \mu_{1}^{\{1\}}(B) \ar@{-}[d] \ar@{{*}}[u]\\
        &}}}\! ,$$
that is the element $f_{1}^{\{1\}}$ bounds the failure of $f_{1}^{\emptyset}$ to preserve the unit.
\end{ex1}

Before we end the section, we use the results above to give the
following definition.

\begin{dei}
A $\infty$-morphism of $\Po_{\infty}$-algebras $F:A \rightsquigarrow B$ is a \emph{quasi-isomorphism}
if the chain map $f_{1}^{\emptyset}:A \to B$ induces an isomorphism in homology.
\end{dei}

\subsection{Rectification}

We now prove that for every $\uAsinf$-algebra $A$ there is a
universal $\infty$-quasi-morphism $I_A$ between $A$ and a
$\uAs$-algebra. This universal morphism takes the form of the unit of
an adjunction. We make use of the
bar and cobar constructions of algebras over Koszul operads (Sections
\ref{barconstruction}, \ref{cobarconstruction}) for $\uAsinf$-algebras and $\uAs$-algebras.\\

The twisting morphisms $\iota : \uAs^{\ash} \rightarrow \Omega \uAs^{\ash} = \uAsinf$ and $\kappa : \uAs^{\ash} \rightarrow \uAs$ are defined in Section \ref{sectbarcobaradjunction} and \ref{Koszuldef}.

\begin{lem}\label{projbarcobarres}
Let $A$ be $\uAsinf$-algebra. The morphism of dg $\So$-modules $A
\mono \Omega_{\kappa}\B_{\iota}A$ is a quasi-isomorphism.
\end{lem}

\begin{pf}
We endow $\uAs \circ_{\kappa} \uAs^{\ash} \circ_{\iota} \Omega \uAs^{\ash}$ with a filtration $F_{p}$ given by
$$F_{p}(\uAs \circ \uAs^{\ash} \circ \Omega \uAs^{\ash}) = \bigoplus_{\omega+m\leq p} (\uAs \circ \uAs^{\ash})^{(\omega)} \circ (\Omega \uAs^{\ash})_{m}.$$
Moreover we endow $\Omega \uAs^{\ash}$ with a filtration given by the homological degree, so that the morphism $\Omega \uAs^{\ash} \mono \uAs \circ_{\kappa} \uAs^{\ash} \circ_{\iota} \Omega \uAs^{\ash}$ preserves the filtrations. Since the weight on $\uAs \circ_{\kappa} \uAs^{\ash}$ is non-negative and $\Omega \uAs^{\ash}$ is non-negatively graded, the filtrations are bounded below. Moreover, the filtrations are exhaustive. Thus, we can apply the classical theorem of convergence of spectral sequences (cf. Theorem 5.5.1 of \cite{Weibel}) to obtain
$$E_{p,q}^{\bullet} \Rightarrow \mathrm{H}_{p+q}(\uAs \circ_{\kappa} \uAs^{\ash} \circ_{\iota} \Omega \uAs^{\ash}) \textrm{ and } {E'}_{p,q}^{\bullet} \Rightarrow \mathrm{H}_{p+q}(\Omega \uAs^{\ash})$$
and an induced morphism between the spectral sequences. The differential on $E_{p,q}^{0}$ coincides with the differential on $\q\uAs \circ \q\uAs^{\ash}$, so Proposition \ref{KosquAs} shows that $E_{p,q}^{1} \cong {E'}_{p,q}^{1}$. It follows that $E_{p,q}^{r} \cong {E'}_{p,q}^{r}$ for all $r\geq 1$ and we get  that $\Omega \uAs^{\ash} \qiso \uAs \circ_{\kappa} \uAs^{\ash} \circ_{\iota} \Omega \uAs^{\ash}$.

We have $\Omega_{\kappa}\B_{\iota}A \cong (\uAs \circ_{\kappa} \uAs^{\ash} \circ_{\iota} \uAsinf) \circ_{\uAsinf} A$. The short exact sequence
$$(\uAs \circ_{\kappa} \uAs^{\ash} \circ_{\iota} \uAsinf) \circ \uAsinf \circ A \rightarrow (\uAs \circ_{\kappa} \uAs^{\ash} \circ_{\iota} \uAsinf) \circ A \rightarrow (\uAs \circ_{\kappa} \uAs^{\ash} \circ_{\iota} \uAsinf) \circ_{\uAsinf} A$$
induces a long exact sequence in homology. Since we work over a field of characteristic $0$, the ring $\K[\So_{n}]$ is semi-simple by Maschke's theorem, that is every $\K[\So_{n}]$-module is projective. So the K\"unneth formula implies that $\mathrm{H}_{\bullet}((\uAs \circ_{\kappa} \uAs^{\ash} \circ_{\iota} \uAsinf) \circ \uAsinf \circ A) \cong \mathrm{H}_{\bullet}(\uAsinf) \circ \mathrm{H}_{\bullet}(\uAsinf) \circ \mathrm{H}_{\bullet}(A)$ and $\mathrm{H}_{\bullet}((\uAs \circ_{\kappa} \uAs^{\ash} \circ_{\iota} \uAsinf) \circ A) \cong \mathrm{H}_{\bullet}(\uAsinf) \circ \mathrm{H}_{\bullet}(A)$. Finally, this gives that $\mathrm{H}_{\bullet}((\uAs \circ_{\kappa} \uAs^{\ash} \circ_{\iota} \uAsinf) \circ_{\uAsinf} A) \cong \mathrm{H}_{\bullet}(\uAsinf) \circ_{\mathrm{H}_{\bullet}(\uAsinf)} \mathrm{H}_{\bullet}(A) \cong \mathrm{H}_{\bullet}(A)$.
$\cqfd$
\end{pf}

\begin{thm}[Universal rectification]\label{universalrectification}
Let $A$ be a $\uAsinf$-algebra. There is a dg $\uAs$-algebra, $\Omega_{\kappa} \mathrm{B}_{\iota}A$ and an $\infty$-quasi-isomorphism $I_{A}: A \rightsquigarrow^{\hspace{-.3cm} \sim} \hspace{.1cm} \Omega_{\kappa} \mathrm{B}_{\iota} A$
so that for any dg $\uAs$-algebra $B$ and any $\infty$-morphism $F : A \rightsquigarrow B$, there is a unique dg $\uAs$-algebras map $\tilde{f}: \Omega_{\kappa} \B_{\iota} A \to B$ so that
$F = f \circ I_{A}$, that is the following diagram commutes:
\begin{equation*}
\xymatrix{
\Omega_{\kappa} \mathrm{B}_{\iota} A \ar@{-->}[dr]^{\tilde{f}} \\
A \ar@{->}[r]_{F} \ar@{->}[u]^{I_{A}} & B 
}
\end{equation*}
\end{thm}

\begin{pf}
The map $I_{A}$ is defined by
\begin{equation*}
i_{n}^{S}(a_1, \ldots , a_{n - |S|}) =
\muns (a_1, \ldots , a_{n - |S|}) \in
\textrm{B}_{\iota} A \hookrightarrow \Omega_{\kappa}
\textrm{B}_{\iota} A.
\end{equation*}
By direct computation, this map is a $\infty$-morphism between the $\uAsinf$-algebras $A$ and $\Omega_{\kappa} \B_{\iota}A$. To see that this map is a quasi-isomorphism, observe that
$i_{1}^{\emptyset}$ is equal to the inclusion map in Lemma \ref{projbarcobarres}. To define the map $\tilde{f}$, we note that the $\infty$-morphism of $\uAsinf$-algebras $F$ is determined by the collection
of maps $f_{n}^{S}$, or by the collection of elements $\{f_{n}^{S}(a_1,\, \ldots ,\, a_{n- |S|})\}$
in $B$. We define the module map
$\textrm{B}_{\iota} A \to B$ by
\begin{equation*}
\muns (a_1,\, \ldots ,\, a_{n- |S|}) \mapsto f_{n}^{S}(a_1,\, \ldots ,\, a_{n- |S|}).
\end{equation*}
This map is a dg module map if and only if $F$ is a $\infty$-morphism. Since the $\uAs$-algebra $\Omega_{\kappa} \textrm{B}_{\iota}A$ is freely generated by  $\{ \muns (a_1,\, \ldots ,\, a_{n- |S|})\}$, we define the map $\tilde{f}$ to be the lift of
the above dg map to a $\uAs$-algebras map $\Omega_{\kappa} \textrm{B}_{\iota}A
\to B$. By construction we have $\tilde{f} \circ I_{A} = F$.
$\cqfd$
\end{pf}

Let us interpret the result above in terms of the categories of
algebras. \label{rectificationadjoint} Since we have an operad map $\uAsinf \twoheadrightarrow
\uAs$, we have an inclusion functor (one-to-one on objects and on morphisms) $\uAs \textsf{-alg}
\hookrightarrow \uAsinf \textsf{-alg}$, which we denote by $i$. Let's denote by $R$ the object-function that takes each
$\uAsinf$-algebra $A$ to the $\uAs$-algebra $R(A) = \Omega_{\kappa}
B_{\iota} (A)$. Because the arrow $A \xrightarrow{I_A} iR(A)$ is
universal, $R$ can be extended to morphisms so that it becomes a
functor from $\uAsinf \text{-alg} \to \uAs \text{-alg}$:
\begin{equation*}
\xymatrix{
  R(A) \ar@{.>}[r]^{R(F)} & R(B) \\
  A \ar@{->}[u] \ar@{->}[r]^{F} & B. \ar@{->}[u]
}
\end{equation*}

We summarize in the following proposition.
\begin{prop}
The functor $i$, the object-function $R$, and the universal morphisms $A \xrightarrow{I_A} iR(A)$
determine the extension of $R$ to a functor $R: \uAsinf \textsf{-alg}
\to \uAs \textsf{-alg}$ so that $I: \id \to iR$ is the unit of an adjunction:
$$\xymatrix{     *{ \uAs \textsf{-alg} \ \ } \ar@<.8ex>[r]^(.45){i} & *{\
\uAsinf \textsf{-alg}}  \ar@<.8ex>[l]^(.43){R}}.$$ 
\end{prop}

\begin{pf}
See Mac Lane \cite{Maclane} chapter 4, theorem 1.
\end{pf}

It is tempting to try to put a model category structure on the
right-hand side so that this pair of functors becomes some kind of Quillen
equivalence, as Lefevre-Hasegawa \cite{Lefevre-Hasegawa} did for $A_{\infty}$-algebras
and $\As$-algebras. (Actually, $A_{\infty}$-algebras are not quite a model category, see the referenced paper for more details) . Instead we
observe that each functor takes quasi-isomorphisms to
quasi-isomorphisms, and so each functor induces a functor between the homotopy
categories
(localizations of each category by its quasi-isomorphisms). We claim these induced functors are an
adjoint-equivalence of the homotopy categories.

\subsection{Transfer formulae}\label{transfersection}

In this section we provide formulae, based on labelled trees, for the pullback of a
$\uAsinf$-structure along a strong deformation retract.

For this entire section, suppose $V, A$ are chain complexes, and%
$$\xymatrix{     *{ V \ \ } \ar@<.8ex>[r]^(.45){i} & *{\
A \quad} \ar@(dr,ur)[]_h \ar@<.8ex>[l]^(.43){p}},$$
is a strong deformation retract, i.e., $p$ and $i$ are chain maps, where $p \circ i = \id_V$ and $\, d_Ah+hd_A =
\id_A - i \circ p$. Moreover, suppose $A$ is a $\uAsinf$-algebra, with
structure map $\mu_A$.

\begin{dei}
Let $n \geq 2, S \subset \underline{n}$, we define the set 
$\mathcal{T}_n^S$ be the set of planar, rooted
trees, with $n$ leaves, and a cork above each $i$th leaf if $i \in S$
which is labelled by either the word ``connected'' or
``disconnected.'' 
We define $\mathcal{T}_{1}^{\emptyset} = \{ \, |\, \}$ and $ \mathcal{T}_{1}^{\{1\}} = \{ \un^{\text{\emph{connected}}}\}$. 
\end{dei}

\begin{dei}
Let $T \in \mathcal{T}_n^S$, and let $v$ be any internal vertex in
$T$. We denote by $\inc(v)$ the ordered (left-to-right) set of incoming edges to the
vertex $v$. For each element $i \in \inc(v)$, we define $l_i$ and
$c_i$ as follows:
\begin{enumerate}
    \item $l_i$ is the total number of leaves without attached corks
      in the tree $T$ whose (unique)
      path to the root passes through edge $i$
    \item $c_i$ is the total number of incoming edges to $v$ without
      attached corks to the \emph{right} of edge $i$. 
\end{enumerate}
\end{dei}

\begin{dei}
For any $T \in \mathcal{T}_n^S$ and any internal
vertex $v \in T$, we define
\begin{equation*}
\epsilon(v) = \sum_{1 \leq i < j \leq
  \lvert \inc(v) \lvert} \left( l_{i} + 1 \right)l_{j} \quad + \sum_{\substack{i \in
    \inc(v)\\ \text{with a connected} \\ \text{cork on it}}}
c_i.
\end{equation*}

For any tree $T \in \mathcal{T}_{n}^{S}$,
we set
\begin{equation} \label{transfersigns}
\epsilon(T) = \sum_{\substack{\text{internal vertices} \\ v
    \in T}} \epsilon(v).
\end{equation}
\end{dei}

\begin{dei}
Let $g_{\text{structure}}: \mathcal{T}_{n}^{S} \to \mathrm{Hom}(V^{\otimes (n - |S|)},\, V)$ be the set map that takes an element $T \in \mathcal{T}_{n}^{S}$ and assigns to
each vertex $v$ the operation $\mu_{\inc(v)}^{S(v)}(A)$ where $S(v)$ are
the positions of the \emph{attached} corks, and the homotopy $h$ to
each internal edge (that is not the outgoing edge of an attached
cork), and the map $i$ to each leaf without a cork above it, and the
map $p$ to the root of the tree. After this assignment, one composes
the operations as indicated by internal edges to arrive at an
operation $V^{\otimes (\inc(v) - |S|)} \to V$. Let
$g_{\text{morphism}} : \mathcal{T}_{n}^{S} \to \mathrm{Hom}(V^{\otimes
  (n - |S|)},\, A)$ be the set map that takes an element $T \in
\mathcal{T}_{n}^{S}$ and assigns to the tree the same element as
$g_{\text{structure}} \left( T \right)$, but with a homotopy assigned
to the root, rather than the map $p$.
\end{dei}

\begin{ex}
Let $T$ be the element of $\mathcal{T}_{5}^{\{1,4\}}$ that looks like
$$T= \vcenter
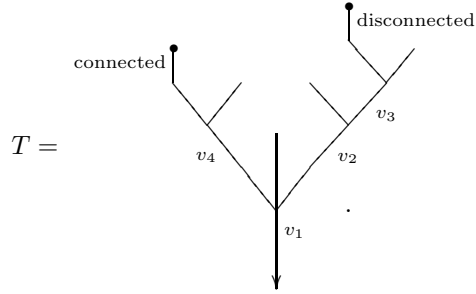
{
    \xymatrix@C=7pt@R=10pt{
      & && & & & & &  & & & & \\
      & & & & & &  *{} \ar@{-}[dr] \ar@{{*}}[u]_{\text{disconnected}} & &\ar@{-}[dl]  \\
      & *{} \ar@{-}[dr] \ar@{{*}}[u]^{\small \text{connected}}&  &*{} \ar@{-}[dl]
      & &  *{} \ar@{-}[dr] & & *{} \ar@{-}[dl]^{v_3} &  & & \\
      & & *{} \ar@{-}[dr]_{v_4} &  &\ar@{-}[d] & & 
     *{} \ar@{-}[dl]^{v_2} & & & & \\
      % &&*{} \ar@{-}[dr] & & & &*{} \ar@{-}[dl] & & & & & \\
      && &*{} \ar@{-}[dr] & *{} \ar@{-}[d] &*{} \ar@{-}[dl] & & & & & & \\
      && & &*{} \ar@{-}[d]^{v_1} & &. & & & & & \\
      && & &*{} \ar@{->}[d] & & & & & & & \\
      && & & & & & & & & &
    }}$$
The sign $(-1)^{\epsilon(T)}$ for this tree is given by
$$\begin{array}{lcl}
\epsilon(T) &=& \epsilon(v_1) + \epsilon(v_2) + \epsilon(v_3)
+\epsilon(v_4)\\
 &=& [(1+1) \cdot 1 + (1+1)\cdot 3 + (1+1)\cdot 3]+ [(1+1)\cdot 2]+ [(1+1)\cdot 1 + 1]+ [0]\\
 & \equiv & 1 \text{ mod } 2.
 \end{array}$$
The operation assigned to the tree $T, g_{\text{structure}}(T)$, is given by the
following composition of operations:
$$g_{\text{structure}}(T)= \vcenter
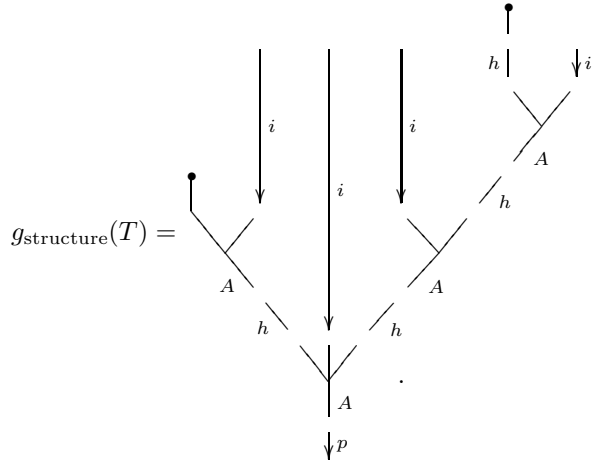
{
    \xymatrix@C=7pt@R=10pt{
      & & & & & & & & & & & & \\
      & &\ar@{->}[dddd]^i & &\ar@{->}[ddddddd]^i & &\ar@{->}[dddd]^i &
      &
      &\ar@{-}[d]_h \ar@{{*}}[u] & &\ar@{->}[d]^i \\
      & & & & & & & & & \ar@{-}[dr] & &\ar@{-}[dl]  \\
      && & & & & & & & & *{} \ar@{-}[dl]^A & \\
      && & & & & & & &\ar@{-}[dl]^h  & & \\
      *{} \ar@{-}[dr] \ar@{{*}}[u] &  &\ar@{-}[dl] & & & & \ar@{-}[dr] & & \ar@{-}[dl] & & & \\
      & *{} \ar@{-}[dr]_A & & & & & & *{} \ar@{-}[dl]^A & & & & \\
      &&\ar@{-}[dr]_h & & & &\ar@{-}[dl]^h & & & & & \\
      && &\ar@{-}[dr] &\ar@{-}[d] &\ar@{-}[dl] & & & & & & \\
      && & &*{} \ar@{-}[d]^A & &. & & & & & \\
      && & &\ar@{->}[d]^p & & & & & & & \\
      && & & & & & & & & & }}$$
while the morphism assigned to the tree $T, g_{\text{morphism}}\left(
  T \right)$ is given by:
$$g_{\text{morphism}}(T)= \vcenter{
    \xymatrix@C=7pt@R=10pt{
      & & & & & & & & & & & & \\
      & &\ar@{->}[dddd]^i & &\ar@{->}[ddddddd]^i & &\ar@{->}[dddd]^i &
      &
      &\ar@{-}[d]_h \ar@{{*}}[u] & &\ar@{->}[d]^i \\
      & & & & & & & & & \ar@{-}[dr] & &\ar@{-}[dl]  \\
      && & & & & & & & & *{} \ar@{-}[dl]^A & \\
      && & & & & & & &\ar@{-}[dl]^h  & & \\
      *{} \ar@{-}[dr] \ar@{{*}}[u] &  &\ar@{-}[dl] & & & & \ar@{-}[dr] & & \ar@{-}[dl] & & & \\
      & *{} \ar@{-}[dr]_A & & & & & & *{} \ar@{-}[dl]^A & & & & \\
      &&\ar@{-}[dr]_h & & & &\ar@{-}[dl]^h & & & & & \\
      && &\ar@{-}[dr] &\ar@{-}[d] &\ar@{-}[dl] & & & & & & \\
      && & &*{} \ar@{-}[d]^A & &. & & & & & \\
      && & &\ar@{->}[d]^h & & & & & & & \\
      && & & & & & & & & & }}$$
\end{ex}

\begin{prop} \label{transferformulae}\label{transferthm}
The maps defined by
\begin{equation} \label{transfer}
\muns(V) := \sum_{T \in \mathcal{T}_{n}^{S}} (-1)^{\epsilon(T)} g_{\text{structure}} \left( T \right).
\end{equation}
give $V$ the structure of a $\uAsinf$-algebra. Moreover, the maps
defined by
\begin{equation}
i_{n}^{S} := \sum_{T \in \mathcal{T}_{n}^{S}} (-1)^{\epsilon(T)}
g_{\text{morphism}} \left( T \right).
\end{equation}
provide a $\infty$-quasi-isomorphism of $\uAsinf$-algebras $I: V \rightsquigarrow^{\hspace{-.3cm} \sim} \hspace{.1cm} A$.
\end{prop}

\begin{pf}
  A combinatorial argument similar to the argument for transferring
  $A_{\infty}$-structures \cite{Markl2} will suffice.  $\cqfd$
\end{pf}

\begin{ex}
For small values of $n$, the transferred structure is given by
$$\begin{array}{l}
    \mu_{1}^{\{1 \}}(V) := p \circ \mu_{1}^{\{1\}}(A) =\, \vcenter{
      \xymatrix@C=7pt@R=10pt{
        & & \\
        \ar@{{*}}[u]_A  \ar[d]^{p} & \\
        & } }\\
    \mu_{2}^{\emptyset}(V) := p \circ \mu_{2}^{\emptyset}(A) \circ i^{\otimes 2} =
    \vcenter{ \xymatrix@C=6pt@R=8pt{
        \ar[d]^{i}  & & \ar[d]^{i} \\
        \ar@{-}[rd] & & \ar@{-}[ld] \\
        & *{} \ar@{-}[d]^A &&&& \\
        & \ar[d]^p & &  \\
        & & & } } \\
 \mu_{2}^{\{1\}}(V) := 
\vcenter{\xymatrix@C=6pt@R=8pt{
    & & &\ar@{->}[dd]^i & \\
    \ar@{{*}}[u]_A \ar@{-}[dr]_{h}& & & & \\
    & \ar@{-}[dr] & & \ar@{-}[dl] & \\
    & &*{} \ar@{-}[d]^A & \\
    & & \ar@{->}[d]^{p} & & \\
    & & & &
  }}
- \vcenter{\xymatrix@C=7pt@R=10pt{
    & & & \\
    & &\ar@{->}[d]^{i} \\
    *{} \ar@{{*}}[u] \ar@{-}[dr] & &  \ar@{-}[dl] \\
    &*{} \ar@{-}[d]^A  & \\
    & \ar@{->}[d]^{p} & \\
    & & & 
  }} \\
 i_{2}^{\{2\}}(V) := 
\vcenter{\xymatrix@C=6pt@R=8pt{
    &\ar@{->}[dd]^i & & & \\
    & & & &\ar@{{*}}[u]_A \ar@{-}[dl]_{h} \\
    & \ar@{-}[dr] & & \ar@{-}[dl] & \\
    & &*{} \ar@{-}[d]^A & \\
    & & \ar@{->}[d]^{h} & & \\
    & & & &
  }}
+ \vcenter{\xymatrix@C=7pt@R=10pt{
    & & & \\
   \ar@{->}[d]^{i} & & \\
     \ar@{-}[dr] & & *{}  \ar@{-}[dl]  \ar@{{*}}[u] \\
    &*{} \ar@{-}[d]^A  & \\
    & \ar@{->}[d]^{h} & \\
    & & & 
  }}

\end{array}$$
\end{ex}

For the reader familiar with transfer of $A_{\infty}$-structures,
restricting attention to the operations $\mu_{n}^{\emptyset}(V)$
recovers the familiar transfer formulae \cite{Kadeishvili, Merkulov, Kontsevich, Markl2, LodayVallette}.

\begin{rem}
Though our signs differ from \cite{Markl2}, we use his ideas to
develop a coherent sign convention for our transfer formulae. The
reader should note that our function $\epsilon(v)$ differs from the
$\thetaop (v)$ in \cite{Markl2} even on the operations $\mu_{n}^{\emptyset}(V)$, in small
ways, such as right-to-left orientation of trees instead of left-to-right.
\end{rem}

\subsection{Comparing Unital-(Associative-infinity) and
  (Unital-Associative)-infinity}
In previous sections, we have developed the definition of the operad
$\uAsinf$ whose algebras are \emph{homotopy unital}
$A_{\infty}$-algebras. There have been definitions of homotopy unital $A_{\infty}$-algebras \cite{Fukaya, Kontsevich, Lyubashenko2}, and these notions have been compared in \cite{Lyubashenko}. There have also been definitions of
\emph{strictly unital} $A_{\infty}$-algebras \cite{Kontsevich,
  Fukaya}---we will refer to these as $\suAinf$-algebras throughout
this section---they may be thought of as unital-(associative-infinity)
algebras as opposed to our (unital-associative)-infinity algebras. We will compare $\uAsinf$-algebras to
$\suAinf$-algebras. This comparison includes Theorem \ref{qisostrict},
which states that every $\uAsinf$-algebra has an equivalent
unital-$A_{\infty}$-structure on its homology. We demonstrate that
this theorem is fairly general, and applies to many algebraic
structures with units, including unital commutative associative
algebras, unital Batalin-Vilkovisky algebras, and
co-algebraic versions of these structures.\\

First we define $\suAinf$-algebras and their $\infty$-morphisms.

\begin{dei}\label{suainf}
  An \emph{$\suAinf$-algebra $(A, \{\mu_n\}_{n\geq 1},\, u)$} is an
  $A_{\infty}$-algebra $(A, \{\mu_n\}_{n\geq 1})$ with $u \in A$ such
  that $d_A (u)=0$ and $u$ is a left and right unit for $\mu_2$, and
  $u$ annihilates $\mu_n$ for $n \geq 3$ \cite{Kontsevich}.
\end{dei}

\begin{rems}
  \begin{enumerate}
	\item There exists a dg-operad whose algebras are precisely
  $\suAinf$-algebras, and we denote it by $\suAinf$. Furthermore, the
  operad $\suAinf$ is the quotient of $\uAsinf$ by the ideal generated
  by $\left\{ \muns \right\}_{n \geq 2, \, |S| \geq 1}$. A quick computation yields that this map is a quasi-isomorphism.
	\item The operad $\suAinf$ is not cofibrant. If it were, the lifting property would imply that it is a retract of $\uAsinf$ by the quotient map $\uAsinf \qiso \suAinf$, which a computation shows is impossible.
  \end{enumerate}	
\end{rems}

We now describe a diagram of categories of algebras. We will use the
following notation
\begin{itemize}
\item $\As$\textsf{-alg}: the category of associative algebras with
  algebra homomorphisms
\item $\uAs$\textsf{-alg}: the category of unital associative algebras
  with algebra homomorphisms that preserve the unit
\item $\infty$-$A_{\infty}$-\textsf{alg}: the category of
  $A_\infty$-algebras with $\infty$-morphisms
\item $\infty$-$\uAsinf$-\textsf{alg}: the category of $\uAsinf$
  algebras with $\infty$-morphisms
\item $\suAinf$-\textsf{alg}: the category of $\suAinf$-algebras with
  the $A_{\infty}$ $\infty$-morphisms for which $f_1$ preserves the
  unit and $f_n$ annihilates it (for $n \geq 2$)
\end{itemize}

First, we have the following diagram of operads:
\begin{equation*}
  \xymatrix@H=15pt{
    \uAs & \ar@{->>}[l]_{\sim} \textsf{su}A_{\infty} & \ar@{->>}[l]_{\quad
      \quad \sim} \uAsinf \\
    \As \ar@{>->}[u] & & \ar@{->>}[ll]_{\sim} \ar@{>->}[u] A_{\infty} &
  }
\end{equation*}

On the categories of algebras, the diagram becomes:
\begin{equation*}
  \xymatrix@H=15pt{
    \uAs \textsf{-alg} \quad  \ar@{->}[r] \ar@{->}[d] &
    \textsf{su}A_{\infty} \textsf{-alg} \quad
    \ar@{->}[r] & \infty \textsf{-}\uAsinf \textsf{-alg} \quad \ar@{->}[d]
    \\
    \As \textsf{-alg} \quad \ar@{->}[rr] & & \infty\textsf{-}A_{\infty} \textsf{-alg} \quad &
  }
\end{equation*}

We proved earlier (Section \ref{rectificationadjoint}) that the first
horizontal inclusion
\begin{equation*}
  \left\{ \begin{array}{lcl}
      \uAs \textsf{-alg} & \to & \infty \textsf{-}\uAsinf \textsf{-alg},\\
      \As \textsf{-alg} & \to & \infty \textsf{-} A_{\infty}\textsf{-alg}
    \end{array} \right.
\end{equation*}
has a left-adjoint, $\Omega_{\iota} B_{\kappa}$, which we called the
universal rectification (it is known that the second has a similarly
defined left-adjoint). Each of the (vertical) inclusions,
\begin{equation*}
  \left\{ \begin{array}{lcl}
      \uAs \textsf{-alg} & \to & \As \textsf{-alg},\\
      \textsf{su-}A_{\infty} \textsf{-alg} & \to & A_{\infty} \textsf{-alg}
    \end{array} \right.
\end{equation*}
has a left-adjoint as well, given by adjoining an element $u$ and
extending the product(s) to make $u$ a strict unit (with appropriate
annihilation conditions, in the case of $\suAinf$-algebras).

We now analyze the relationship between $\uAsinf$ and
unital-$A_{\infty}$ via our transfer formulae.

\begin{thm} \label{stricttransfer} Let $\xymatrix{ *{ V \ \ }
    \ar@<.8ex>[r]^(.45){i} & *{\ A \quad} \ar@(dr,ur)[]_h
    \ar@<.8ex>[l]^(.45){p}},$ be a strong deformation retract, and $\{
  \un^A, \smtwo_A \}$ a strict $\uAs$-structure on $A$. Suppose
  further that $h \left( \un^A \right) = 0.$ Then the operations
  $\muns(V)$ given by the transfer formulae (see definition in
  Proposition \ref{transferformulae}) have the property that
$$\muns(V) = 0$$
whenever $n \geq 2$ and $|S| \geq 1$. Furthermore, the
$\uAsinf$-morphism structure $J$ on the chain map $i$ has the property
that whenever $|S| \geq 1$,
$$J_{n}^{S} = 0,$$
even when $n=1$. In particular this means that the transferred
$\uAsinf$ structure is an $\suAinf$-algebra, and the
$\uAsinf$-$\infty$ quasi-isomorphism is an $\suAinf$-$\infty$
quasi-isomorphism.
\end{thm}

\begin{pf}
  For $n \geq 2, |S| \geq 1$, each summand in $\muns(V)$ contains as
  some part of the diagram (of compositions) the following composite:
$$\begin{array}{l}
  \vcenter{ \xymatrix@C=6pt@R=8pt{
      \\
      \ar@{{*}}[u]_A  \ar@{-}[d]^{h}\\
      *{}
    } }
  = h ( \un^A) = 0,
\end{array}$$
so each of those operations is itself 0. The same fact gives the result for $J$, along with the fact that
$$\begin{array}{l}
  J_{1}^{\{1\}} = \vcenter{ \xymatrix@C=6pt@R=8pt{
      \\
      \ar@{{*}}[u]_A  \ar@{-}[d]^{h}\\
      *{}
    }} = 0.
\end{array}$$

The vanishing of these higher operations and morphisms implies that
the transferred $\uAsinf$ structure and morphism are strictly-unital,
because the operad $\suAinf$ is the quotient of $\uAsinf$ by precisely
these operations.  $\cqfd$
\end{pf}

\begin{rem}
  We point out that since we are working over a field, and $d ( \un^A
  ) = 0$, it is always possible to choose a strong deformation retract
  between $V$ and $A$ so that $h \left( \un^A \right) = 0$ (provided,
  of course, $V$ is equivalent to $A$).
\end{rem}

\begin{cor}\label{qisostrict}
  Let $A$ be a $\uAsinf$-algebra. Then there exists a $\uAs$-algebra $R$, and an $\suAinf$-algebra structure on $H_{\bullet}(A)$ so that 
  $A \qiso R$ and $H_{\bullet}(A) \qiso R$. That is, for an arbitrary $\uAsinf$-algebra $A$, there is a minimal model for $A$ which is an $\suAinf$-algebra.
\end{cor}

\begin{pf}
  By Theorem \ref{universalrectification}, we have $I_{A}: A \rightsquigarrow^{\hspace{-.3cm} \sim} \hspace{.1cm} \Omega_{\kappa} \B_{\iota} A = R(A)$. Note that in particular,
  $H_{\bullet}\left( A \right) \simeq_{i} H_{\bullet} \left(
    R(A) \right) $. We will denote both by $H$.

Since there exist strong
  deformation retracts $\xymatrix{ *{ H \ \ } \ar@<.8ex>[r]^(.45){i} &
    *{\ R(A)  \quad} \ar@(dr,ur)[]_h
    \ar@<.8ex>[l]^(.45){p}}$ where $h$ annihilates the unit,
  transferring the $\uAs$ structure on $\Omega_{\kappa}
  \textrm{B}_{\iota} A $ along any such strong deformation
retract, by Theorem \ref{stricttransfer}, gives an equivalent
$\suAinf$-algebra structure on $H$.
$\cqfd$
\end{pf}

We point out that this method of proof (that one can always obtain
``strict units'') can be generalized to many cases in the following
way.

First, we axiomatize the relationship between the operads $\As$ and
$\uAs$.

\begin{dei}
  Let $\Po$ be a properad. We say an operad $u\Po$, together with an isomorphism:
  \begin{equation}
    { \q u\Po } \simeq \un \oplus \q \Po
  \end{equation}
is a \emph{unital version} of $\Po$.
\end{dei}

\begin{rem}
  For a unital version $u\Po$ of a Koszul properad $\Po$, we may define
  $\texttt{su}\Po_{\infty}$-algebras and morphisms, the same way we
  defined $\suAinf$-algebras using the unital version $\uAs$ of $\As$.
\end{rem}

\begin{thm}\label{generalstrictunittransfer}
  Let $\Po$ be a Koszul, inhomogeneous quadratic properad, and suppose
  $u\Po$ is a Koszul, inhomogeneous-quadratic properad which is a
  unital version of $\Po$. Then given any $u\Po$-algebra $A$ and a
  strong deformation retract $\xymatrix{ *{ V \ \ }
    \ar@<.8ex>[r]^(.45){i} & *{\ A \quad} \ar@(dr,ur)[]_h
    \ar@<.8ex>[l]^(.45){p}},$ where the homotopy $h$ satisfies
  $h(\un^A)=0$, the transferred $(u\Po)_{\infty}$-algebra is a
  strictly unital $\Po_{\infty}$-algebra structure, and the
  $(u\Po)_{\infty}$ $\infty$-morphism structure on $J$ is a strictly
  unital $\Po_{\infty}$ $\infty$-morphism.
\end{thm}

\begin{cor}
  Suppose we have properads $\Po, u\Po$ as in Theorem
  \ref{generalstrictunittransfer}, and suppose $A$ is a
  $u\Po$-algebra. Then there is an
  $\texttt{su}\Po_{\infty}$-algebra structure on the homology of $A$
  and an $\texttt{su}\Po_{\infty}$ $\infty$-quasi isomorphism $H \rightsquigarrow^{\hspace{-.3cm} \sim} \hspace{.1cm} A$.
\end{cor}

\begin{pf}
  It is a corollary of the proof for $\uAs$ and the transfer formulae
  for arbitrary Koszul inhomogeneous quadratic properads $u\Po$ (whose
  formulae are not made explicit in this paper).  $\cqfd$
\end{pf}

\begin{cor}
  The following $u\Po$-algebras can always be transferred to equivalent
  $\texttt{su}\Po_{\infty}$-algebra structures on homology in the above sense:
  \begin{enumerate}
  \item $u\Com$ --- the operad governing unital commutative and
    associative algebras --- to unital $\Com_{\infty}$
  \item $u\BV$ --- the operad governing unital $\BV$-algebras, where
    the unit is annihilated by the $\BV$-operator and the bracket ---
    to unital $\BV_{\infty}$ (see \cite{Ga-CaToVa} for a model of
    $\BV_{\infty}$).
  \end{enumerate}
\end{cor}

\begin{rem}
  Treating $\ucFrob$, the properad governing
  Frobenius algebras with unit and counit, would be interesting to the authors.
\end{rem}

\subsection{Cohomology theory for unital associative algebra}\label{cohomology}

In this section, we define the André-Quillen cohomology theory for unital associative algebras following the general definition of \cite{Milles}. We prove that the cohomology can be written as an Ext-functor and we compare this definition to the Hochschild cohomology theory.

\subsubsection{\bf André-Quillen cohomology theory}

We consider now the operad $\Po = \uAs$ and the curved cooperad $\C = \uAs^{\ash} = (\q\uAs^{\ash},\, 0,\, \theta)$. The Koszul morphism between $\uAs$ and $\uAs^{\ash}$ is given by
$$\kappa : \uAs^{\ash} \epi \vcenter{
      \xymatrix@M=0pt@R=4pt@C=3pt{&&\\
        & \ar@{{*}}[u] \ar@{-}[d]\\
        &}} \oplus  \vcenter{
\xymatrix@M=0pt@R=4pt@C=4pt{\ar@{-}[dr] && \ar@{-}[dl]\\
& \ar@{-}[d] &  \\ & &}} \mono \uAs.$$

Let $A$ be a $\uAs$-algebra. Following Sections $1$ and $2$ of \cite{Milles}, we use the cofibrant resolution
$$\Omega_{\kappa} \B_{\kappa} A = \uAs \circ_{\kappa} \uAs^{\ash} (A) \qiso A$$
of Section \ref{resofalg} to compute the André-Quillen cohomology of $A$ thanks to the cotangent complex
$${\small A \otimes^{\uAs}\uAs^{\ash}(A) \cong \vcenter{
\xymatrix@M=0pt@R=4pt@C=4pt{A & \uAs^{\ash}(A) & A\\
\ar@{-}[dr] & \ar@{-}[dd] & \ar@{-}[dl]\\
& &  \\ & &}}} \cong A \otimes \uAs^{\ash}(A) \otimes A.$$

We denote an element in $A \otimes^{\uAs}\uAs^{\ash}(A)$ by $a \otimes (\munsc \otimes b_{1} \cdots b_{n-|S|}) \otimes c$, where $a$, $b_{t}$ and $c$ are in $A$ and where $\munsc $ is in $\uAs^{\ash}(n-|S|)$. Following the end of Section $2$ of \cite{Milles}, we compute the differential on $A \otimes^{\uAs}\uAs^{\ash}(A)$, which is given by
$$d_{\varphi} := d_{A \otimes^{\uAs}\uAs^{\ash}(A)} - \delta_{\varphi}^{l} + \delta_{\varphi}^{r}.$$
The differential $d_{A \otimes^{\uAs}\uAs^{\ash}(A)}$ depends only on $d_{A}$ (since $d_{\uAs} = 0$, $d_{\uAs^{\ash}} = 0$), the map $\varphi : \uAs^{\ash} (A) \epi A$ is the projection and the terms $\delta_{\varphi}^{l}$ and $\delta_{\varphi}^{r}$ are given by the following proposition. 

%%%%%%%%%%%%%%%%%%%%%%%%%%%%%%%%%%%%%%%%

%%%%%%%%%%%%%%%%%%%%%%%%%%%%%%%%%%%%%%%%

\begin{prop}
To simplify the signs, we assume that $A$ is concentrated in degree $0$. We have\\
$\delta_{\varphi}^{l}(a \otimes (\munsc \otimes b_{1} \cdots b_{n-|S|}) \otimes c) :=$
$$\epsilon_{1} a\cdot b_{1} \otimes (\overline{\mu}_{n-1}^{S-1} \otimes b_{2} \cdots b_{n-|S|}) \otimes c + (-1)^{n} \epsilon_{2} a \otimes (\overline{\mu}_{n-1}^{S} \otimes b_{1} \cdots b_{n-|S|-1}) \otimes b_{n-|S|}\cdot c,$$
where $\epsilon_{1} := \left\{\begin{array}{ll}
1 & \textrm{ if } 1 \notin S,\\
0 & \textrm{ otherwise,}
\end{array} \right.$ and $\epsilon_{2} := \left\{\begin{array}{ll}
1 & \textrm{ if } n \notin S,\\
0 & \textrm{ otherwise,}
\end{array} \right.$
and\\
$\delta_{\varphi}^{r}(a \otimes (\overline{\mu}_{n}^{S} \otimes b_{1} \cdots b_{n-|S|}) \otimes c) := (\delta_{\small{\vcenter{
  \xymatrix@M=0pt@R=3pt@C=3pt{&&\\
    & \ar@{{*}}[u] \ar@{-}[d]\\
    &}}}} + \delta_{\gamma})(a \otimes (\overline{\mu}_{n}^{S} \otimes b_{1} \cdots b_{n-|S|}) \otimes c) =$
$$\hspace{-1.8cm} -\sum_{S = S_{1} \sqcup \{ u\} \sqcup S_{1}'}
(-1)^{n+|S_{1}|} a \otimes (\overline{\mu}_{n}^{S\backslash u} \otimes b_{1} \cdots 1_{A} \cdots b_{n-|S|}) \otimes c$$
$$\hspace{1.6cm} - \sum_{\{t,\, t+1\} \sqcup S = S_{2} \sqcup \{ t,\, t+1\} \sqcup S_{2}'} (-1)^{t+|S|} a \otimes (\overline{\mu}_{n-1}^{S_{2} \sqcup \{ S_{2}' - 1\}} \otimes b_{1} \cdots b_{t}\cdot b_{t+1} \cdots b_{n-|S|}) \otimes c,$$
where $\max S_{1} < u < \min S_{1}'$ and $\max S_{2} < t < t+1 < \min S_{2}'$ and $\delta_{\small{\vcenter{
  \xymatrix@M=0pt@R=3pt@C=3pt{&&\\
    & \ar@{{*}}[u] \ar@{-}[d]\\
    &}}}}$ holds for the first sum and $\delta_{\gamma}$ for the second. Moreover, $d_{\varphi} (\vcenter{
      \xymatrix@M=0pt@R=4pt@C=3pt{&&\\
        & \ar@{{*}}[u] \ar@{-}[d]\\
        &}}) = 0$.
\end{prop}

\begin{pf}
The differential on the cotangent complex is given following the end of Sections 2 of \cite{Milles}. We make the computations explicit thanks to the infinitesimal decomposition map of $\uAs^{\ash}$, described in Corollary \ref{quAs}.
$\cqfd$
\end{pf}

\begin{prop}
The André-Quillen cohomology groups of a $\uAs$-algebra $A$ with coefficients in a unital $A$-bimodule $M$ are given by
$$\mathrm{H}_{\Po}^{\bullet} (A, \, M) := (\mathrm{Hom}_{\textsf{A-bimod.}}(A \otimes^{\uAs}\uAs^{\ash}(A),\, M),\, \partial),$$
where $\partial (f) := d_{M} \circ f - (-1)^{|f|} f \circ d_{\varphi}$ and $\textsf{A-bimod.}$ is the category of unital $A$-bimodules.
\end{prop}

\subsubsection{\bf Ext-functor and comparison with the Hochschild cohomology theory}

To a unital associative algebra, we can associate two abelian groups: the Hochschild cohomology groups of $A$ (as defined in \cite{Hochschild}, or \cite{Loday}, chap. 1, for a modern reference), that is the André-Quillen cohomology groups of the associative algebra $A$ (forgetting the unit), or the André-Quillen cohomology groups of $A$ seen as a unital associative algebra (previous section). We show that the cohomology groups coincide.

\begin{thm}\label{Ext}
Let $A$ be a $\uAs$-algebra and let $M$ be a unital $A$-bimodule. We have
$$\mathrm{H}_{\uAs}^{\bullet} (A,\, M) \cong \mathrm{Ext}_{A \otimes^{\uAs} \K}^{\bullet} (\Omega_{\uAs}(A),\, M),$$
where $\Omega_{\uAs}(A)$ is the unital $A$-bimodule of K\"ahler differential forms (see \cite{Milles} for more details).
\end{thm}

\begin{pf}
Similarly to the case of Hochschild cohomology theory, we define the map $h$ on $A \otimes \uAs^{\ash}(A) \otimes A$ by
$$h(a \otimes (\munsc \otimes b_{1} \cdots b_{n-|S|}) \otimes c) = - 1 \otimes (\overline{\mu}_{n+1}^{S+1} \otimes a b_{1} \cdots b_{n-|S|}) \otimes c.$$
It satisfies $dh + hd = id$ on $A \otimes \overline{\uAs^{\ash}}(A) \otimes A$. Thus the chain complex
$$A \otimes \overline{\uAs^{\ash}}(A) \otimes A \xrightarrow{d_{\varphi}} A \otimes A \otimes A \epi \Omega_{\uAs}(A) \rightarrow 0$$
is acyclic since we derive the left-adjoint functor of K\"ahler differential forms to obtain the cotangent complex, and the cohomology is an Ext-functor.
$\cqfd$
\end{pf}

We use this theorem to compare this cohomology theory to the Hochschild cohomology theory.

\begin{prop}\label{projres}
There is a quasi-isomorphism of unital $A$-bimodules
$$A \otimes^{\uAs} \As^{\ash}(A) \cong A \otimes \As^{\ash}(A) \otimes A \qiso A \otimes^{\uAs} \uAs^{\ash}(A) \cong A \otimes \uAs^{\ash}(A) \otimes A.$$
\end{prop}

\begin{pf}
First, we endowed $A \otimes \uAs^{\ash}(A) \otimes A$ with a filtration given by the number of corks, denoted by
$$F_{p}(A \otimes \uAs^{\ash}(A) \otimes A) := \bigoplus_{S \subseteq \underline{n},\, |S| \leq p} A \otimes (\uAs^{\ash}(n-|S|) \otimes_{\So_{n-|S|}} A^{\otimes (n-|S|)}) \otimes A.$$
We have $\delta_{\varphi}^{l} : F_{p} \rightarrow F_{p}$, $\delta_{\small{\vcenter{
  \xymatrix@M=0pt@R=3pt@C=3pt{&&\\
    & \ar@{{*}}[u] \ar@{-}[d]\\
    &}}}} : F_{p} \rightarrow F_{p-1}$ and $\delta_{\gamma} : F_{p} \rightarrow F_{p}$. Thus the filtration is a filtration of chain complexes. It is bounded below and exhaustive so we can apply the classical theorem of convergence of spectral sequences (cf. Theorem 5.5.1 of \cite{Weibel}) and we obtain a spectral sequence $E_{p,q}^{\bullet}$ such that
$$E_{p,q}^{\bullet} \Rightarrow \mathrm{H}_{p+q}(A \otimes \uAs^{\ash}(A) \otimes A).$$
The differential $d^{0}$ on $E_{p,q}^{0} := F_{p}(A \otimes \uAs^{\ash}(A) \otimes A)_{p+q}/F_{p-1}(A \otimes \uAs^{\ash}(A) \otimes A)_{p+q}$ is given by $d^{0} = \delta_{\varphi}^{l} + \delta_{\small{\vcenter{
  \xymatrix@M=0pt@R=3pt@C=3pt{&&\\
    & \ar@{{*}}[u] \ar@{-}[d]\\
    &}}}} + \delta_{\gamma}$. There is an inclusion of chain complexes
$$i : A \otimes \As^{\ash}(A) \otimes A \mono \oplus_{p,\, q} E_{p,q}^{0} \cong A \otimes \uAs^{\ash}(A) \otimes A,$$
where the last isomorphism is only of vector spaces. The projection $p : \oplus_{p,\, q} E_{p,q}^{0} \cong A \otimes \As^{\ash}(A) \otimes A \oplus C_{\geq 1} \epi A \otimes \As^{\ash}(A) \otimes A$, where $C_{\geq 1}$ is given by elements with at least one cork, is a chain complexes map. We define the map $h$ by
$$h(a\otimes (\munsc \otimes b_{1} \cdots b_{n-|S|}) \otimes c) := -(-1)^{\min S} a \otimes (\overline{\mu}_{n+1}^{S+1} \otimes b_{1} \cdots b_{(\min S)-1} 1_{A} b_{\min S} \cdots b_{n-|S|}) \otimes c.$$
With these definitions, we have $p \circ i = id_{A \otimes \uAs^{\ash}(A) \otimes A}$ and $id_{\oplus_{p,\, q} E_{p,q}^{0}} - i \circ p = dh + hd$. Hence, we have a deformation retract
$$\xymatrix{     *{ \quad \ \  \quad A \otimes \As^{\ash}(A) \otimes A\ } \ar@<.8ex>[r]^(.55){i} & *{\
\oplus_{p,\, q} E_{p,q}^{0}\quad \ \  \ \quad } \ar@(dr,ur)[]_h \ar@<.8ex>[l]^(.43){p}}$$
and the inclusion $i$ is a quasi-isomorphism. It follows that $E_{p,q}^{1} = 0$ when $p \neq 0$ and the spectral sequence collapses. Considering the filtration $F'_{p}(A \otimes \As^{\ash}(A) \otimes A) = A \otimes \As^{\ash}(A) \otimes A$ for all $p \geq 0$ (bounded below and exhaustive), the inclusion induces a map of spectral sequences which is a quasi-isomorphism on the $E^{1}$-pages and higher. Since ${E'}_{p,q}^{\bullet}$ converges to $\mathrm{H}_{p+q}(A \otimes \As^{\ash}(A) \otimes A)$ and $E_{p,q}^{\bullet}$ converges to $ \mathrm{H}_{p+q}(A \otimes \uAs^{\ash}(A) \otimes A)$, we get the proposition.
$\cqfd$
\end{pf}

\begin{cor}\label{hochschild}
Let $A$ be a unital associative algebra. For $\bullet \geq 1$, we have
$$\mathrm{H}_{\uAs}^{\bullet} (A,\, M) \cong \mathrm{HH}^{\bullet +1} (A,\, M).$$
\end{cor}

\begin{pf}
The cohomology of $\uAs$-algebras is given by the Ext-functor $\mathrm{Ext}_{A \otimes^{\uAs} \K}^{\bullet} (\Omega_{\uAs}(A),\, M)$ (Theorem \ref{Ext}) and we have the projective resolution $A \otimes \uAs^{\ash}(A) \otimes A \qiso \Omega_{\uAs}(A)$. By Proposition \ref{projres}, the projective (quasi-free) $A$-bimodule $A \otimes \As^{\ash}(A) \otimes A$ is also a projective resolution of $\Omega_{\uAs}(A)$ and computes the Hochschild cohomology (see the definition 1.1.3 in \cite{Loday}).
$\cqfd$
\end{pf}

\section*{Acknowledgments}
The authors thank Bruno Vallette for his time, his help and all his
comments. The first author is grateful to Gabriel Drummond-Cole for his
tutelage in Koszul Duality while this project was in its early stages,
and for many helpful discussions throughout. He would also like to
thank John Terilla for his insights and encouragement. The second author is grateful to Éric Hoffbeck for his comments on the preliminary version of the paper. He is also grateful to the Max-Planck-Institut f\"ur Mathematik in Bonn for the excellent working conditions he found there between September 2009 and January 2010.

% To include all references
%\nocite{*}
\bibliographystyle{alpha}
\bibliography{bib}
\end{document}